\documentclass[review,hidelinks,onefignum,onetabnum]{siamart220329}
\usepackage{lipsum}
\usepackage{amsfonts}
\usepackage{graphicx,subfigure}
\usepackage{epstopdf}
\usepackage{algorithmic}
\usepackage{xr-hyper}
\usepackage{cleveref}
\usepackage{nicefrac}
\usepackage{bm}
\usepackage{bbm}
\usepackage{amssymb}
\usepackage{amsmath}
\usepackage{mathabx}
\usepackage{mathrsfs}
\usepackage{mathtools}
\usepackage{empheq}
\usepackage{color}
\usepackage[normalem]{ulem}
\usepackage{stmaryrd}
\usepackage{multirow}
\usepackage[colorinlistoftodos,prependcaption,textsize=tiny]{todonotes}
\usepackage{extarrows}
\usepackage[]{hyperref}

\ifpdf
  \DeclareGraphicsExtensions{.eps,.pdf,.png,.jpg}
\else
  \DeclareGraphicsExtensions{.eps}
\fi

\def\softd{{\leavevmode\setbox1=\hbox{d}%
          \hbox to 1.05\wd1{d\kern-0.4ex{\char039}\hss}}}

\newcommand{\dx}{\Delta x}
\newcommand{\dt}{\Delta t}
\newcommand{\mbu}{\mathbf{u}}

\newcommand{\mbf}{\mathbf{f}}
\newcommand{\mbS}{\mathbf{S}}
\newcommand{\mbE}{\mathbf{E}}
\newcommand{\mbG}{\mathbf{G}}
\newcommand{\mbR}{\mathbf{R}}
\newcommand{\hmbG}{\widehat{\mathbf{G}}}
\newcommand{\wmbG}{\widetilde{\mathbf{G}}}
\newcommand{\R}{\mathbb{R}}
\newcommand{\jph}{{j+\frac{1}{2}}}
\newcommand{\jmh}{{j-\frac{1}{2}}}
\newcommand{\mbPhil}{\overleftarrow{\mathbf{\Phi}}}
\newcommand{\mbPhir}{\overrightarrow{\mathbf{\Phi}}}
\newcommand{\hmbPhil}{\widehat{\overleftarrow{\mathbf{\Phi}}}}
\newcommand{\hmbPhir}{\widehat{\overrightarrow{\mathbf{\Phi}}}}
\newcommand{\wmbPhil}{\widetilde{\overleftarrow{\mathbf{\Phi}}}}
\newcommand{\wmbPhir}{\widetilde{\overrightarrow{\mathbf{\Phi}}}}

\newcommand{\avg}[1]{\left\{\!\!\left\{\smash{#1}\right\}\!\!\right\}}

\newcommand{\jump}[1]{\llbracket #1 \rrbracket}

\newcommand*\xbar[1]{%
  \hbox{%
    \vbox{%
      \hrule height 0.5pt % The actual bar
      \kern0.4ex%         % Distance between bar and symbol
      \hbox{%
        \kern-0.05em%      % Shortening on the left side
        \ensuremath{#1}%
        \kern-0.00em%      % Shortening on the right side
      }%
    }%
  }%
}

\hypersetup{
    colorlinks=true,
    linkcolor=green,     % equations
    citecolor=red,      % references
    urlcolor=blue    % URL
}

\newsiamthm{thm}{Theorem}
\newsiamremark{hypothesis}{Hypothesis}
\crefname{hypothesis}{Hypothesis}{Hypotheses}
\newsiamthm{claim}{Claim}
\newsiamthm{prop}{Proposition}

\headers{Positivity-Preserving Well-Balanced PAMPA Schemes}{R. Abgrall, Y. Liu, and M. Ricchiuto}

\graphicspath{{figures/}{./}}

\title{Positivity-preserving Well-balanced PAMPA Schemes %with Flux Globalization 
with Global Flux quadrature
for One-dimensional Shallow Water Models\thanks{Submitted to the editors DATE.
\funding{The work of Y. Liu was supported by UZH Postdoc Grant, 2024 / Verf\"{u}gung Nr. FK-24-110 and SNSF grant 200020$\_$204917.}}}

\author{R\'emi Abgrall\thanks{Institute of Mathematics, University of Z\"{u}rich, 8057 Z\"{u}rich, Switzerland} \and Yongle Liu\thanks{Corresponding author. Institute of Mathematics, University of Z\"{u}rich, 8057 Z\"{u}rich, Switzerland
  (\email{yongle.liu@math.uzh.ch; liuyl2017@mail.sustech.edu.cn}).} \and Mario Ricchiuto\thanks{Univ. Bordeaux, CNRS, Bordeaux INP, IMB, UMR 5251, F-33400 Talence, France.}}

\usepackage{amsopn}

\begin{document}

\maketitle

\nolinenumbers

% REQUIRED
\begin{abstract}
%We present a novel hydrostatic and non-hydrostatic equilibria preserving Point-Average-Moment PolynomiAl-interpreted (PAMPA) method for solving the one-dimensional hyperbolic balance laws, with applications to the shallow water models including the Saint--Venant system with the Manning friction term and rotating shallow water equations. The idea is based on a global flux formulation in which the primitive of the source term is added as an additional flux function. The reformulated system is quasi-conservative with global integral terms computed using Gauss--Lobatto quadrature nodes, recursively. The resulting method is capable of preserving a large family of smooth moving equilibria: supercritical and subcritical flows, in a super-convergent manner. We also show that, by an appropriate quadrature strategy for the source, we can exactly preserve the still water states for the Saint--Venant system and geostrophic equilibrium for the rotating shallow water equations. Moreover, to guarantee the positivity of water depth and eliminate the spurious oscillations near shocks, we blend the high-order PAMPA scheme with the first order local Lax--Fridrichs scheme using the method developed in \cite{Abgrall2025}. The first-order scheme is designed to preserve the still water equilibria and positivity of water height. Extensive numerical experiments are tested to validate the advantages and robustness of the proposed scheme. 
We present a novel hydrostatic and non-hydrostatic equilibria preserving Point-Average-Moment PolynomiAl-interpreted (PAMPA) method for solving the one-dimensional hyperbolic balance laws, with applications to the shallow water models including the Saint--Venant system with the Manning friction term and rotating shallow water equations. The idea is based on a global flux  quadrature formulation, in which
the discretization of the source terms is obtained from the derivative of
and additional flux function computed via high order quadrature  of the source term.
 The reformulated system is quasi-conservative with global integral terms computed using Gauss--Lobatto quadrature nodes. %, recursively. 
 The resulting method is capable of preserving a large family of smooth moving equilibria: supercritical and subcritical flows, in a super-convergent manner. We also show that, by an appropriate quadrature strategy for the source, we can exactly preserve the still water states. Moreover, to guarantee the positivity of water depth and eliminate the spurious oscillations near shocks, we blend the high-order PAMPA schemes with the first order local Lax--Friedrichs schemes using the method developed in \cite{Abgrall2025}. The first-order schemes are designed to preserve the still water equilibria and positivity of water height, as well as to deal with wet-dry fronts. Extensive numerical experiments are tested to validate the advantages and robustness of the proposed scheme. 

\end{abstract}

% REQUIRED

% REQUIRED

\section{Introduction}\label{sec1}
This paper focus on development of high-order and robust numerical methods for one-dimensional (1-D) shallow water models, which play a fundamental role in various phenomena in natural science and engineering. Such models are given by the following form:
\begin{equation}\label{1.1}
  \frac{\partial \mbu}{\partial t}+~\mbf(\mbu)_x=\mbS(\mbu,x),\quad x\in\Omega\subset\R,\quad t\geq0, 
\end{equation}
where $\mbu$ contains the conserved variables, $\mbf$ is the flux function, and $\mbS$ denotes the source term that depends on the solution as well as some external data. For example, the influence of spatial changes of the height of the bottom topography, of the Manning friction effect, as well as other physical effects such as Coriolis force, are accounted for by means of appropriately defined this source term. 

The addition of the source term results in a system of balance laws which may admit quite a large number of equilibrium solutions dictated by the interaction between different components of the flux variation in space and the forcing terms. Many of such equilibria have some interest in themselves and many physical applications involve small perturbations of such equilibria. The ability of a numerical method to resolve with enhanced accuracy such steady states is a unanimously acclaimed design criterion usually referred to as well-balanced (WB). Even in 1-D case, the challenge of devising WB numerical approximations of shallow water models agnostic of the form of the steady state is still a very active and open research domain. There is already quite a large literature on the subject, with several different approaches to manage this issue. It is quite impossible to mention all the relevant literatures, we mention a few numerical methods having been proposed in many different settings: finite difference (FD) methods \cite{Li2015,Li2020,Xu2025,Wang2022,Zhang2023,Ren2024}, finite volume (FV) methods \cite{Audusse2004,Castro2008,Cheng2016,Diaz2013,Noelle2007,Kurganov2020,Hwang2023,Guo2023,Ciallella2022,Busto2022,Ciallella2025,CaballeroCardenas2024,Chertock2022,Ciallella2023,Cao2022,Zhang2025}, Discontinuous Galerkin (DG) methods \cite{Liu2022,Zhang2023,Mantri2024,Zhang2025,Xu2024a,Yang2021,Haidar2022,Xing2014}, continuous finite element methods \cite{Azerad2017,Guermond2018,Behzadi2020, Knobloch2024,Micalizzi2024}, residual distribution methods \cite{Ricchiuto2009,Chang2023,Ricchiuto2015,Ricchiuto2011,Arpaia2018,Arpaia2020}, Active Flux methods \cite{Barsukow2024,Abgrall2024,Liu2025}, and so on. 

Concerning the preservation of steady states, if the steady-state solutions are explicitly known, one can use a simple and efficient idea that consists in evolving a discrete error with respect to the given equilibrium. This boils down to removing from the discrete equations of any scheme the discrete expression corresponding to the application of the scheme itself to the given solution. A more intricate scenario arises when the relevant equilibrium is not explicitly known, but it can be characterized by a set of (generally nonlinear) algebraic relations that define constant invariants. To put it differently, in this scenario, the equilibrium state can be described through the collection of relationships expressed as $\mbE(\mbu)=\mbE_0$, where $\mbE$ represents a comprehensive set of variables that, ideally, can be leveraged to derive the conserved quantities $\mbu$ and is referred to as equilibrium variable. In this case, one usually has to solve nonlinear equations at each cell and time to proceed the transformation from equilibrium to conserved variables and a nontrivial unisolvent root-finding mechanism should be carefully developed. This idea, known as generalized hydrostatic reconstruction, has been adapted to many discretization approaches; see, e.g., \cite{Xing2006,Noelle2007,Xing2014,Cheng2016,Cheng2019,Chertock2022,Xu2024a,Micalizzi2024}.

An alternative WB numerical method, based on the semi-discrete Active Flux approach introduced in \cite{Abgrall2023a}, was recently proposed in \cite{Abgrall2024} for solving the Saint-Vanent system of shallow water equations with/without friction term. This method is later referred to as Point-Average-Moment PolynomiAl-interpreted (PAMPA) method in \cite{Liu2025a}. The PAMPA method uses degrees of freedom (DoFs) from quadratic polynomials---all lying on element boundaries---along with the traditional average values within each element. The solution is a global continuous representation of these DoFs (point and average values), which are evolved simultaneously according to two formulations of the same PDEs. The point values can be evolved according to several writings of the same PDEs and can even be nonconservative. For instance, in \cite{Abgrall2024, Liu2025}, a primitive formulation with the flux derivative given by equilibrium variable $\mbE$ is used to evolve the point values and the WB property is automatically attained through consistent finite difference approximations of the flux derivative. Whereas, the average values should be updated following the conservative formulation \eqref{1.1}, due to the well-known Lax--Wendroff theorem. Given the cell boundary point values, the flux is easily computed and the WB evolution of average values only depends on an appropriate approximation of the source term. This not only requires the approximation should be WB but also satisfy exact conservation property or consistent approximation when the source term vanishes. To achieve this, in \cite{Abgrall2024} a cutoff function is used to switch off the inconsistent moving-water equilibrium correction term when the solution is away from steady states. A more elegant way that subtracts a local reference steady-state solution obtained from solving nonlinear equations derived from the relationship $\mbE(\mbu)=\mbE_0$, was used in \cite{Liu2025} to ensure the exact conservation property.

In this paper, without assuming either the steady states or local reference steady states a priori, we present a novel flux globalization based positivity-preserving (PP) WB PAMPA scheme for the 1-D shallow water models. This high-order scheme inherits the advantages of PAMPA methods in eliminating the need for Riemann solvers and free of any special WB (generalized) hydrostatic reconstructions in solving hyperbolic balance laws \cite{Abgrall2024,Liu2025}. Beyond the positivity-preserving, the low-order schemes are also capable of preserving the still water steady states, which differs from the scheme that we have used in the previous works \cite{Abgrall2024,Liu2025}. Following the idea from \cite{Chertock2018,Cheng2019,Ciallella2023,Mantri2024,Kazolea2025}, we incorporate the source term into the flux and rewrite \eqref{1.1} in the following equivalent quasi-conservative form:
\begin{equation}\label{1.2}
  \frac{\partial \mbu}{\partial t}+~\partial_x\mbG(\mbu,x)=0,
\end{equation}  
where the global operator (flux) is defined as
\begin{equation}\label{1.2a}
  \mbG(\mbu,x)=\mbf(\mbu)-\int_{\widehat{x}}^{x}\mbS(\mbu,\xi)\,\mathrm{d}\xi,
\end{equation} 
which is the natural operator to represent  solutions of the steady state ordinary differential equation (ODE)
$$
\partial_x \mbf = \mbS.
$$
Note that in \eqref{1.2a}  $\widehat{x}$ is a given initial integration point, for example the left-end of the domain. We can observe that now   any steady-state solution can be naturally expressed by
\begin{equation}\label{1.2b}
  \mbG\equiv{\rm Const},
\end{equation} 
thus \eqref{1.2} provides a reasonable path to obtain WB schemes capable of preserving a wider class of equilibria. The above condition alone however does not necessarily provide a precise characterization
of the discrete equilibria. Here, following \cite{Mantri2024,Kazolea2025},
we will pursue an approach that allows such characterization in the context of ODE integration methods.  To this end  a modified version of the the PAMPA method will be developed, using   the quasi-conservative form \eqref{1.2}  to update both the cell average values within each element, and point values at the cell boundaries.

We apply the proposed flux globalization based PP WB PAMPA scheme for solving two shallow water models. The first one is the Saint--Venant system with the Manning friction term:
\begin{equation}\label{1.3}
  \begin{aligned}
  &h_t+(hu)_x=0,\\
  &(hu)_t+(hu^2+\frac{g}{2}h^2)_x=-ghB_x-\frac{gn^2}{h^{\frac{7}{3}}}\vert hu\vert hu,
  \end{aligned}
\end{equation} 
where $h$ is the water depth, $u$ is the velocity, $B$ is bottom topography, $g$ is the constant gravitational acceleration, $n$ is the Manning friction coefficient. The global flux is given by
\begin{equation*}
  \mbG=\begin{pmatrix}
  &hu\\
  &hu^2+\frac{g}{2}h^2+\int_{\widehat{x}}^{x}\big[ghB_x+\frac{gn^2}{h^{\frac{7}{3}}}\vert hu \vert hu\big]\;\mathrm{d}\xi
  \end{pmatrix}.
\end{equation*} 

The second model is the shallow water equations with the Coriolis forces, which reads as
\begin{equation}\label{1.3a}
  \begin{aligned}
  &h_t+(hu)_x=0,\\
  &(hu)_t+(hu^2+\frac{g}{2}h^2)_x=-ghB_x+fhv,\\
  &(hv)_t+(huv)_x=-fhu,
  \end{aligned}
\end{equation}
where $u$ and $v$ are the two components of horizontal velocity, $f(x)=f_0+\beta x$ is the Coriolis parameter with positive constants $f_0$ and $\beta$. The global flux is given by
\begin{equation*}
  \mbG=\begin{pmatrix}
  &hu\\
  &hu^2+\frac{g}{2}h^2+\int_{\widehat{x}}^{x}\big[ghB_x-fhv\big]\;\mathrm{d}\xi\\
  &huv+\int_{\widehat{x}}^{x}fhu\;\mathrm{d}\xi
  \end{pmatrix}.
\end{equation*} 

For other 1-D shallow water models, such as the thermal rotating shallow water equations, the shallow water flows in an open channel, and the blood flow models in veins and arteries, we believe that the proposed numerical schemes are extendable, and we will pursue these extensions in future research.

Finally, we point out that the high-order PAMPA scheme is only linear stable under certain CFL (Courant--Friedrichs--Lewy) constraint \cite{Abgrall2023a}. When the solution develops discontinuities, the high-order PAMPA scheme will be prone to numerical oscillations. Additionally, the high-order PAMPA scheme cannot guarantee the preservation of water depth positivity, which is crucial due to not only the physical reason but also the numerical reason (code may crash when water depth is negative). We, therefore, require the developed scheme should also be positive-preserving. This can be achieved by considering a convex combination of high-order and low-order fluxes/residuals. We provide PP WB low order schemes based on hydrostatic reconstruction \cite{Audusse2004}. The optimal combination coefficients are then derived following the idea introduced by bound-preserving PAMPA method developed in \cite{Abgrall2025}. We also impose different conditions on the blending coefficients for the concern of spurious oscillations.
 
The main findings of this work can be summarized as follows:
\begin{itemize}
       \item We propose a high-order fully WB PAMPA scheme via the global flux Gauss--Lobatto quadrature. We prove, both theoretically and numerically, that the proposed scheme verifies a nodal super convergence property for moving water steady states. Compared to some existing fully WB flux globalization based finite-volume and discontinuous Galerkin methods, the proposed method avoids solving nonlinear algebraic equations (as required, for example, in \cite{Chertock2018,Chertock2022,Zhang2025}). We also show that the high-order WB PAMPA scheme is capable of exactly preserving the still water equilibria.

     \item We present first-order schemes that preserve simultaneously both still water equilibria and positivity of water height. Positivity is guaranteed by the use of the local Lax--Friedrichs' flavor fluxes and residuals, while still water equilibria are maintained through the hydrostatic reconstructions used in \cite{Audusse2004,Xing2011,Xing2010}. % or the topography correction parameter following \cite{Audusse2015,Hajduk2022}. 
         Numerical tests on perturbed still water steady states confirm these properties and further highlight the superior performance of the high-order schemes in capturing small perturbation dynamics. The proposed PP WB first-order schemes improve upon the purely PP first-order schemes used in our previous works \cite{Liu2025,Abgrall2024}.
         
     \item We introduce new convex blending methods of the global fluxes and residuals computed by the high-order PAMPA and low-order local Lax–Friedrichs' flavor schemes, following the idea in \cite{Abgrall2025}. By examining the intermediate solution states, we establish a sufficient condition for the PP property and derive the optimal blending coefficients. %For the first type of low order schemes, similar blending coefficients are obtained as in \cite{Abgrall2025}. For the second type of low order schemes, we derive the new form of blending coefficients based on the adopted hydrostatic reconstruction. 
         To reduce oscillations near strong discontinuities, when the solution is away from a local discrete steady state, we incorporate an oscillating-elimination coefficient based on the damping term given in \cite{Peng2025}, which is then combined with the PP blending coefficient. This approach differs from our previous WB PAMPA methods in \cite{Liu2025,Abgrall2024}, where a posteriori MOOD paradigm is used. %A CPU-time comparison on a test problem demonstrates that the proposed monolithic convex limiting strategy is more efficient.
           
\end{itemize}

The rest of the paper is organized as follows. In \cref{sec2}, we introduce the novel flux globalization based PP WB PAMPA scheme. We provide some theoretical analysis of the proposed schemes in \cref{sec3}, which includes the positivity-preserving and still water equilibria maintaining for the first-order scheme, as well as the positivity and still- and moving- water steady states preserving for the high-order scheme. In \cref{sec4}, extensive numerical examples are presented to illustrate the high-order accuracy in smooth regions for general solutions, the positivity-preserving and the well-balanced property over non-flat bottom topography, and the robustness in capturing shocks over flat and non-flat bottom topography, of the proposed schemes. 

\section{Numerical scheme}\label{sec2}
In this section, we describe the flux globalization based PP WB PAMPA scheme for the general quasi-conservative system \eqref{1.2}. To this end, we consider a tessellation of the 1-D spatial domain $\Omega$ in non overlapping elements $K_\jph=[x_j,x_{j+1}]$ with uniform size $\dx=x_{j+1}-x_j$. For reference, we list other commonly used notations and variables in Table \ref{tab_notation}.
\begin{table}[ht!]
\caption{\sf Notations and variables to be used in numerical methods.\label{tab_notation}}
\begin{center}
\begin{tabular}{ l l }
\hline
$t$  &The current time level \\
$\dt$  &The adaptive time step determined by the CFL condition\\
$\dx$  &The spatial grid size\\
$\lambda$  &The time-space ratio $\frac{\dt}{\dx}$\\ 
$\xbar\mbu_\jph$ &The cell averages of the numerical solution in $K_\jph$ at time $t$\\
$\mbu_j$   &The point values of the numerical solution at $x_j$ and time $t$\\
$\xbar\mbu_\jph^{t+\dt}$  &The cell averages of the numerical solution in $K_\jph$ at time $t+\dt$\\
$\mbu_j^{t+\dt}$   &The point values of the numerical solution at $x_j$ and time $t+\dt$\\
$\jump{\cdot}_{l,k}$ &The difference of a function/variable evaluated at two points: $(\cdot)_l-(\cdot)_k$\\
$\avg{\cdot}_{l,k}$ &The average of a function/variable evaluated at two points: $\frac{(\cdot)_l+(\cdot)_k}{2}$\\
\hline
\end{tabular}
\end{center}
\end{table}

We now assume that the solution given in terms of average values $\xbar{\mbu}_\jph$ and point values $\mbu_{j,j+1}$ for each element $K_\jph$ are available at the certain time level $t$. On each element $K_\jph$, we consider a global continuous finite element approximation $\mbu_{\mathrm{h}}(x)$\footnote{Following the finite element convention, we use the subscript $\mathrm{h}$ to indicate a finite element approximation of a variable. For example, $\mbu_{\mathrm{h}}$ denotes a finite element approximation of $\mbu$. Meanwhile, we have also used the letter $h$ for the water depth in this paper. Therefore, $h_{\mathrm{h}}$ would stand for a finite element approximation of the water depth $h$. Moreover, we have omit the dependence of $\mbu_{\mathrm{h}}$ and other quantities on $t$ for the sake of brevity from here on.} spanned by the quadratic polynomial basis functions:
\begin{equation}\label{2.1}
  \mbu_{\mathrm{h}}(x)=\varphi_j(x)\mbu_j+\varphi_{j+1}(x)\mbu_{j+1}+\xbar{\varphi}_\jph(x)\xbar{\mbu}_\jph,\quad x\in K_\jph,
\end{equation}
where $\mbu_{j}\approx\mbu(x_{j})$ is the point value of $\mbu$ at the element boundary $x_j$ and $\xbar{\mbu}_\jph$ is the cell average defined as
\begin{equation*}
  \xbar{\mbu}_\jph=\frac{1}{\vert K_\jph\vert}\int_{K_\jph}\mbu_{\mathrm{h}}(x)\;\mathrm{d}x.
\end{equation*}
Based on the following three conditions of the finite element approximation space:
\begin{itemize}
  \item conservation condition:
  \begin{equation*}
    \int_{K_\jph}\xbar{\varphi}_\jph(x)\,\mathrm{d}x=\vert K_\jph\vert,\quad \int_{K_\jph}\varphi_j(x)\,\mathrm{d}x=0,
  \end{equation*}
  \item interpolation condition:
  \begin{equation*}
    \xbar{\varphi}_\jph(x_j)=0,\quad \varphi_j(x_i)=\delta_{ij},
  \end{equation*}
  \item accuracy condition: $\varphi_j(x)$ and $\xbar{\varphi}_\jph(x)$ are quadratic polynomials,
\end{itemize}
we can easily obtain the explicit form of the basis functions (\cite{Abgrall2023a,Abgrall2024}):
\begin{equation*}
  \varphi_j=(1-\xi)(1-3\xi),\quad \xbar{\varphi}_\jph=6\xi(1-\xi),\quad \varphi_{j+1}=\xi(3\xi-2),
\end{equation*}
where $\xi=\frac{x-x_j}{\dx}\in[0,1]$.

The next question is to update the DoFs $\xbar{\mbu}_\jph$ and $\mbu_j$ in a PP and WB manner, which can be done by solving the following two semi-discrete forms:
\begin{subequations}\label{2.2}
  \begin{equation}\label{2.2a}
  \frac{\mathrm{d}}{\mathrm{d}t}\xbar{\mbu}_\jph=-\frac{\hmbG_{j+1}-\hmbG_j}{\dx}
\end{equation} 
and
\begin{equation}\label{2.2b}
  \frac{\mathrm{d}}{\mathrm{d}t}{\mbu}_j=-\big(\hmbPhir_\jmh+\hmbPhil_\jph\big).
\end{equation} 
\end{subequations}

\subsection{PP WB update for cell averages}
For each element $K_\jph$, we can define the PP WB blended fluxes $\hmbG_j$ and $\hmbG_{j+1}$ in \eqref{2.2a} as:
\begin{equation}\label{2.3}
\begin{aligned}
  \hmbG_j&={\wmbG}_j\big(\xbar\mbu_\jph,\xbar\mbu_\jmh\big)+\theta_j(\mbG_j-{\wmbG}_j):={\wmbG}_j+\theta_j\Delta\mbG_j,\\
   \hmbG_{j+1}&={\wmbG}_{j+1}\big(\xbar\mbu_{j+\frac{3}{2}},\xbar\mbu_\jph\big)+\theta_{j+1}(\mbG_{j+1}-{\wmbG}_{j+1}):={\wmbG}_{j+1}+\theta_{j+1}\Delta\mbG_{j+1},
\end{aligned}
\end{equation}
where $\mbG_j$ and $\mbG_{j+1}$ are the high-order fluxes based on the Gauss--Lobatto quadrature rule, $\theta_j, \theta_{j+1}\in[0,1]$ are the blending coefficients, ${\wmbG}_j={\wmbG}_j\big(\xbar\mbu_\jph,\xbar\mbu_\jmh\big)$ and ${\wmbG}_{j+1}={\wmbG}_{j+1}\big(\xbar\mbu_{j+\frac{3}{2}},\xbar\mbu_\jph\big)$ are the first-order local Lax--Friedrichs numerical fluxes, on the left and right interfaces, respectively. 

\subsubsection{Low order fluxes}\label{sec211}
We begin with constructing the first-order numerical fluxes, which is positivity-preserving, still water state preserving, and wet-dry fronts capturing.  We proceed as follows. First we write
\begin{equation}\label{eq.Gtil}
\begin{aligned}
  {\wmbG}_{j}&=\bm{\mathcal{F}}^{\mathrm{LLF}}\big(\mbu_{j}^+,\mbu_{j}^-\big) -\mbR_{\jph} +\int_{x_j}^{x_\jph} \mbS \;\mathrm dx, \\
  {\wmbG}_{j+1}&=\bm{\mathcal{F}}^{\mathrm{LLF}}\big(\mbu_{j+1}^+,\mbu_{j+1}^-\big)-\mbR_{\jph} -\int_{x_\jph}^{x_{j+1}} \mbS \;\mathrm dx,
  \end{aligned}
\end{equation}
with $\bm{\mathcal{F}}^{\mathrm{LLF}}$ denoting the local Lax-Friedrich's flux. 
Then, to compute the left and right states we use the hydrostatic reconstruction idea  \cite{Audusse2004}. In particular, we define:
\begin{equation*}
\begin{aligned}
 & h_j^{+}=\xbar h_\jph+\xbar B_\jph-B_j^{+},&& h_j^-=\xbar h_\jmh+\xbar B_\jmh-B_j^-,\\
 & u_j^+=\frac{\xbar{(hu)}_\jph}{\xbar h_\jph},&& u_j^-=\frac{\xbar{(hu)}_\jmh}{\xbar h_\jmh},\\
 & v_j^+=\frac{\xbar{(hv)}_\jph}{\xbar h_\jph},&& v_j^-=\frac{\xbar{(hv)}_\jmh}{\xbar h_\jmh},
  \end{aligned}
\end{equation*} 
with modified reconstructions of the topography at  cell interfaces  given by
\begin{equation*}
\begin{aligned}
  &B_j^+=\min(\xbar h_\jph+\xbar B_\jph, \max(\xbar B_\jph, \xbar B_\jmh)),\\
  &B_j^-=\min(\xbar h_\jmh+\xbar B_\jmh, \max(\xbar B_\jph, \xbar B_\jmh)).
\end{aligned}
\end{equation*}

For the Saint-Venant system \eqref{1.3}, we  end with  the left and right low order fluxes as:
\begin{equation}\label{flux2Al}
\begin{aligned}
  {\wmbG}_j&=\bm{\mathcal{F}}^{\mathrm{LLF}}\big(\mbu_j^+,\mbu_j^-\big)
  -g\begin{pmatrix}
                0 \\
                \avg{h}_{\jph,j^+}\jump{B}_{\jph,j^+}
              \end{pmatrix}\\
  &\qquad -\frac{gn^2}{2}\dx\begin{pmatrix}
                0 \\
                \avg{{\vert hu\vert hu}/{h^{\frac{7}{3}}}}_{\jph,j^+} 
              \end{pmatrix} -\mbR_\jph
\end{aligned}
\end{equation}
and
\begin{equation}\label{flux2Ar}
\begin{aligned}
  {\wmbG}_{j+1}&=\bm{\mathcal{F}}^{\mathrm{LLF}}\big(\mbu_{j+1}^+,\mbu_{j+1}^-\big)
  +g\begin{pmatrix}
                0 \\
                \avg{h}_{{j+1}^-,\jph}\jump{B}_{{j+1}^-,\jph}
              \end{pmatrix}\\
   &\qquad +\frac{gn^2}{2}\dx\begin{pmatrix}
                0 \\
                \avg{{\vert hu\vert hu}/{h^{\frac{7}{3}}}}_{{j+1}^-,\jph} 
              \end{pmatrix} -\mbR_\jph.
\end{aligned}
\end{equation}

For the rotating shallow water equations \eqref{1.3a}, we can define the left and right low order fluxes as:
\begin{equation}\label{flux2Al2}
  {\wmbG}_j=\bm{\mathcal{F}}^{\mathrm{LLF}}\big(\mbu_j^+,\mbu_j^-\big)
  -\begin{pmatrix}
                0 \\
                g\avg{h}_{\jph,j^+}\jump{B}_{\jph,j^+}-\frac{\dx}{2}\avg{fhv}_{\jph,j^+}\\
                \frac{\dx}{2}\avg{fhu}_{\jph,j^+}
              \end{pmatrix}-\mbR_\jph,
\end{equation}
and
\begin{equation}\label{flux2Ar2}
\begin{split}
  {\wmbG}_{j+1}=&\bm{\mathcal{F}}^{\mathrm{LLF}}\big(\mbu_{j+1}^+,\mbu_{j+1}^-\big)
  \\+&\begin{pmatrix}
                0 \\
                g\avg{h}_{{j+1}^-,\jph}\jump{B}_{{j+1}^-,\jph}-\frac{\dx}{2}\avg{fhv}_{{j+1}^-,\jph}\\
                \frac{\dx}{2}\avg{fhu}_{{j+1}^-,\jph}
              \end{pmatrix}-\mbR_\jph.
              \end{split}
\end{equation}

%\begin{remark}
%In the computation of \eqref{fluxAl}, \eqref{fluxAr}, \eqref{alpha_jB}, \eqref{flux2Al}, and \eqref{flux2Ar}, the values used for the low-order schemes at cell center $x=x_\jph$ are the cell averages. We have also used
%\begin{equation*}
%  \xbar u_\jph=\frac{\xbar{(hu)}_\jph}{\xbar h_\jph},\quad \forall j.
%\end{equation*}
%\end{remark}

\subsubsection{High order fluxes} Thanks to the global continuous approximation in the PAMPA framework, we do not need to use any Riemann solver and directly compute the high order fluxes at the cell interfaces instead:
\begin{equation}\label{2.4}
\begin{aligned}
  &\mbG_j=\mbf(\mbu_j)-\mbR_j,&& 
  \mbG_{j+1}=\mbf(\mbu_{j+1})-\mbR_{j+1}\\
  &\mbR_j=\int_{\widehat{x}}^{x_j}\mbS(\mbu_{\mathrm h}, x)\,\mathrm{d}x && \mbR_{j+1}=\int_{\widehat{x}}^{x_{j+1}}\mbS(\mbu_{\mathrm h}, x)\;\mathrm{d}x.
  \end{aligned}
\end{equation}
In order to complete the computations in \eqref{eq.Gtil} and \eqref{2.4}, we need to compute the source primitive terms $\{\mbR_j\}$ and $\{\mbR_\jph\}$. This can be done in several ways. In \cite{Mantri2024} the authors have used  a nodal Lagrange expansion of the source term on the three Gauss--Lobatto points $[x_j,\;x_{j+1/2},\; x_{j+1}]$. This leads to
\begin{equation}\label{2.5}
\begin{array}{ll}
\mbR_{j+1}=&\!\!\!\mbR_j+\int_{x_j}^{x_{j+1}}\mbS_h(\mbu_{\mathrm h},x)\;\mathrm{d}x
  =\mbR_j+\sum_{i=1}^3\omega_i\mbS\big(\mbu_{\mathrm h}(x_i),x_i\big)\\[10pt]
  \mbR_\jph=&\!\!\!\mbR_j+\int_{x_j}^{x_\jph}\mbS_h(\mbu_{\mathrm h},x)\;\mathrm{d}x
  =\mbR_j+\sum_{i=1}^3\omega'_i\mbS\big(\mbu_{\mathrm h}(x_i),x_i\big)
  \end{array}\!\!\!\!,\;
    \forall\,j 
\end{equation}
with  $\omega_i$ and $\omega'_i$ obtained integrating the Lagrange polynomials on the interval and half-interval respectively.
This formulation can be written in compact form as 
\begin{equation}\label{eq.collocation1}
\left(\begin{array}{c}
\mbR_{j}\\
\mbR_{\jph}\\
\mbR_{j+1}
\end{array}\right) =
\mbR_{j} + \mathrm{I}_{\text{LIIIA}} \left(\begin{array}{c}
\mbS_{j}\\
\mbS_{\jph}\\
\mbS_{j+1}
\end{array}\right)
\end{equation}
where, as  noted in \cite{Mantri2024},  \emph{by construction} the matrix  
\begin{equation}\label{eq.collocation12}
 \mathrm{I}_{\text{LIIIA}} = \left[\begin{array}{ccc}
    0  & 0& 0 \\
   \omega'_1 & \omega'_2 & \omega'_{3}  \\
      \omega_1 & \omega_2 & \omega_{3}  
 \end{array}
 \right
 ] =\dx \left[\begin{array}{ccc}
    0  & 0& 0 \\
   5/24 & 1/3 & -1/24   \\
      1/6 & 2/3  & 1/6
 \end{array}
 \right
 ] 
\end{equation}
is the Butcher table of the well--known  fully implicit Runge--Kutta collocation method LobattoIIIA \cite{Hairer1993}. For autonomous problems, this method can be  shown to have a convergence order $p+2$ for the midpoint, and $2p$ for the end-points (here $p=2$ is the polynomial degree in the finite element approximation space). These estimates have been confirmed on the steady ODE in \cite{Mantri2024},  and lead in our case to  a fourth order method.  Despite of this, note that the second quadrature  rule
associated to the second line of \eqref{eq.collocation1} (or equivalently 
 the second in  \eqref{2.5}) is not symmetric and only exact for quadratic polynomials.

A different approach is to apply the usual three-points Gauss--Lobatto formula within  both the half, and full intervals. In this case, both formulas are exact for cubic polynomials, and we obtain
a sub-cell LobattoIII method which can be  written as  
\begin{equation}\label{eq:scLIII1}
\left(\begin{array}{c}
\mbR_{j}\\
\mbR_{\jph}\\
\mbR_{j+1}
\end{array}\right) =
\mbR_{j} + \mathrm{I}_{\text{sc-LIII}} \left(\begin{array}{c}
\mbS(\mbu_{j},x_j)\\[2pt]
\mbS(\mbu_{j+1/4},x_{j+1/4})\\[2pt]
\mbS(\mbu_{\jph},x_{\jph})\\[2pt]
\mbS(\mbu_{j+1}, x_{j+1})
\end{array}\right)\;,
\end{equation}
and
\begin{equation}\label{eq:scLIII2}
\mathrm{I}_{\text{sc-LIII}} 
  =\dx \left[\begin{array}{cccc}
    0  & 0& 0 & 0\\
   1/12 & 1/3 & 1/12 & 0   \\
      1/6 & 0 & 2/3  & 1/6
 \end{array}
 \right
 ] \,,\quad
 \mbu_{j+1/4} = \dfrac{3}{16} \mbu_j +  \dfrac{9}{8} \xbar\mbu_{\jph}-
  \dfrac{5}{16} \mbu_{j+1},
 %\mbu_{j+1/4} = \dfrac{3}{8} \mbu_j +  \dfrac{3}{4} \mbu_{\jph}-\dfrac{1}{8} \mbu_{j+1}
\end{equation}
where the last expression is obtained from quadratic interpolation within the cell \eqref{2.1}. We can show the following result.
\begin{proposition}[sc-LIII consistency\label{pp:sc-lobattoIII}] The implicit sub-cell LobattoIII RK me\-thod defined by  \eqref{eq:scLIII1}-\eqref{eq:scLIII2} is fourth order accurate for smooth enough solutions.
\begin{proof} 
The result can be easily shown by truncated Taylor series developments. Details are reported in \cref{app:sc-lobattoIII} for completeness.
\end{proof}
\end{proposition}

In both cases,  the computation of the nodal values of $\mbR$ can be seen as the integration of the steady--state ODE with  the  LobattoIIIA or sc-LobattoIII methods. A very important remark is that in practice \emph{only differences of values of $\mbR$ appear in our scheme}. This, and the continuity of the approximation, make the method  independent on the initial quadrature point. As a consequence, despite the name commonly used in literature, this approach leads to a fully local method. 
If necessary to prepare the solution, one can initiate the computation of these terms from $\widehat{x}=x_0$ (left-end of the domain), and set $\mbR_0=0$. However, as already observed, this initial value plays no role in practice.   

The two methods have similar properties. For the shallow water models, the sc-LobattoIII approach has some advantage as it allows to integrate exactly the hydrostatic term. For this reason we will use this approach in the following.

\subsection{PP WB update for point values}
Analogously to the end--point fluxes, the PP blended residuals are defined as
\begin{equation}\label{2.6}
\begin{aligned}
 \hmbPhir_\jmh&=\wmbPhir_\jmh+\theta_\jmh(\mbPhir_\jmh-\wmbPhir_\jmh)
 :=\wmbPhir_\jmh+\theta_\jmh\Delta\mbPhir_\jmh,\\
  \hmbPhil_\jph&=\wmbPhil_{\jph}+\theta_{\jph}(\mbPhil_{\jph}-\wmbPhil_{\jph})
  :=\wmbPhil_{\jph}+\theta_{\jph}\Delta\mbPhil_{\jph}.
  \end{aligned}
\end{equation} 

\subsubsection{Low order residuals}
Mimicking  the low order fluxes based on the hydrostatic reconstruction, we construct low order point fluctuations. First, we define the following quantities
\begin{equation*}
\begin{aligned}
  &B_{j+\frac{1}{4}}^+=\min(\xbar h_\jph+\xbar B_\jph, \max(\xbar B_\jph, B_j)),~ B_{j+\frac{1}{4}}^-=\min(h_j+B_j, \max(\xbar B_\jph, B_j)),\\
  &h_{j+\frac{1}{4}}^+=\xbar h_\jph+\xbar B_\jph-B_{j+\frac{1}{4}}^+,~ h_{j+\frac{1}{4}}^-=h_j+B_j-B_{j+\frac{1}{4}}^-,\\
  &u_{j+\frac{1}{4}}^+=\frac{\xbar{(hu)}_\jph}{\xbar h_\jph},~ u_{j+\frac{1}{4}}^-=u_j,~v_{j+\frac{1}{4}}^+=\frac{\xbar{(hv)}_\jph}{\xbar h_\jph},~ v_{j+\frac{1}{4}}^-=v_j.
\end{aligned}
\end{equation*}
Then, the low order residuals are given by
\begin{equation}\label{low_res_b1}
\begin{aligned}
  \frac{\dx}{2}\wmbPhir_\jmh&=\mbf(\mbu_j)-\bm{\mathcal{F}}^{\mbox{LLF}}(\mbu_{{j-\frac{1}{4}}^+},\mbu_{{j-\frac{1}{4}}^-})+\begin{pmatrix}
  0 \\
  g\avg{h}_{j,{j-\frac{1}{4}}^+}\jump{B}_{j,{j-\frac{1}{4}}^+} 
  \end{pmatrix} \\
  &~+\frac{gn^2}{4}\dx\begin{pmatrix}
                   0 \\
                   \avg{\vert hu\vert hu/h^{\frac{7}{3}}}_{j,{j-\frac{1}{4}}^+} 
                 \end{pmatrix}
  \end{aligned}
\end{equation}
and
\begin{equation}\label{low_res_b2}
 \begin{aligned}
  \frac{\dx}{2}\wmbPhil_\jph&=\bm{\mathcal{F}}^{\mbox{LLF}}(\mbu_{j+\frac{1}{4}}^+,\mbu_{j+\frac{1}{4}}^-)-\mbf(\mbu_j)
  +\begin{pmatrix}
  0 \\
  g\avg{h}_{{j+\frac{1}{4}}^-,j}\jump{B}_{{j+\frac{1}{4}}^-,j} 
  \end{pmatrix} \\
  &~+\frac{gn^2}{4}\dx\begin{pmatrix}
                   0 \\
                   \avg{\vert hu\vert hu/h^{\frac{7}{3}}}_{{j+\frac{1}{4}}^-,j} 
                 \end{pmatrix}
  \end{aligned}
\end{equation}

For the rotating shallow water equations \eqref{1.3a}, the low order residuals are given by
\begin{equation}\label{low_res2_b1}
\begin{aligned}
  \frac{\dx}{2}\wmbPhir_\jmh&=\mbf(\mbu_j)-\bm{\mathcal{F}}^{\mbox{LLF}}(\mbu_{{j-\frac{1}{4}}^+},\mbu_{{j-\frac{1}{4}}^-})+\begin{pmatrix}
  0 \\
  g\avg{h}_{j,{j-\frac{1}{4}}^+}\jump{B}_{j,{j-\frac{1}{4}}^+}\\
  0
  \end{pmatrix} \\
  &~-\frac{\dx}{4}\begin{pmatrix}
                   0 \\
                   \avg{fhv}_{j,{j-\frac{1}{4}}^+} \\
                   -\avg{fhu}_{j,{j-\frac{1}{4}}^+}
                 \end{pmatrix}
  \end{aligned}
\end{equation}
and
\begin{equation}\label{low_res2_b2}
 \begin{aligned}
  \frac{\dx}{2}\wmbPhil_\jph&=\bm{\mathcal{F}}^{\mbox{LLF}}(\mbu_{j+\frac{1}{4}}^+,\mbu_{j+\frac{1}{4}}^-)-\mbf(\mbu_j)
  +\begin{pmatrix}
  0 \\
  g\avg{h}_{{j+\frac{1}{4}}^-,j}\jump{B}_{{j+\frac{1}{4}}^-,j} \\
  0
  \end{pmatrix} \\
  &~-\frac{\dx}{4}\begin{pmatrix}
                   0 \\
                   \avg{fhv}_{{j+\frac{1}{4}}^-,j} \\
                   -\avg{fhu}_{{j+\frac{1}{4}}^-,j}
                 \end{pmatrix}
  \end{aligned}
\end{equation}
\subsubsection{High order residuals}
The high-order residuals in \eqref{2.6} are given by
\begin{equation}\label{2.6b}
  \mbPhir_\jmh=\widetilde{J}^+(\mbu_j)\partial_x^+\mbG,\quad \mbPhil_{\jph}=\widetilde{J}^-(\mbu_j)\partial_x^-\mbG,
\end{equation}
where $\partial_x^+\mbG$ and $\partial_x^-\mbG$ are the left- and right-biased high-order FD approximations of $\mbG_x$ at $x=x_j$, respectively. They are given by
\begin{equation}\label{2.6c}
  \partial_x^+\mbG=\frac{1}{\dx}\big(\mbG_{j-1}-4\mbG_\jmh+3\mbG_j\big),\quad
  \partial_x^-\mbG=\frac{1}{\dx}\big(-3\mbG_j+4\mbG_{\jph}-\mbG_{j+1}\big),
\end{equation}
where $\mbG_j$ and $\mbG_{j+1}$ for each element $K_\jph$ are computed in \eqref{2.4} and 
\begin{equation*}
  \mbG_\jph=\mbf(\mbu_\jph)-\mbR_\jph
\end{equation*}
with $\mbu_\jph$ being computed by the local scaling PP limiter:
\begin{equation}\label{h_C}
  \mbu_\jph=(1-\eta)\xbar{\mbu}_\jph+\eta\mathring{\mbu}_\jph,\quad \eta\in[0,1],
\end{equation}
where $\mathring{\mbu}_\jph=\frac{3}{2}\xbar{\mbu}_\jph-\frac{1}{4}(\mbu_j+\mbu_{j+1})$ obtained from \eqref{2.1} and $\eta$ is given by
\begin{equation*}
 \eta=\min\Big(1,\frac{\xbar{h}_\jph-\varepsilon_h}{\xbar{h}_\jph-\mathring{h}_\jph}\Big),
 \quad \varepsilon_h:=\min_j\big(10^{-13},\xbar{h}_\jph\big).
\end{equation*} 
Note that $\mathring{\mbu}_\jph$ does not guarantee the positivity of water depth and thus the PP limiter is crucial. Similar strategy has also been applied in the computation of $\mbu_{j+\frac{1}{4}}$ in \eqref{eq:scLIII2}.

The sign matrices are defined as
\begin{equation*}
  \widetilde{J}^\pm(\mbu_j)=J^{-1}J^\pm=\mathcal{R}\Lambda^\pm\Lambda^{-1}\mathcal{R}^{-1}
  =\mathcal{R}\begin{pmatrix}
                \frac{\lambda_1^{\pm}}{\lambda_1} &  &  \\
                 & \frac{\lambda_2^{\pm}}{\lambda_2}  &  \\
                 &  & \ddots 
              \end{pmatrix}\mathcal{R}^{-1},
\end{equation*}
where $\lambda_i$, $i=1,2,\ldots$ are the eigenvalues of the Jacobian $J=\frac{\partial f}{\partial \mbu}$ and $\mathcal{R}$ is the corresponding eigenvectors matrix. $\lambda_i^+=\max(\lambda_i,0)$ and $\lambda_i^-=\min(\lambda_i,0)$.
 
The idea behind defining the high-order residuals in \eqref{2.6b} are based on the following two conditions:
\begin{itemize}
  \item linear stable, which is achieved by considering the upwinding flavor and biased FD approximations of the spatial derivatives,  
  \item consistency: if we have $(\partial_x\mbG)_{x_j}={\rm Const}$, we are able to recover
  \begin{equation*}
    \mbPhir_\jmh+\mbPhil_{\jph}=(\partial_x\mbG)_{x_j}.
  \end{equation*}
\end{itemize}
With these definitions we can easily show the following property.
\begin{proposition}[Compactness\label{compactness}] The positivity preserving PAMPA method with global flux quadrature  is fully local, and has a compact stencil.
\begin{proof}
To show the property, it is enough to check that everywhere we only use differences
of the values of $\mbR$. We can start with the numerical fluxes for the  update 
of the cell average. Combining \eqref{eq.Gtil} and \eqref{flux2Al}--\eqref{flux2Ar} (or \eqref{flux2Al2}--\eqref{flux2Ar2}) we can write
$$
\wmbG_{j+1} - \wmbG_{j} =  \bm{\mathcal{F}}^{\mathrm{LLF}}\big(\mbu_{j+1}^+,\mbu_{j+1}^-\big) -\int_{x_{\jph}}^{x_{j+1}}\mbS\; \mathrm dx- \bm{\mathcal{F}}^{\mathrm{LLF}}\big(\mbu_{j}^+,\mbu_{j}^-\big) 
- \int_{x_{j}}^{x_{\jph}}\mbS \; \mathrm dx
$$
where the integrals are replaced by the approximations in
\eqref{flux2Al}--\eqref{flux2Ar} (or \eqref{flux2Al2}--\eqref{flux2Ar2}).
This expression is independent of $\mbR$. Similarly, we can see that in \eqref{2.3}
\begin{equation*}
    \begin{aligned}
&\Delta \mbG_j = \mbf_j - \mbR_j +\mbR_\jph -
 \bm{\mathcal{F}}^{\mathrm{LLF}}\big(\mbu_{j}^+,\mbu_{j}^-\big)
- \int_{x_{j}}^{x_{\jph}}\mbS\;\mathrm dx,\\
&\Delta \mbG_{j+1} = \mbf_{j+1} - \mbR_{j+1} +\mbR_\jph -
 \bm{\mathcal{F}}^{\mathrm{LLF}}\big(\mbu_{j+1}^+,\mbu_{j+1}^-\big)
+\int_{x_{\jph}}^{x_{j+1}}\mbS\;\mathrm dx,
    \end{aligned}
\end{equation*}
where the integrals are again replaced by the approximations in \eqref{flux2Al}--\eqref{flux2Ar} (or \eqref{flux2Al2}--\eqref{flux2Ar2}) 
and where by definition $\mbR_\jph - \mbR_j$ is obtained from 
the second in \eqref{2.5}, or from \eqref{eq:scLIII1}--\eqref{eq:scLIII2}, which  both only involve computations within the cell (similar to $\mbR_{j+1}-\mbR_\jph=\mbR_{j+1}-\mbR_j+\mbR_j-\mbR_\jph$). This shows that the limited update for the cell average is local. 
 
Concerning the nodal updates, we just need to show that the residuals  are local. 
We can recast \eqref{2.6c} as 
$$
\partial_x^+\mbG =
\dfrac{4}{\Delta x} ( \mbG_j - \mbG_\jmh) -
\dfrac{1}{\Delta x} ( \mbG_j - \mbG_{j-1}).
$$
And similarly for $\partial_x^-\mbG$ (omitted for brevity). For the source part this only involves the differences 
$$
\begin{aligned}
\dfrac{4}{\Delta x} ( -\mbR_j + \mbR_\jmh)
-&\dfrac{1}{\Delta x} ( -\mbR_j + \mbR_{j-1}) =
-\dfrac{3}{\Delta x} ( \mbR_j- \mbR_{j-1} ) + \dfrac{4}{\Delta x}(  \mbR_\jmh-\mbR_{j-1} ) 
\end{aligned}
$$
which again can be obtained from local evaluations given by \eqref{eq.collocation1}--\eqref{eq.collocation12}, or  \eqref{eq:scLIII1}--\eqref{eq:scLIII2}.
So  the  high order fluctuations are also fully local. 
\end{proof}\end{proposition}

\subsection{Convex combination coefficients}\label{sec23}
Here, we derive the convex combination coefficients that used for the blended fluxes and residuals in \eqref{2.3} and \eqref{2.6} respectively, to close this section. As we mentioned in \cref{sec1}, the proposed PP WB PAMPA scheme will be applied for solving the Saint--Venant system \eqref{1.3} and the rotating shallow water equations \eqref{1.3a}. Consequently, we only need to ensure the PP property for the water depth $h$. 

We apply the low order schemes introduced in \cref{sec211} and consider the forward Euler time stepping for the average value $\xbar h_\jph$ from $t$ to $t+\dt$:
\begin{equation*}
  \begin{aligned}
  \xbar h_\jph^{t+\dt}&=\xbar h_\jph-\lambda\Big(\hmbG_{j+1}^{(1)}-\hmbG_j^{(1)}\Big)\pm\frac{\lambda}{2}(\alpha_{j+1}+\xbar u_\jph)h_{j+1}^-\pm\frac{\lambda}{2}(\alpha_j-\xbar u_\jph)h_j^+\\
  &=\xbar h_\jph-\frac{\lambda}{2}(\alpha_{j+1}+\xbar u_\jph)h_{j+1}^--\frac{\lambda}{2}(\alpha_j-\xbar u_\jph)h_j^+\\
  &\quad +\lambda\alpha_j\Big(\frac{1}{2}\big(1-\frac{\xbar u_\jph}{\alpha_j}\big)h_j^++\frac{\hmbG_j^{(1)}}{\alpha_j}\Big)
  +\lambda\alpha_{j+1}\Big(\frac{1}{2}\big(1+\frac{\xbar u_\jph}{\alpha_{j+1}}\big)h_{j+1}^--\frac{\hmbG_{j+1}^{(1)}}{\alpha_{j+1}}\Big)\\
  &=\widetilde{h}_\jph^C+\lambda\alpha_j\widetilde{h}_\jph^L+\lambda\alpha_{j+1}\widetilde{h}_\jph^R,
  \end{aligned}
\end{equation*}
where
\begin{equation}\label{hC_v2}
\begin{aligned}
  \widetilde{h}_\jph^C&=\xbar h_\jph-\frac{\lambda}{2}(\alpha_{j+1}+\xbar u_\jph)h_{j+1}^--\frac{\lambda}{2}(\alpha_j-\xbar u_\jph)h_j^+\\
  &\geq\xbar h_\jph\Big(1-\lambda\alpha_{j+1}\frac{h_{j+1}^-}{\xbar h_\jph}-\lambda\alpha_j\frac{h_j^+}{\xbar h_\jph}\Big)\geq0.
\end{aligned}
\end{equation}
The last inequality due to the fact that $0\leq h_{j+1}^-\leq\xbar h_\jph$, $0\leq h_j^+\leq\xbar h_\jph$, and the CFL condition $\lambda\alpha\leq\frac{1}{2}$ with $\alpha=\max_j\{\alpha_j\}$. Moreover,
\begin{equation*}
\begin{aligned}
  \widetilde{h}_\jph^L&=\frac{1}{2}\big(1-\frac{\xbar u_\jph}{\alpha_j}\big)h_j^++\frac{\hmbG_j^{(1)}}{\alpha_j}
=\frac{1}{2}h_j^-(1+\frac{\xbar u_\jmh}{\alpha_j})+\theta_j\frac{\Delta\mbG^{(1)}_{j}}{\alpha_{j}},\\
\widetilde{h}_\jph^R&=\frac{1}{2}\big(1+\frac{\xbar u_\jph}{\alpha_{j+1}}\big)h_{j+1}^--\frac{\hmbG_{j+1}^{(1)}}{\alpha_{j+1}}
=\frac{1}{2}h_{j+1}^+(1-\frac{\xbar u_{j+\frac{3}{2}}}{\alpha_{j+1}})-\theta_{j+1}\frac{\Delta\mbG^{(1)}_{j+1}}{\alpha_{j+1}}.
\end{aligned}
\end{equation*}
Since $\frac{1}{2}h_j^-(1+\frac{\xbar u_\jmh}{\alpha_j})\geq0$ and $\frac{1}{2}h_{j+1}^+(1-\frac{\xbar u_{j+\frac{3}{2}}}{\alpha_{j+1}})\geq0$, by requiring $\widetilde{h}_\jph^L\geq0$ and $\widetilde{h}_\jph^R\geq0$, we therefore obtain the value of blending coefficient, i.e.,
\begin{equation}\label{theta_j2}
  \theta_j=\min\Big(1, \frac{\alpha_j\frac{h_j^+}{2}(1-\frac{\xbar u_\jph}{\alpha_j})}{\vert\Delta\mbG_j^{(1)}\vert}, \frac{\alpha_j\frac{h_j^-}{2}(1+\frac{\xbar u_\jmh}{\alpha_j})}{\vert\Delta\mbG_j^{(1)}\vert}\Big).
\end{equation}  

Similarly, we can  write down the forward Euler time stepping for the point values as
\begin{equation*}
\begin{aligned}
  h_j^{t+\dt}&=h_j-2\lambda(\hmbPhir_\jmh^{(1)}+\hmbPhil_\jph^{(1)})\pm\lambda(\alpha_{j+\frac{1}{4}}+u_j)h_{j+\frac{1}{4}}^-\pm\lambda(\alpha_{j-\frac{1}{4}}-u_j)h_{j-\frac{1}{4}}^+\\
  &=h_j\Big(1-\lambda\big(\alpha_{j+\frac{1}{4}}+u_j\big)\frac{h_{j+\frac{1}{4}}^-}{h_j}-\lambda\big(\alpha_{j-\frac{1}{4}}-u_j\big)\frac{h_{j-\frac{1}{4}}^+}{h_j}\Big)\\
  &\quad +\lambda(\alpha_{j-\frac{1}{4}}+u_j)h_{j-\frac{1}{4}}^+-2\lambda\hmbPhir_\jmh^{(1)}+\lambda(\alpha_{j+\frac{1}{4}}+u_j)h_{j+\frac{1}{4}}^--2\lambda\hmbPhil_\jph^{(1)}\\
  &=\widetilde h_j^C+2\lambda\alpha_{j-\frac{1}{4}}\widetilde h_j^L+2\lambda\alpha_{j+\frac{1}{4}}\widetilde h_j^R.
  \end{aligned}
\end{equation*}
Since $0\leq h_{j-\frac{1}{4}}^+\leq h_j$, $0\leq h_{j+\frac{1}{4}}^-\leq h_j$, and the CFL condition $\lambda\alpha\leq\frac{1}{4}$, we immediately obtain $\widetilde{h}_j^C\geq0$. Then, we observe 
\begin{equation*}
\begin{aligned}
  \widetilde h_j^L&=\frac{1}{2}\big(1+\frac{u_j}{\alpha_{j-\frac{1}{4}}}\big)h_{j-\frac{1}{4}}^+-\frac{\hmbPhir_\jmh^{(1)}}{\alpha_{j-\frac{1}{4}}}
  =\frac{h_{j-\frac{1}{4}}^+}{2}\big(1+\frac{\xbar u_\jmh}{\alpha_{j-\frac{1}{4}}}\big)-\theta_\jmh\frac{\Delta \mbPhir_\jmh^{(1)}}{\alpha_{j-\frac{1}{4}}},\\
  \widetilde h_j^R&=\frac{h_{j+\frac{1}{4}}^-}{2}\big(1-\frac{u_j}{\alpha_{j+\frac{1}{4}}}\big)-\frac{\hmbPhil_\jph^{(1)}}{\alpha_{j+\frac{1}{4}}}
  =\frac{h_{j+\frac{1}{4}}^-}{2}\big(1-\frac{\xbar u_\jph}{\alpha_{j+\frac{1}{4}}}\big)-\theta_\jph\frac{\Delta \mbPhil_\jph^{(1)}}{\alpha_{j+\frac{1}{4}}}.
\end{aligned}
\end{equation*} 
Knowing that $\frac{1}{2}\big(1+\frac{\xbar u_\jmh}{\alpha_{j-\frac{1}{4}}}\big)h_{j-\frac{1}{4}}^+\geq0$ and $\frac{1}{2}\big(1-\frac{\xbar u_\jph}{\alpha_{j+\frac{1}{4}}}\big)h_{j+\frac{1}{4}}^-\geq0$, we can ask $\widetilde h_j^L\geq0$ and $\widetilde h_j^R\geq0$. Thus, we obtain
\begin{equation}\label{theta_jh2}
  \theta_\jmh=\min\Big(1,\alpha_{j-\frac{1}{4}}\frac{\frac{h_{j-\frac{1}{4}}^+}{2}\big(1+\frac{\xbar u_\jmh}{\alpha_{j-\frac{1}{4}}}\big)}{\vert \Delta \mbPhir_\jmh^{(1)} \vert}\Big),~
  \theta_\jph=\min\Big(1, \alpha_{j+\frac{1}{4}}\frac{\frac{h_{j+\frac{1}{4}}^-}{2}\big(1-\frac{\xbar u_\jph}{\alpha_{j+\frac{1}{4}}}\big)}{\vert \Delta \mbPhil_\jph^{(1)} \vert }\Big)
\end{equation}
for any $K_\jph$.

%\begin{remark}\label{compactness}
%We note that the low order flux given in \eqref{flux2Al} and \eqref{flux2Ar} are local values. When we blend the high-order and low-order fluxes, we should ensure that these fluxes are compatible. Therefore, we recast the global fluxes into the local ones as
%\begin{equation*}
%  \mbG_j=\mbf(\mbu_j)+(\mbR_\jph-\mbR_j),\quad \mbG_{j+1}=\mbf(\mbu_{j+1})-(\mbR_{j+1}-\mbR_\jph)
%\end{equation*} 
%in the blended numerical fluxes. Indeed, considering the semi-discretization of \eqref{1.2}--\eqref{1.2a}, we have
%\begin{equation*}
%\begin{aligned}
%  \frac{\mathrm d}{\mathrm dt}\xbar\mbu_\jph&=-\frac{1}{\dx}\big(\mbG_{j+1}-\mbG_j\big)\\
%  &=-\frac{1}{\dx}\Big(\mbf(\mbu_{j+1})-(\mbR_{j+1}-\mbR_\jph)-\big(\mbf(\mbu_j)+(\mbR_\jph-\mbR_j)\big)\Big)\\
%  &=-\frac{1}{\dx}\Big(\mbf(\mbu_{j+1})-\mbR_{j+1}-\big(\mbf(\mbu_j)-\mbR_j\big)\Big).
%  \end{aligned}
%\end{equation*} 
%\end{remark}

\begin{remark}
There are several versions of hydrostatic reconstruction, see, e.g., \cite{Chen2017}. We think that many of them can be applied in the proposed monolithic convex limiting framework and the combination coefficients derived in \eqref{theta_j2} and \eqref{theta_jh2} will keep the same form.
\end{remark}

\subsubsection{Suppress spurious oscillations}
In order to suppress spurious oscillations  as much as possible when  approximating   discontinuous solutions, we revisit the combination weights determined by the high-order derivative jumps at interfaces. The blending coefficient is designed as (inspired by \cite{Abgrall2025b})
\begin{equation}\label{theta_OE}
    \theta^{\text{OE}}_\jph=
     \exp\Big(-\frac{\alpha_\jph\dt_n\sigma_\jph(h_{\mathrm h})}{\dx}\Big),
\end{equation}
where $\dt_n$ is the $n$-th adaptive time step computed using a suitable CFL condition, $\alpha_\jph$ is the estimate of the local maximum wave speed around $K_\jph$. The key of the method is  $\sigma_\jph$, which here is computed based on the solution of water depth, and plays a critical role in balancing accuracy and stability. Ideally, $\sigma_\jph$ remains very small in smooth regions to preserve high-order accuracy, and becomes sufficiently large near discontinuities to suppress spurious oscillations. Getting inspiration from \cite{Peng2025}, we define it as
\begin{equation*}%\label{jump0}
  \sigma_\jph(h_{\mathrm h})=\sum_{m=0}^{2}\sigma_\jph^m(h_{\mathrm h})
\end{equation*}
with
\begin{equation*}
  \sigma_\jph^m=\left\{
  \begin{aligned}
  &0,&&\mbox{if}~ h_{\mathrm h}\equiv \langle h_{\mathrm h}\rangle_{\Omega} ,\\%\text{avg}_\Omega(h_{\mathrm h}),\\
  &\dx^m%\frac{\dx^m}{3m!}
  %\frac{\vert \llbracket \partial_x^mh_{\mathrm h}\rrbracket_j \vert+\vert \llbracket \partial_x^mh_{\mathrm h}\rrbracket_{j+1} \vert}{2\Vert h_{\mathrm h}-\text{avg}_\Omega(h_{\mathrm h})\Vert_{L^\infty(\Omega)}},
  \frac{\vert \llbracket \partial_x^mh_{\mathrm h}\rrbracket_j \vert+\vert \llbracket \partial_x^mh_{\mathrm h}\rrbracket_{j+1} \vert}{2\Vert h_{\mathrm h}-\langle h_{\mathrm h}\rangle_{\Omega}\Vert_{L^\infty(\Omega)}},&& \mbox{otherwise},
  \end{aligned}\right.
\end{equation*}
% \begin{equation*}
%   \sigma_\jph^m=\left\{
%   \begin{aligned}
%   &0,&&\mbox{if}~ h_{\mathrm h}\equiv \text{avg}_\Omega(h_{\mathrm h}),\\
%   &\frac{(2m+1)\dx^m}{3m!}\frac{\vert \llbracket \partial_x^mh_{\mathrm h}\rrbracket_j \vert+\vert \llbracket \partial_x^mh_{\mathrm h}\rrbracket_{j+1} \vert}{2\Vert h_{\mathrm h}-\text{avg}_\Omega(h_{\mathrm h})\Vert_{L^\infty(\Omega)}}, && \mbox{otherwise},
%   \end{aligned}\right.
% \end{equation*}
where $\llbracket ~\cdot~ \rrbracket_j$ depicts the jump quantity at $x=x_j$ and $\langle h_{\mathrm h}\rangle_{\Omega}=\frac{1}{\vert\Omega\vert}\int_\Omega h_{\mathrm h}{\mathrm d}x$ denotes the global average of $h_{\mathrm h}(x)$ over the entire computational domain $\Omega$.

\begin{figure}[ht!]
\centerline{\includegraphics[trim=0.2cm 0.05cm 0.5cm 0.05cm,clip,width=5.5cm]{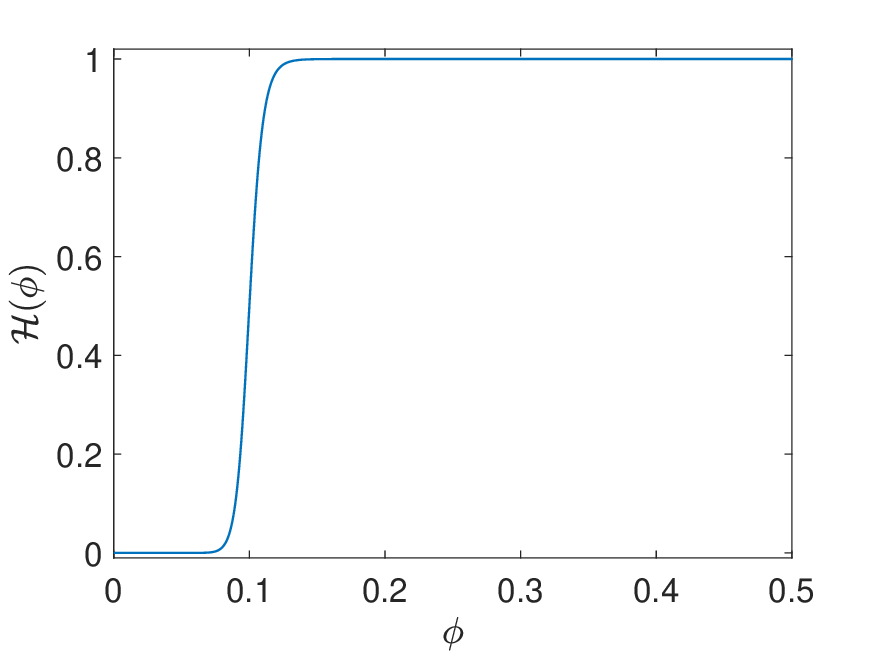}}
\caption{\sf Sketch of the cut-off function $\mathcal{H}(\phi)$.\label{H_cutoff}}
\end{figure}

As currently formulated, this correction would 
alter the fully well-balanced property. For this reason,
in the numerical implementation  we only activate this oscillation eliminating procedure when the numerical solution is far from equilibrium. 
In particular, we say the solution is near/at a local steady state around cell $K_\jph$ if $\mathcal H(\phi_\jph)\leq10^{-3}$, where
\begin{equation}\label{Hfunc}
  \mathcal H(\phi_\jph)=\frac{(C\phi_\jph)^\kappa}{1+(C\phi_\jph)^\kappa}
\end{equation}
with the constants $C=10$ and $\kappa=20$ being used in all of the numerical experiments, and
\begin{equation}\label{phi_s}
 \phi_\jph=%\text{avg}_{K_\jph}\!\big(\partial\mbG^{(2)}_x\big)\cdot\frac{\vert\Omega\vert}{\max\{\mbG^{(2)}_j,\mbG^{(2)}_\jph,\mbG^{(2)}_{j+1}\}}.
 \langle  \partial\mbG^{(2)}_x\rangle_{K_\jph}\cdot\frac{\vert\Omega\vert}{\max\{\mbG^{(2)}_j,\mbG^{(2)}_\jph,\mbG^{(2)}_{j+1}\}}.
\end{equation}
We plot a sketch of this function in Figure \ref{H_cutoff}, where one can clearly see that if $\phi$ is small (near/at a steady state), then the value of $\mathcal{H}$ is still close to 0. When $\phi$ increases, the values of $\phi$ rapidly approaches 1.

In the end, the combination coefficients for the fluxes and residuals are taken as the minimal value between \eqref{theta_OE} and \eqref{theta_j2}, \eqref{theta_jh2}, respectively.

\begin{remark}
In the near dry regions (water height $h$ closes to zero), the potentially harmful division by zero, for example when we need to compute the velocity $u=\frac{hu}{h}$, should be controlled. We achieve this following the idea in \cite{Dumbser2024}: if $h\leq10^{-14}$, we set $u=0$. Otherwise, if $h>h_0=10^{-4}$, the velocity is computed normally as $u=\frac{hu}{h}$. If instead $h<h_0$, then the velocity is computed by a filtered division as
\begin{equation*}
  u=hu\frac{h}{h^2+f(h)\varepsilon},
\end{equation*}
where $\varepsilon=5\times 10^{-9}$ and $f(h)$ is a filter function defined as
\begin{equation*}
  f(h)=2(h/h_0)^3-3(h/h_0)^2+1.
\end{equation*}
For more details, we refer to \cite[Section 3.4]{Dumbser2024}.
\end{remark}

\section{Positivity-preserving and well-balanced properties}\label{sec3}
In this section, we provide some basic properties of the proposed flux globalization based WB PAMPA schemes.

\begin{prop}[Positivity-preserving of the first-order schemes]\label{propa1}
%In the first-order numerical schemes with the numerical fluxes given by \eqref{fluxAl}--\eqref{fluxAr} or \eqref{flux2Al}--\eqref{flux2Ar} and the residuals given by \eqref{2.6a}--\eqref{2.6a2} or \eqref{low_res_b1}--\eqref{low_res_b2}, if we set
In the first-order numerical schemes with the numerical fluxes given by \eqref{flux2Al}--\eqref{flux2Ar} or \eqref{flux2Al2}--\eqref{flux2Ar2} and the residuals given by \eqref{low_res_b1}--\eqref{low_res_b2} or \eqref{low_res2_b1}--\eqref{low_res2_b2}, if we set
\begin{equation}\label{CFL}
  \lambda\alpha\leq\frac{1}{4}, \quad \lambda=\frac{\dx}{\dt},\quad \alpha=\max_{\forall K_\jph}\big(\alpha_j,\alpha_\jph\big),
\end{equation}
the schemes preserve the following property: if $\xbar h_\jph\geq0$ and $h_j\geq0$, then $\xbar h_\jph^{t+\dt}\geq0$, $h_j^{t+\dt}\geq0$, for any $j$.
\begin{proof}
If the low order fluxes (\eqref{flux2Al}--\eqref{flux2Ar} or \eqref{flux2Al2}--\eqref{flux2Ar2}) and residuals (\eqref{low_res_b1}--\eqref{low_res_b2} or \eqref{low_res2_b1}--\eqref{low_res2_b2}) are taken into account. We have
\begin{equation*}
\begin{aligned}
  \xbar h_\jph^{t+\dt}&=\xbar h_\jph-\lambda\Big(\wmbG_{j+1}^{(1)}-\wmbG_j^{(1)}\Big)\\
  &=\xbar h_\jph-\frac{\lambda}{2}(\alpha_{j+1}-\xbar u_\jph)h_{j+1}^--\frac{\lambda}{2}(\alpha_j-\xbar u_\jph)h_j^+\\
  &\quad +\frac{\lambda\alpha_j}{2}h_j^-(1+\frac{\xbar u_\jph}{\alpha_j})+\frac{\lambda\alpha_{j+1}}{2}h_{j+1}^+(1-\frac{\xbar u_{j+\frac{3}{2}}}{\alpha_{j+1}}).
\end{aligned}
\end{equation*}
From \eqref{hC_v2}, the non-negativity $h_j^-\geq0$ and $h_{j+1}^+\geq0$, and the fact of $1+\frac{\xbar u_\jph}{\alpha_j}\geq0$, $1-\frac{\xbar u_{j+\frac{3}{2}}}{\alpha_{j+1}}\geq0$, we have $\xbar h_\jph^{t+\dt}\geq0$. 

Similarly, for the point values, we have
\begin{equation*}
\begin{aligned}
h_j^{t+\dt}&=h_j-2\lambda(\wmbPhir_\jmh^{(1)}+\wmbPhil_\jph^{(1)})\\
&\xlongequal{\eqref{low_res_b1},\eqref{low_res_b2}}\Big(1-\lambda\alpha_{j+\frac{1}{4}}\bigg(1+\frac{u_j}{\alpha_{j+\frac{1}{4}}}\bigg)\frac{h_{j+\frac{1}{4}}^-}{h_j}
-\lambda\alpha_{j-\frac{1}{4}}\bigg(1-\frac{u_j}{\alpha_{j-\frac{1}{4}}}\bigg)\frac{h_{j-\frac{1}{4}}^+}{h_j}\Big)h_j\\
&\quad +\lambda\alpha_{j-\frac{1}{4}}h_{j-\frac{1}{4}}^-\Big(1+\frac{\xbar{u}_\jmh}{\alpha_{j-\frac{1}{4}}}\Big)
+\lambda\alpha_{j+\frac{1}{4}}h_{j+\frac{1}{4}}^+\Big(1-\frac{\xbar u_\jph}{\alpha_{j+\frac{1}{4}}}\Big)\geq0,
\end{aligned}  
\end{equation*}
since $h_{j\pm\frac{1}{4}}^{\pm}\geq0$, $-\alpha_{j-\frac{1}{4}}\leq\pm u_j,\pm\xbar u_\jmh\leq\alpha_{j-\frac{1}{4}}$, and $-\alpha_{j+\frac{1}{4}}\leq\pm u_j,\pm\xbar u_\jph\leq\alpha_{j+\frac{1}{4}}$.
\end{proof}

\end{prop}

%\begin{remark}\label{rmk:beta}
%If the discrete solutions satisfy a still water steady state with $u=0$, $w=h+B=w_0 (\mbox{constant})$ and with arbitrary topography, for the parameters defined in \eqref{alpha_jB}, \eqref{alpha_JB1}, and \eqref{alpha_JB2}, we have
%\begin{equation*}
%  \beta_j=1,\quad \beta_{\jmh}=1,\quad \beta_\jph=1,\quad \forall K_\jph.
%\end{equation*}
%\end{remark}
%\begin{proof}
%At the still water equilibrium, from \eqref{alpha_jB}, \eqref{alpha_JB1}, and \eqref{alpha_JB2}, we obtain
%\begin{equation*}
%  \begin{aligned}
%  \beta_j&=\min\Big(1,\frac{\xbar h_\jmh+\xbar h_\jph}{\vert \xbar B_\jph-\xbar B_\jmh\vert}\Big)=\min\Big(1,\frac{\xbar h_\jmh+\xbar h_\jph}{\vert \xbar h_\jmh-\xbar h_\jph\vert}\Big)=1,\\
%   \beta_\jmh&=\min\Big(1,\frac{\xbar h_\jmh+h_j}{ \vert B_j-\xbar B_\jmh\vert }\Big)=\min\Big(1,\frac{\xbar h_\jmh+h_j}{\vert \xbar h_\jmh-h_j\vert}\Big)=1,\\
%    \beta_\jph&=\min\Big(1,\frac{h_j+\xbar h_\jph}{ \vert B_j-\xbar B_\jph\vert}\Big)=\min\Big(1,\frac{h_j+\xbar h_\jph}{ \vert \xbar h_\jph-h_j\vert}\Big)=1.
%  \end{aligned}
%\end{equation*} 
%This completes the proof.
%\end{proof}

\begin{prop}[Preserving still water of the first-order schemes]\label{prop:still_low}
If the discrete solutions satisfy the still water steady state with $u=v=0$, $w=h+B=w_0 (\mbox{constant})$ and with arbitrary topography, the first-order schemes can exactly preserve this equilibrium for any time.
\begin{proof}
  At the steady state, we have $h_j+B_j=w_0$, $\xbar h_\jph+\xbar B_\jph=w_0$, $(hu)_j=0$, $\xbar{(hu)}_\jph=0$, $(hv)_j=0$, $\xbar{(hv)}_\jph=0$, for any $j$. It is easy to obtain $\frac{\mathrm d}{\mathrm dt}\xbar\mbu_\jph^{(3)}=0$ and $\frac{\mathrm d}{\mathrm dt}\mbu_j^{(3)}=0$. We now focus on the first- and second-components. 

For the average values, we apply \eqref{flux2Al}--\eqref{flux2Ar} and have
\begin{equation*}
  \frac{\mathrm d}{\mathrm d t}\xbar\mbu_\jph^{(1)}=\frac{\alpha_{j+1}}{2\dx}\big(h_{j+1}^+-h_{j+1}^-\big)
  -\frac{\alpha_{j}}{2\dx}\big(h_j^+-h_j^-\big).
\end{equation*}
Since $h_{j}^+=\xbar w_\jph-\max(\xbar B_\jph,\xbar B_\jmh)$, $h_j^-=\xbar w_\jmh-\max(\xbar B_\jph,\xbar B_\jmh)$, and $\xbar w_\jph=\xbar w_\jmh$ for any j, we obtain $h_j^+=h_j^-$ for any $j$ and thus $\frac{\mathrm d}{\mathrm d t}\xbar\mbu_\jph^{(1)}=0$. For the second component, we have
 \begin{equation*}
\begin{aligned}
  \frac{\mathrm d}{\mathrm dt}\xbar\mbu_\jph^{(2)}&=-\frac{g}{2\dx}\Big(\frac{1}{2}(h_{j+1}^+)^2+\frac{1}{2}( h_{j+1}^-)^2+(h_{j+1}^-+\xbar h_\jph)(B_{j+1}^--\xbar B_\jph)\\
  &\quad -\frac{1}{2}(h_j^+)^2-\frac{1}{2}(h_j^-)^2+(h_j^++\xbar h_\jph)(\xbar B_\jph-B_j^+)\Big)\\
  &=-\frac{g}{4\dx}\Big((h_{j+1}^+)^2-(h_{j+1}^-)^2+(h_j^+)^2-(h_j^-)^2\Big)=0,
\end{aligned}
\end{equation*}
as $h_j^+=h_j^-$ at the steady state, for all $j$. For the point values, we apply \eqref{low_res_b1}--\eqref{low_res_b2} and also note that $h_{j\pm\frac{1}{4}}^+=h_{j\pm\frac{1}{4}}^-$, $h_{j-\frac{1}{4}}^+=h_{j+\frac{1}{4}}^-$, and obtain
\begin{equation*}
  \begin{aligned}
  \frac{\mathrm d}{\mathrm dt}\mbu_j^{(1)}&=-\frac{1}{\dx}\Big(\alpha_{j-\frac{1}{4}}(h_{j-\frac{1}{4}}^+-h_{j-\frac{1}{4}}^-)
  +g(h_j+h_{j-\frac{1}{4}}^+)(B_j-B_{j-\frac{1}{4}}^+)\\
  &\quad -\alpha_{j+\frac{1}{4}}(h_{j+\frac{1}{4}}^+-h_{j+\frac{1}{4}}^-)
  +g(h_j+h_{j+\frac{1}{4}}^-)(B_{j+\frac{1}{4}}^--B_j)\Big)\\
  &=-\frac{1}{\dx}\Big(g(h_j+h_{j-\frac{1}{4}}^+)(h_{j-\frac{1}{4}}^+-h_j)
  +g(h_j+h_{j+\frac{1}{4}}^+)(h_j-h_{j+\frac{1}{4}}^-)\Big)\\
  &=-\frac{g}{\dx}((h_{j-\frac{1}{4}}^+)^2-(h_{j+\frac{1}{4}}^-)^2)=0.
  \end{aligned}
\end{equation*}
For the second component, we have
\begin{equation*}
  \begin{aligned}
  \frac{\mathrm d}{\mathrm dt}\mbu_j^{(2)}&=-\frac{g}{\dx/2}\Big(\frac{1}{4}(h_{j-\frac{1}{4}}^+)^2+
  \frac{1}{4}(h_{j-\frac{1}{4}}^-)^2+\frac{1}{2}(h_j+h_{j-\frac{1}{4}}^+)(B_j-B_{j-\frac{1}{4}}^+)\\
  &\quad \frac{1}{4}(h_{j+\frac{1}{4}}^+)^2+\frac{1}{4}(h_{j+\frac{1}{4}}^-)^2+\frac{1}{2}(h_j+h_{j+\frac{1}{4}}^-)(B_{j+\frac{1}{4}}^--B_j)\Big)\\
  &=-\frac{g}{\dx}\Big(-(h_{j-\frac{1}{4}}^+)^2+(h_{j+\frac{1}{4}}^-)^2+(h_j+h_{j-\frac{1}{4}}^+)(h_{j-\frac{1}{4}}^+-h_j)\\
  &\quad +(h_j+h_{j+\frac{1}{4}}^-)(h_j-h_{j+\frac{1}{4}}^-)\Big)=0.
  \end{aligned}
\end{equation*}
\end{proof}

\end{prop}

\begin{remark}\label{rmk1}
If the discrete solutions satisfy the steady state characterized by $\mbG={\rm Const}$, namely, 
\begin{equation*}
 \mbG_j=\mbG_{j+1}=\mbG_\jph,\quad \forall K_\jph,
\end{equation*}
we then have the blending coefficients $\theta_j=1$ and $\theta_\jph=1$. This implies that, at a steady state, the scheme reduces to the high-order WB PAMPA method. 
\begin{proof}
From \eqref{theta_OE}, \eqref{Hfunc}, and \eqref{phi_s}, it follows directly that $\theta_\jph^{\text{OE}}=1$ as the solution is at a discrete steady state. Furthermore, using \eqref{flux2Al}--\eqref{flux2Ar} and \eqref{low_res_b1}--\eqref{low_res_b2}, we obtain $\Delta\mbG_j^{(1)}=0$, $\Delta \mbPhir_\jmh^{(1)}=0$, and $\Delta \mbPhil_\jph^{(1)}=0$. Then, from \eqref{theta_j2} and \eqref{theta_jh2}, we obtain
\begin{equation*}
  \theta_j=1,\quad \theta_\jmh=1,\quad \theta_\jph=1.
\end{equation*}
\end{proof}
\end{remark}

\begin{prop}[Super-consistency for arbitrary moving equilibria]\label{prop3}
The WB PAMPA scheme with Gauss--Lobatto flux globalization preserves exactly discrete moving steady states obtained by integrating the nonlinear ODE
\begin{equation}\label{eq:flux-ode}
    \mbf'(\mbu)(x) = \mbS(\mbu(\mbf),x)
\end{equation}
with the LobattoIIIA collocation method when using \eqref{2.5},
or with the sc-LobattoIII method when using \eqref{eq:scLIII1}-\eqref{eq:scLIII2}. Provided $\mbf(\mbu)$ is a one to one mapping and the inverse $\mbu(\mbf)$ is bounded and uniquely defined, these discrete solutions are formally fourth order accurate.
\begin{proof}
We write the proof for the case of the LobattoIIIA method.
The one for the  sc-LobattoIII is almost identical.
Steady data obtained integrating \eqref{eq:flux-ode} with the  RK-LobattoIIIA method,
verify within each element the nonlinear algebraic system
\begin{equation*}
\begin{aligned}
\mbf_\jph - \mbf_{j} = \sum_{i=1}^3 \omega_i'\mbS(\mbu(\mbf_{j + \frac{i-1}{2} }),x_{j + \frac{i-1}{2} })\\
\mbf_{j+1} - \mbf_{j} = \sum_{i=1}^3 \omega_i\mbS(\mbu(\mbf_{j + \frac{i-1}{2} }),x_{j + \frac{i-1}{2} })\\
\end{aligned}
\end{equation*}
with $\{\omega_i',\omega_i\}_{i=1,3}$ the LobattoIIIA Butcher coefficients, also used in \eqref{2.5} to integrate the source. 
Combining \eqref{2.4} and \cref{rmk1} this leads to
\begin{equation*} 
\hmbG_{j+1}= \mbG_{j+1} = \mbG_{j} =\hmbG_j,  \quad \forall K_\jph \quad \Rightarrow    \quad
\frac{\mathrm{d}}{\mathrm{d}t}\xbar{\mbu}_\jph =0.
\end{equation*}
For the nodal updates, using \eqref{2.6b}--\eqref{2.6c} and \cref{rmk1}, 
we find easily that  
\begin{equation*}
\hmbPhir_\jmh=\hmbPhil_\jph=0 \quad \Rightarrow    \quad \frac{\mathrm{d}}{\mathrm{d}t}{\mbu}_{j}=0,\quad \forall j.
\end{equation*}
Thus the LobattoIIIA discrete solutions are in the kernel of the flux globalization based WB PAMPA scheme. The invertibility of $\mbf(\mbu)$  is enough to claim that the nodal accuracy of the LobattoIIIA is preserved, which is order 4  (cf. \cite{Hairer1993}). The proof for the case of  the  sc-LobattoIII method follows the same ideas, with the obvious changes in quadrature coefficients and points.
\end{proof}
\end{prop}

\begin{remark}[Shallow water: invertibility of $\mbf(\mbu)$ and critical point detection] Let us consider the shallow water flux and set $\mbf(\mbu)=(f_1,f_2,f_3)^\top$.  
For moving-water solutions, so $u\neq 0$,  and strictly positive water depth, the invertibility of the flux require detecting critical points where $u=\sqrt{gh}$. It is quite simple to verify that $f_2 >0$, and
\begin{equation*}
  \kappa:=g\dfrac{f_1^4}{f_2^3} =\dfrac{8 \mathrm{Fr}^4}{(2 \mathrm{Fr}^2 +1)^3}\;,\quad \mathrm{Fr}:= \dfrac{u}{\sqrt{gh}}.
\end{equation*}
This gives an immediate criterion to verify if a flux value corresponds to a critical point where  $\kappa=8/27$.
\end{remark}
As discussed  in \cite{Mantri2024,Kazolea2025} the static case is a particular one which requires a reformulation  of the consistency analysis.
As already said, for the shallow water equations, the sc-LobattoIII approach
has the advantage of always integrating exactly the hydrostatic term.
For this reason we now focus on this approach, which is also the one used in the numerical results. 

\begin{prop}[Preservation of still water states]\label{prop2}
The WB PAMPA scheme with Gauss-Lobatto global flux quadrature using the sc-LobattoIII approach 
\begin{enumerate}
\item preserves exactly still water states with $u=v=0$, $w=h+B=w_0$ constant, and with arbitrary topography;
\item preserves exactly the fourth order accurate discrete states at rest computed with the sc-LobattoIII collocation method applied to
\begin{equation*}
  f_2'(h) = -ghB_x+ fh v
\end{equation*}
with given $v(x)\neq 0$, and for arbitrary  topography.
\end{enumerate}
\begin{proof} To prove the first part we note that the relation
\begin{equation*}
\frac{\mathrm{d}}{\mathrm{d}t}\xbar{\mbu}_\jph=\frac{\mathrm{d}}{\mathrm{d}t}{\mbu}_j=0
\end{equation*}
is trivial for the first component and last component (appear in \eqref{1.3a}). We thus focus on the expression of $\mbG$ for the second component for which we can write
\begin{equation*}
  \mbG_{j}^{(2)}-\mbG_{\jmh}^{(2)} =\int_{x_{\jmh}}^{x_j} \Big(g(\frac{h^2}{2})_x + ghB_x \Big)\mathrm dx=g\dfrac{h^2_j - h^2_\jmh }{2} + \int_{x_\jmh}^{x_j }g h B_x\; \mathrm dx
\end{equation*}
and similarly for the $\mbG_{\jph}$. For the still-water steady-state solutions satisfying $g(\frac{h^2}{2})_x + ghB_x=0$. So, once the last integral is evaluated exactly by the Gauss--Lobatto formula, we immediately have $ \mbG_{j}^{(2)}-\mbG_{\jmh}^{(2)}=0$. This is achieved by using the sub-cell LobattoIII method in \eqref{eq:scLIII1}--\eqref{eq:scLIII2}, which is exactly for the polynomial $hB_x$ of degree $3$.

Note that when the LobattoIIIA method is applied, in order to exactly preserve the still-water states, one can use the same modifications of the hydrostatic source
proposed in \cite{Xing2005,Mantri2024,Kazolea2025} to replace the part of the source terms as
\begin{equation*}
    -ghB_x=-gwB_x+g\big(\frac{B^2}{2}\big)_x.
\end{equation*}
So we can write \emph{with no approximation}
\begin{equation*}
  \int_{x_\jmh}^{x_j }g w B_x \; \mathrm dx -\int_{x_\jmh}^{x_j }g  (\frac{B^2}{2})_x \;\mathrm dx  =
\int_{x_\jmh}^{x_j }g w B_x \;\mathrm dx  - g\dfrac{B^2_j- B^2_\jmh }{2}
\end{equation*}
One easily checks that for $w=h+B=w_0$ constant
\begin{equation*}
  \begin{aligned}
\mbG_{j}^{(2)}-\mbG_{\jmh}^{(2)} = &g\dfrac{h^2_j - h^2_\jmh }{2} + \int_{x_\jmh}^{x_j }g w_0 B_x \;\mathrm dx
 - g\dfrac{B^2_j - B^2_\jmh }{2}\\
 =& g\dfrac{h^2_j - h^2_\jmh}{2} + gw_0 (B_j - B_\jmh) - g\dfrac{B^2_j - B^2_\jmh }{2} =0.
 \end{aligned}
\end{equation*}
We  invoke  \cref{rmk1} to deduce that $\hmbG_{j+1}^{(2)}=\hmbG_j^{(2)}=\hmbG_\jph^{(2)}=\mbox{constant}$, and obtain the result sought.

The second statement of the proposition is a special case of \cref{prop3}, considering only the second
component of the flux $f_2=gh^2/2$ which is always invertible. The proof is identical to that of \cref{prop3}.
\end{proof}
\end{prop}

\section{Numerical examples}\label{sec4}
We have thoroughly tested the flux globalization based PP WB PAMPA scheme proposed on the Saint--Venant system \eqref{1.3} to assess convergence property, verify the PP and WB properties presented in the previous section, as well as to investigate its application to resolution of Riemann problems and to test its robustness. The benchmarks used are quite classical in literature. They allow to verify the correct implementation and convergence of the proposed method. 

The time-dependent ODE systems \eqref{2.2a} and \eqref{2.2b} are integrated using the three-stage third-order strong stability preserving (SSP) Runge-Kutta method (see, e.g.,\cite{Gottlieb2011,Gottlieb2001}). All numerical simulations have been run at CFL=0.2. Unless specified differently, we use the zero-order extrapolation boundary conditions and take the acceleration due to gravity $g=9.812$. In the following examples, the numerical results obtained using the developed high-order schemes are denoted as ``HO-PAMPA'' and the corresponding low-order schemes alone are denoted as ``LO-PAMPA''.

In Examples 1--7, the PP WB PAMPA scheme is applied to solve the Saint--Venant system with/without Manning friction term \eqref{1.3}, while in Examples 8--11, the proposed scheme is applied to solve the rotating shallow water equations \eqref{1.3a}.

\subsection*{Example 1---Accuracy Test}
In the first example, a slight modification from \cite{Abgrall2024} with the Manning coefficient $n=0.05$, we study a problem with smooth solutions in a finite time frame to check the experimental rate of convergence. The initial conditions and non-flat bottom topography are:
\begin{equation*}
  h(x,0)=0.3\Big[1+e^{-\frac{(x-0.5)^2}{0.05^2}}\Big]-0.2\cos(6\pi x),\quad u(x,0)=0,\quad B(x)=0.2(1+\cos(6\pi x)),
\end{equation*}
prescribed in a computational domain $[0,1]$ and subject to the periodic boundary conditions.

We compute the numerical solution until the final time $t =0.03$ using the high-order PP WB PAMPA scheme on a sequence of uniform meshes with $\dx=1/64$, $1/128$, $1/256$, $1/512$, $1/1024$, $1/2048$, and $1/4096$. After that, we measure the discrete $L^1$-errors for the cell boundary point and average values of variables $h$ and $hu$ using the Runge formulae provided in \cite{Abgrall2024}. Equipped with these errors, we compute the corresponding convergence rates and report the results in Table \ref{tab1a} for the point values and cell averages, respectively. From these obtained results, we can clearly see that the expected experimental third order of accuracy is achieved for the proposed schemes.

%\begin{table}[ht!]
%\caption{\sf Example 1: $L^1$-errors and experimental convergence rates of the cell boundary point values and average values of variables $h$ and $hu$ with the low order schemes in \textbf{Choice A}.\label{tab1}}
%\begin{center}
%\begin{tabular}{|c|c|c|c|c|}
%\hline
%\multicolumn{5}{|c|}{errors and rates for $\mbu_j$}\\
%\hline
%\multicolumn{1}{|c|}{$\dx$}&\multicolumn{1}{c|}{$L^1$-error in $h$}&{Rate}&\multicolumn{1}{c|}{$L^1$-error in $hu$}&{Rate}\\
%\hline
%$1/256$&$7.08\times 10^{-4}$&-&$1.35\times 10^{-3}$&-\\
%\hline
%$1/512$&$1.30\times 10^{-4}$&2.44&$2.56\times 10^{-4}$&2.40\\
%\hline
%$1/1024$&$1.86\times 10^{-5}$&2.80&$3.70\times 10^{-5}$&2.79\\
%\hline
%$1/2048$&$2.09\times 10^{-6}$&3.15&$4.15\times 10^{-6}$&3.15\\
%\hline
%$1/4096$&$2.30\times 10^{-7}$&3.18&$4.54\times 10^{-7}$&3.19\\
%\hline
%\multicolumn{5}{|c|}{errors and rates for $\xbar \mbu_\jph$}\\
%\hline
%\multicolumn{1}{|c|}{$\dx$}&\multicolumn{1}{c|}{$L^1$-error in $h$}&{Rate}&\multicolumn{1}{c|}{$L^1$-error in $hu$}&Rate\\
%\hline
%$1/256$&$8.42\times 10^{-4}$&-&$1.70\times 10^{-3}$&-\\
%\hline
%$1/512$&$1.36\times 10^{-4}$&2.63&$2.92\times 10^{-4}$&2.54\\
%\hline
%$1/1024$&$1.75\times 10^{-5}$&2.96&$3.75\times 10^{-5}$&2.96\\
%\hline
%$1/2048$&$1.88\times 10^{-6}$&3.21&$4.04\times 10^{-6}$&3.21\\
%\hline
%$1/4096$&$2.02\times 10^{-7}$&3.22&$4.33\times 10^{-7}$&3.22\\
%\hline
%\end{tabular}
%\end{center}
%\end{table}

\begin{table}[ht!]
\caption{\sf Example 1: $L^1$-errors and experimental convergence rates of the cell boundary point values and average values of variables $h$ and $hu$.\label{tab1a}}
\begin{center}
\begin{tabular}{|c|c|c|c|c|}
\hline
\multicolumn{5}{|c|}{errors and rates for $\mbu_j$}\\
\hline
\multicolumn{1}{|c|}{$\dx$}&\multicolumn{1}{c|}{$L^1$-error in $h$}&{Rate}&\multicolumn{1}{c|}{$L^1$-error in $hu$}&{Rate}\\
\hline
$1/256$&$4.16\times 10^{-4}$&-&$7.37\times 10^{-4}$&-\\
\hline
$1/512$&$9.69\times 10^{-5}$&2.10&$1.80\times 10^{-4}$&2.04\\
\hline
$1/1024$&$1.34\times 10^{-5}$&2.85&$2.59\times 10^{-5}$&2.79\\
\hline
$1/2048$&$1.66\times 10^{-6}$&3.01&$3.24\times 10^{-6}$&3.00\\
\hline
$1/4096$&$2.02\times 10^{-7}$&3.04&$3.94\times 10^{-7}$&3.04\\
\hline
\multicolumn{5}{|c|}{errors and rates for $\xbar \mbu_\jph$}\\
\hline
\multicolumn{1}{|c|}{$\dx$}&\multicolumn{1}{c|}{$L^1$-error in $h$}&{Rate}&\multicolumn{1}{c|}{$L^1$-error in $hu$}&Rate\\
\hline
$1/256$&$4.72\times 10^{-4}$&-&$1.01\times 10^{-3}$&-\\
\hline
$1/512$&$8.21\times 10^{-5}$&2.52&$1.78\times 10^{-4}$&2.51\\
\hline
$1/1024$&$1.14\times 10^{-5}$&2.85&$2.46\times 10^{-5}$&2.86\\
\hline
$1/2048$&$1.40\times 10^{-6}$&3.02&$3.01\times 10^{-6}$&3.03\\
\hline
$1/4096$&$1.70\times 10^{-7}$&3.04&$3.66\times 10^{-7}$&3.04\\
\hline
\end{tabular}
\end{center}
\end{table}

\subsection*{Example 2---Still Water Equilibrium and Its Small Perturbation}
In the second example, taken from \cite{Abgrall2024}, we examine the still water preservation and positivity-preserving properties of the proposed high-order and first-order schemes. We consider a non-flat bottom topography,
\begin{equation*}
  B(x)=\left\{\begin{aligned}
  &2(\cos(10\pi(x+0.3))+1),\quad && x\in[-0.4,-0.2],\\
  &0.5(\cos(10\pi(x-0.3))+1),\quad && x\in[0.2,0.4],\\
  &0,\quad &&\mbox{otherwise},
  \end{aligned}\right.
\end{equation*}
and the following initial data
\begin{equation*}
  h(x,0)=4.000001-B(x), \quad u(x,0)=0,
\end{equation*}
in the computational domain $[-1,1]$. The Manning friction coefficient is $n=0$. We run the simulations until a long final time $t=10$ using the schemes proposed with $50$ uniform cells and report the discrete $L^1$-, $L^2$-, and $L^\infty$-errors of the average values of $h$ and $hu$ in Table \ref{tab2}. All the calculated errors are within the level of round-off error, clearly validating \cref{prop:still_low} and \cref{prop2}.   

\begin{table}[!ht]
\caption{\sf Example 2: Errors in $h$ and $hu$ computed by high-order and first-order schemes.\label{tab2}}
\begin{center}
\begin{tabular}{c| c| c c c}\hline
\multicolumn{1}{c|}{\multirow{1}{*}{Schemes}} &\multicolumn{1}{c|}{\multirow{1}{*}{Variables}} &$L^1$-error &$L^2$-error&$L^\infty$-error\\ \hline
 \multicolumn{1}{c|}{\multirow{2}{*}{HO-PAMPA}} &$h$ &$2.77\times10^{-14}$ & $2.44\times10^{-14}$  & $3.02\times10^{-14}$ \\ \cline{2-5}
 &$hu$ &$2.20\times10^{-13}$ & $5.96\times10^{-13}$  & $2.93\times10^{-12}$ \\ \hline
  \multicolumn{1}{c|}{\multirow{2}{*}{LO-PAMPA}} &$h$ &$1.29\times10^{-15}$ & $3.51\times10^{-15}$  & $1.02\times10^{-14}$ \\  \cline{2-5}
 &$hu$ &$1.98\times10^{-14}$ & $1.77\times10^{-14}$  & $2.71\times10^{-14}$  \\ \hline
\end{tabular}
\end{center}
\end{table}

%\begin{table}[!ht]
%\caption{\sf Example 2: Errors in $h$ and $hu$ computed by high-order and first-order schemes.\label{tab2}}
%\begin{center}
%\begin{tabular}{c| c| c c c}\hline
%\multicolumn{1}{c|}{\multirow{1}{*}{Schemes}} &\multicolumn{1}{c|}{\multirow{1}{*}{Variables}} &$L^1$-error &$L^2$-error&$L^\infty$-error\\ \hline
% \multicolumn{1}{c|}{\multirow{2}{*}{HO-PAMPA-A}} &$h$ &$4.88\times10^{-16}$ & $2.08\times10^{-15}$  & $1.02\times10^{-14}$ \\ \cline{2-5}
% &$hu$ &$5.74\times10^{-14}$ & $1.02\times10^{-13}$  & $4.79\times10^{-13}$ \\ \hline
% \multicolumn{1}{c|}{\multirow{2}{*}{HO-PAMPA-B}} &$h$ &$6.02\times10^{-13}$ & $4.25\times10^{-13}$  & $3.10\times10^{-13}$ \\ \cline{2-5}
% &$hu$ &$5.23\times10^{-14}$ & $4.52\times10^{-14}$  & $6.03\times10^{-14}$ \\ \hline
% \multicolumn{1}{c|}{\multirow{2}{*}{LO-PAMPA-A}} &$h$ &$1.29\times10^{-15}$ & $3.51\times10^{-15}$  & $1.02\times10^{-14}$ \\  \cline{2-5}
% &$hu$ &$7.29\times10^{-14}$ & $5.78\times10^{-14}$  & $7.06\times10^{-14}$  \\ \hline
%  \multicolumn{1}{c|}{\multirow{2}{*}{LO-PAMPA-B}} &$h$ &$1.29\times10^{-15}$ & $3.51\times10^{-15}$  & $1.02\times10^{-14}$ \\  \cline{2-5}
% &$hu$ &$1.98\times10^{-14}$ & $1.77\times10^{-14}$  & $2.71\times10^{-14}$  \\ \hline
%\end{tabular}
%\end{center}
%\end{table}

\begin{figure}[ht!]
\centerline{\subfigure[$t=0.02$]{\includegraphics[trim=1.0cm 0.05cm 0.5cm 0.05cm,clip,width=4.25cm]{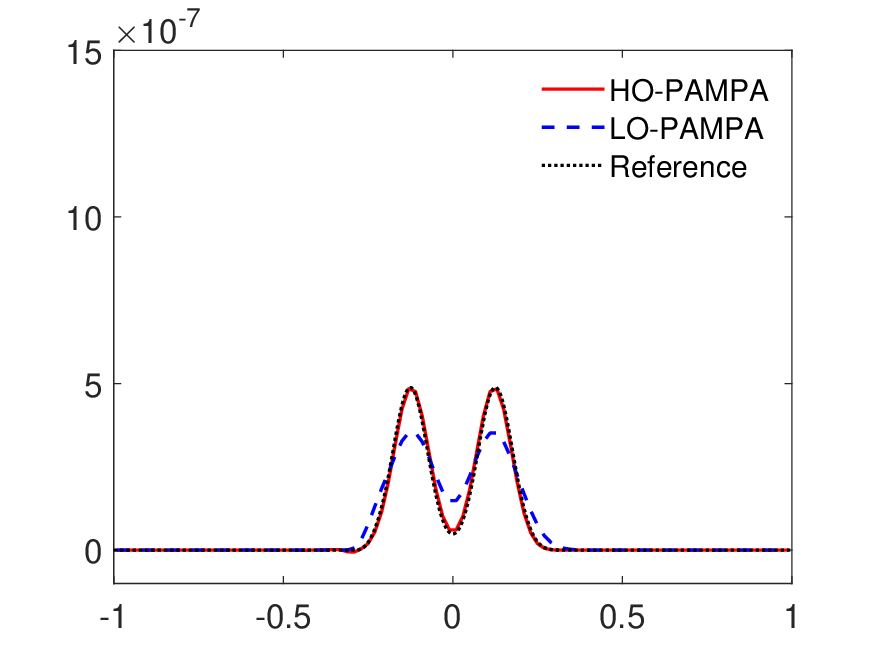}}\hspace*{0.05cm}
	\subfigure[$t=0.04$]{\includegraphics[trim=1.0cm 0.05cm 0.5cm 0.05cm,clip,width=4.25cm]{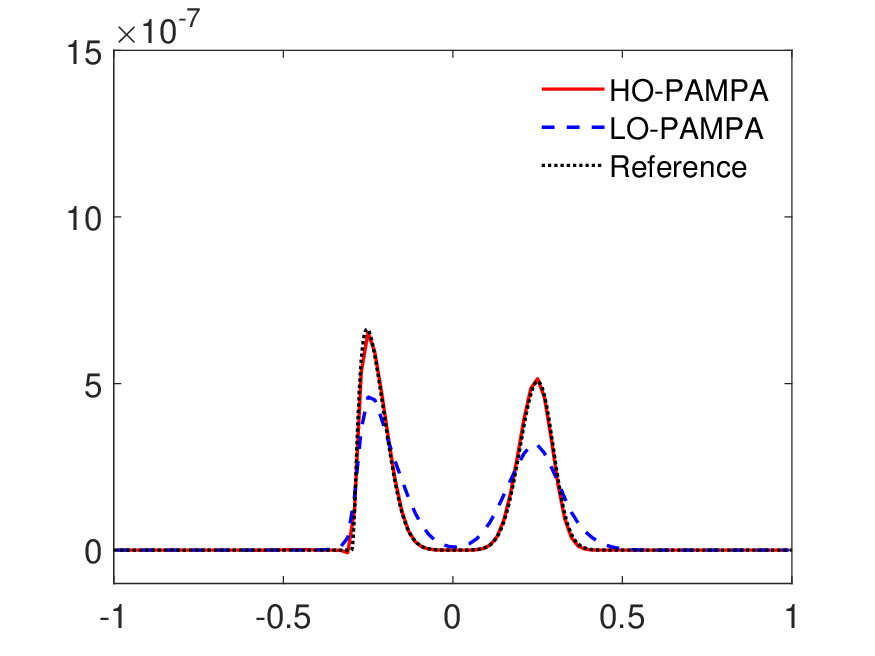}}\hspace*{0.05cm}
	\subfigure[$t=0.06$]{\includegraphics[trim=1.0cm 0.05cm 0.5cm 0.05cm,clip,width=4.25cm]{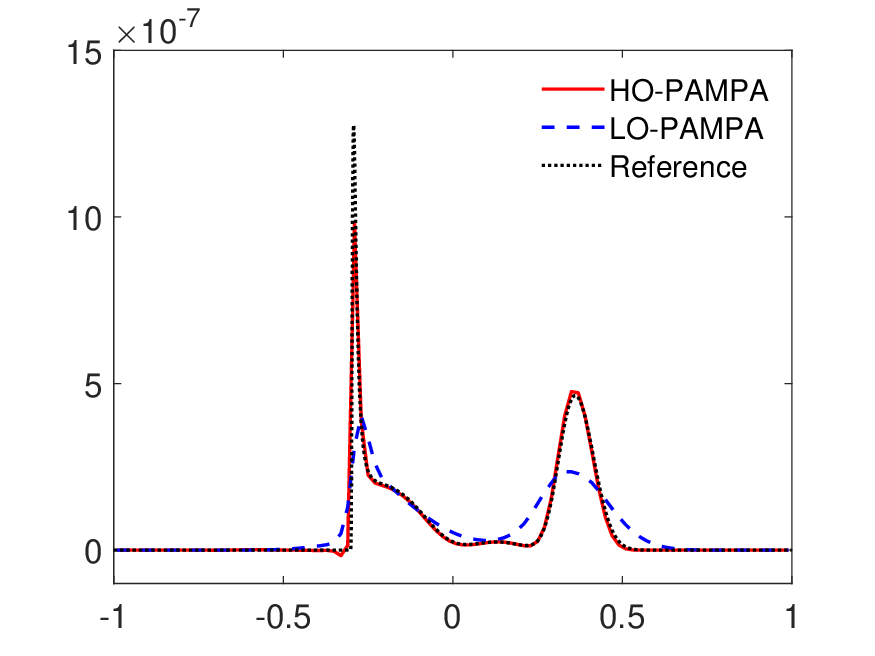}}}
\caption{\sf Example 2: Time snapshots of the small perturbation in $h$.\label{Ex2_dh}}
\end{figure}

To further verify the WB property for preserving still-water equilibrium, we add a small Gaussian-shaped perturbation to the stationary water depth and test the ability of the proposed schemes to accurately capture this small perturbation. The small perturbation of the form $10^{-6}\times e^{-200x^2}$ is superimposed onto the initial cell average of water depth. We compute the numerical solution at three different times: $t=0.02$, $0.04$, and $0.06$, using $100$ uniform cells. A reference solution computed by the first-order scheme on a uniform mesh with $2000$ cells is also provided. The differences between the computed and background cell averages of water depth are plotted in Figure \ref{Ex2_dh}. 

As one can observe, both high-order and low-order PP WB PAMPA schemes can capture the dynamics of the small perturbation, there are no spurious oscillations: the initial perturbation splits into two humps, which are then propagating into opposite directions, respectively. The results generated on coarse and fine meshes are consistent and converge to the same solution profile. %We present only the numerical results obtained using the high-order schemes blended with the low-order schemes of \textbf{Choice B}. When the low-order schemes of \textbf{Choice A} are used, the results coincide with those already reported; hence, they are omitted here for brevity.

%
%\begin{figure}[ht!]
%\centerline{\subfigure[average values: $\xbar h_\jph+\xbar B_\jph$ and $\xbar B_\jph$]{\includegraphics[trim=1.0cm 0.05cm 0.5cm 0.05cm,clip,width=6.0cm]{figures/Ex2_surface_Bottom}}\hspace*{0.25cm}
%	\subfigure[point values: $h_j+B_j$ and $B_j$]{\includegraphics[trim=1.0cm 0.05cm 0.5cm 0.05cm,clip,width=6.0cm]{figures/Ex2_surface_Bottomp}}}
%\caption{\sf Example 2: Values of water surface and topography at time $t=0.06$.\label{Ex2_w_B}}
%\end{figure}

%Finally, we measure the difference in CPU times between the monolithic convex limiting PP WB PAMPA scheme proposed here and a posteriori MOOD paradigm based WB PAMPA scheme of \cite{Abgrall2024} for simulating the perturbation propagation up to the final time $t=0.06$. The results show that the a priori limiting PP WB PAMPA scheme reduces the computational cost by approximately $17.6\%$ compared to the MOOD-based approach. This improvement is expected, since the a priori limiter avoids the repeated admissibility checks and solution recomputations for troubled cells required in the MOOD loop.

\subsection*{Example 3---Super-convergent to Moving Water Equilibria without fri\-ction}
\begin{figure}[ht!]
\centerline{\subfigure[Case (I): surface $\xbar h+\xbar B$]{\includegraphics[trim=0.01cm 0.05cm 0.5cm 0.05cm,clip,width=4.25cm]{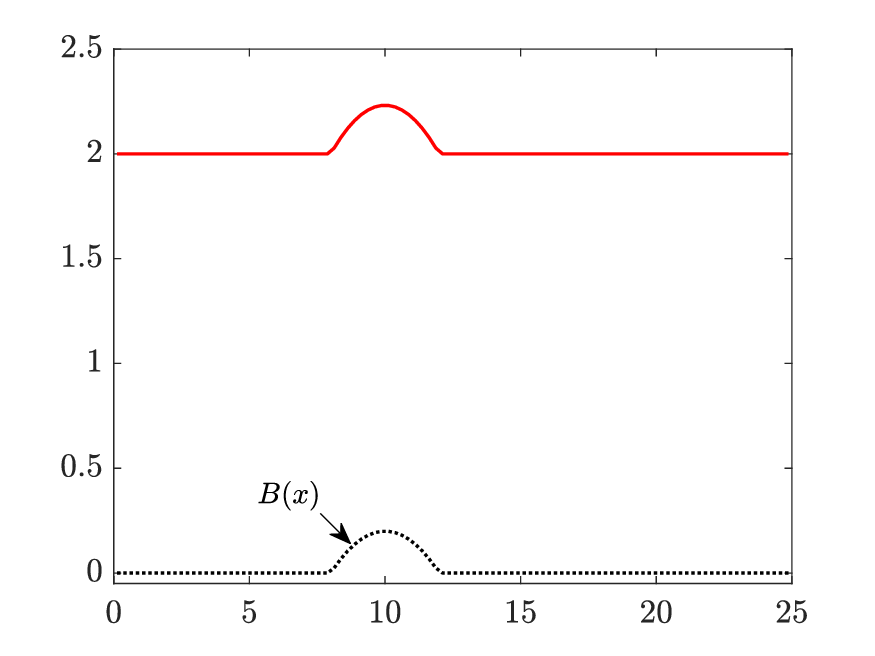}}\hspace*{0.5cm}
	\subfigure[Case (I): momentum $\xbar{(hu)}$]{\includegraphics[trim=0.01cm 0.05cm 0.5cm 0.05cm,clip,width=4.25cm]{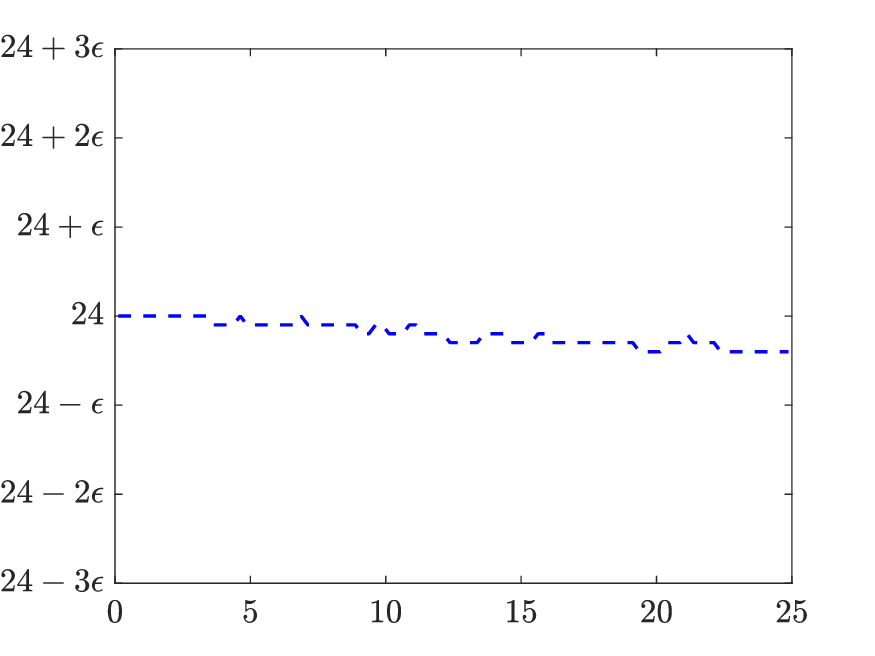}}\hspace*{0.5cm}
	\subfigure[Case (I): Global flux $\mbG^{(2)}$]{\includegraphics[trim=0.01cm 0.05cm 0.5cm 0.05cm,clip,width=4.25cm]{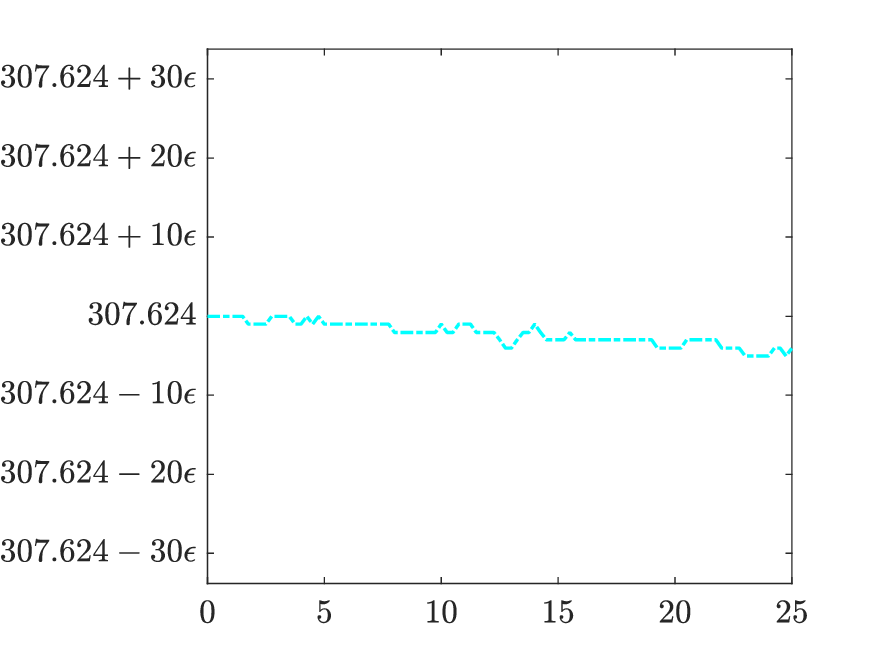}}}
\caption{\sf Example 3: Convergent solutions for supercritical flow. The deviation parameter $\epsilon =10^{-12}$.\label{Ex3a}}
\end{figure}

\begin{figure}[ht!]
\centerline{\subfigure[Case (II): surface $\xbar h+\xbar B$]{\includegraphics[trim=0.01cm 0.05cm 0.5cm 0.05cm,clip,width=4.25cm]{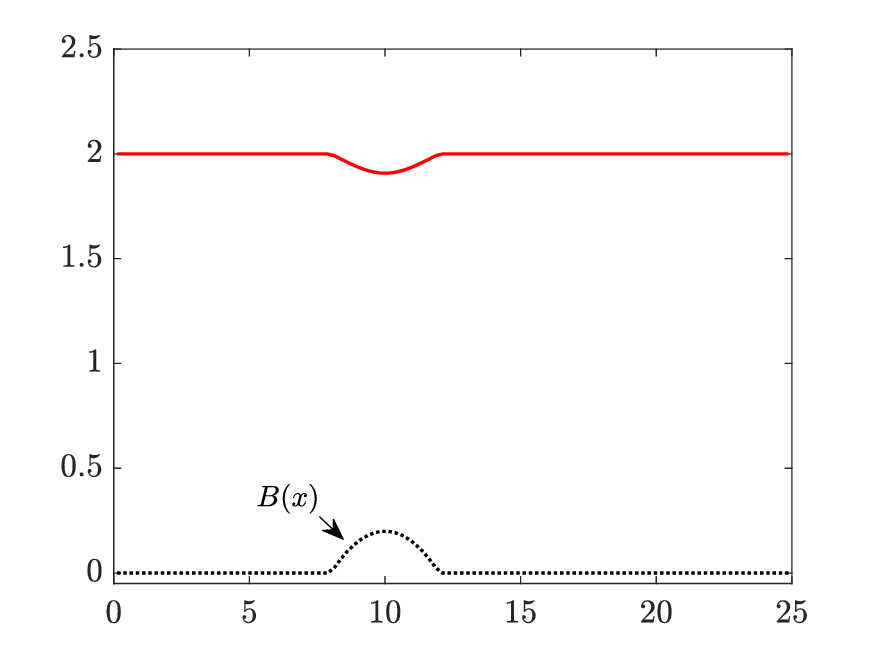}}\hspace*{0.5cm}
	\subfigure[Case (II): momentum $\xbar{(hu)}$]{\includegraphics[trim=0.01cm 0.05cm 0.5cm 0.05cm,clip,width=4.25cm]{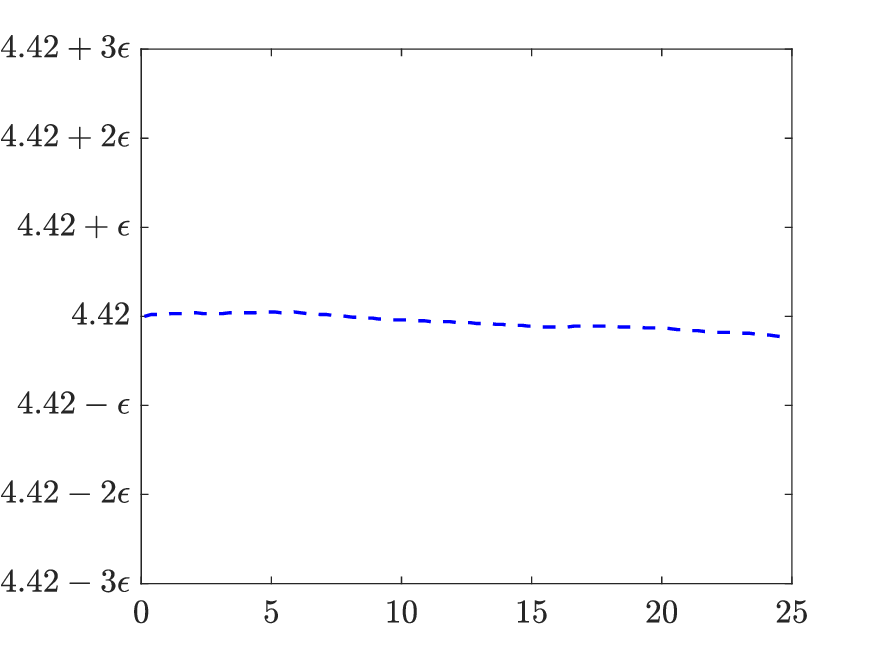}}\hspace*{0.5cm}
	\subfigure[Case (II): Global flux $\mbG^{(2)}$]{\includegraphics[trim=0.01cm 0.05cm 0.5cm 0.05cm,clip,width=4.25cm]{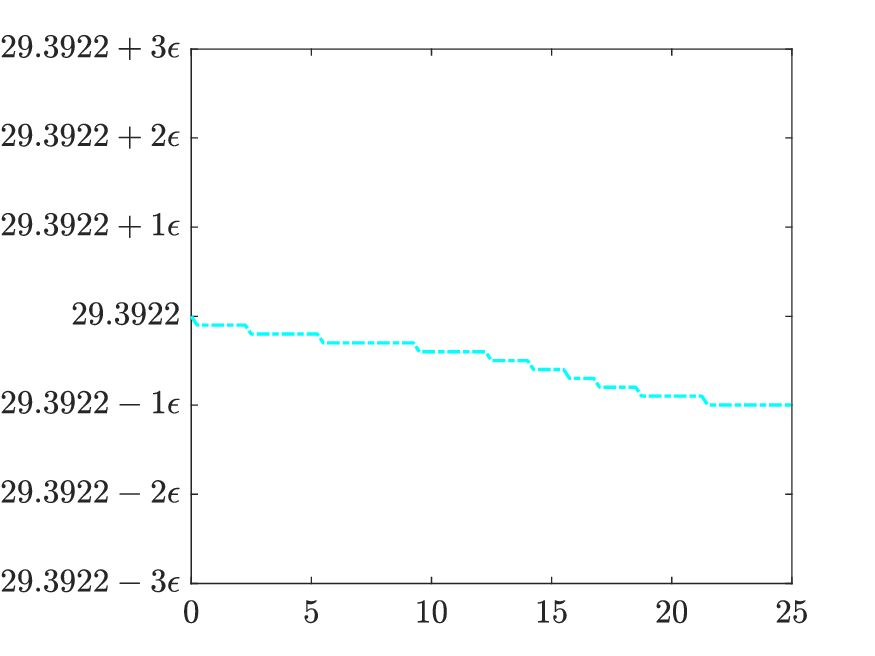}}}
\caption{\sf Example 3: Convergent solutions for subcritical flow. The deviation parameter $\epsilon =10^{-12}$.\label{Ex3b}}
\end{figure}

In the third example, we study the convergence of the numerical solution computed by the proposed PP WB PAMPA scheme towards steady flow over a continuous bottom topography given by
\begin{equation}\label{Ex3_bottom}
  B(x)=\left\{
  \begin{aligned}
  &0.2-0.05(x-10)^2,&&\mbox{if~}8\leq x\leq 10,\\
  &0,&&\mbox{else}.
  \end{aligned}\right.
\end{equation}
The goal here is to verify the theoretical predictions on the super convergence on the smooth moving-water equilibria for the scheme. We compute the convergent flow solutions for two classical smooth states: supercritical and subcritical, which depend on the following initial and boundary conditions:
\begin{equation*}
   \begin{aligned}
&\mbox{Case (I):}&&\left\{\begin{array}{l}h(x,0)=2-B(x),\quad hu(x,0)=0,\\h(0,t)=2,\qquad\qquad~~ hu(0,t)=24;\end{array}\right.\\
&\mbox{Case (II):}&&\left\{\begin{array}{l}h(x,0)=2-B(x),\quad hu(x,0)=0,\\hu(0,t)=4.42,\qquad~~ h(25,t)=2.\end{array}\right.
\end{aligned}
\end{equation*}
We let the scheme evolve until a given finite time $t=500$. The numerical solutions computed on a uniform mesh with $100$ cells are presented in Figures \ref{Ex3a} and \ref{Ex3b} for supercritical and subcritical flows, respectively. As one can see, the water depth obtained by the proposed PP WB PAMPA scheme is very close to the corresponding steady states (reported in \cite[Example 3]{Abgrall2024}). The
global flux $\mbG$ with the first component being the momentum $hu$ have reached to the expected constant states. This indicates that the numerical solutions have converged to the discrete steady state satisfying \eqref{1.2b}. 

\begin{figure}[ht!]
\centerline{\subfigure[supercritical flow]{\includegraphics[trim=1.0cm 0.05cm 0.5cm 0.05cm,clip,width=4.25cm]{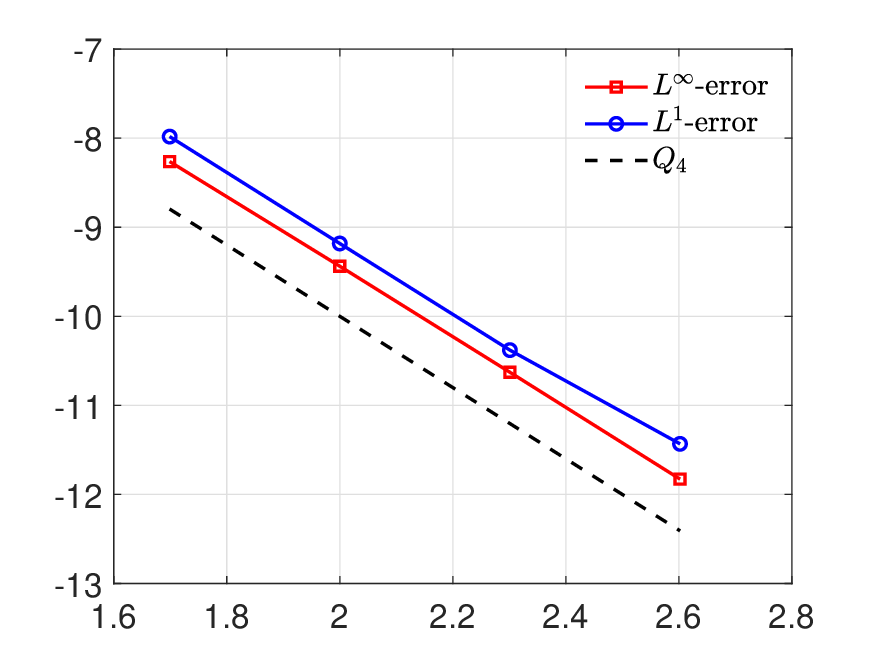}}\hspace*{0.5cm}
	\subfigure[subcritical flow]{\includegraphics[trim=1.0cm 0.05cm 0.5cm 0.05cm,clip,width=4.25cm]{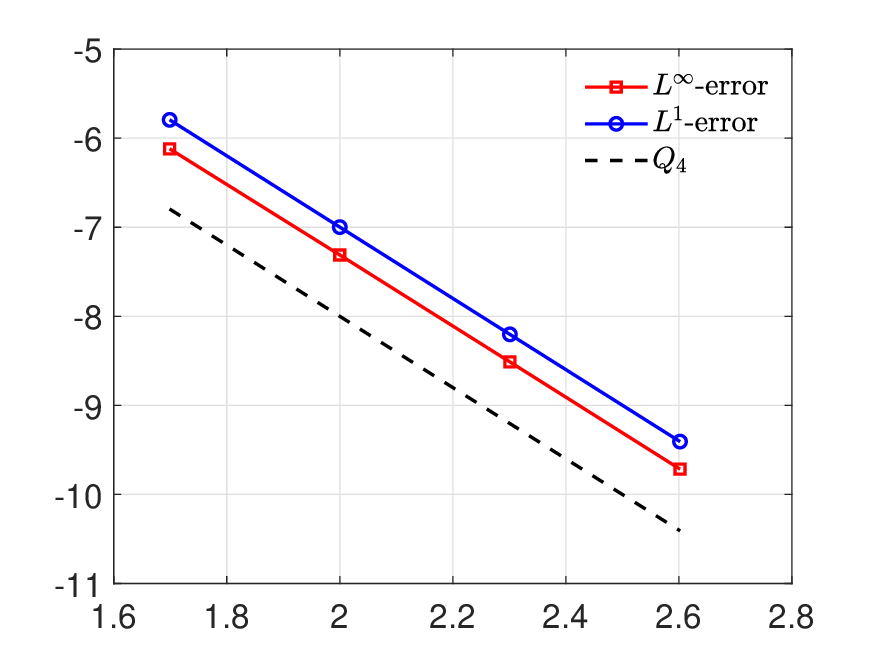}}}
\caption{\sf Example 3: logarithmic plots for errors in water depth versus number of cells.\label{Ex3}}
\end{figure}

Next, we qualify the convergence between the discrete water height obtained by this text and the analytical steady-state water height. To this end, we consider the mesh covered by a set of uniform cells with numbers: 50, 100, 200, and 400. We compute the discrete $L^1$- and $L^\infty$-errors between the discrete solutions and the analytical steady states given in \cite[Example 3]{Abgrall2024}. The logarithmic plots for errors in water depth versus number of cells are as shown in Figure \ref{Ex3}. We can see that the convergence rates are as predicted in \cref{prop3}, i.e., the convergence is of order four for the third-order PP WB PAMPA scheme.

\subsection*{Example 4---Discrete Moving Water Equilibria and Their Small Perturbations}
\begin{figure}[ht!]
\centerline{\subfigure[Case (I): surface $\xbar h+\xbar B$]{\includegraphics[trim=0.01cm 0.05cm 0.5cm 0.05cm,clip,width=4.15cm]{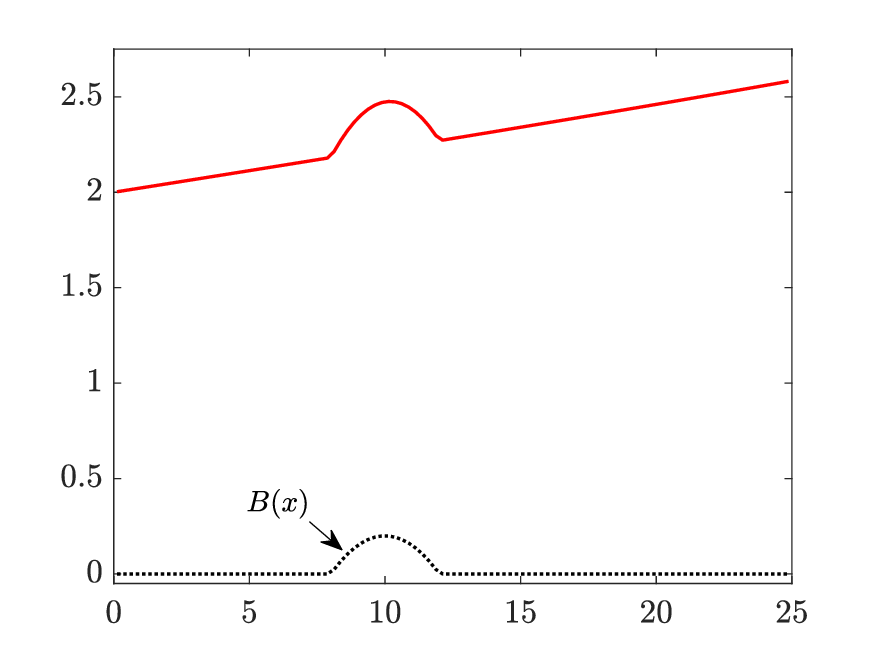}}\hspace*{0.05cm}
	\subfigure[Case (I): momentum $\xbar{(hu)}$]{\includegraphics[trim=0.01cm 0.05cm 0.5cm 0.05cm,clip,width=4.15cm]{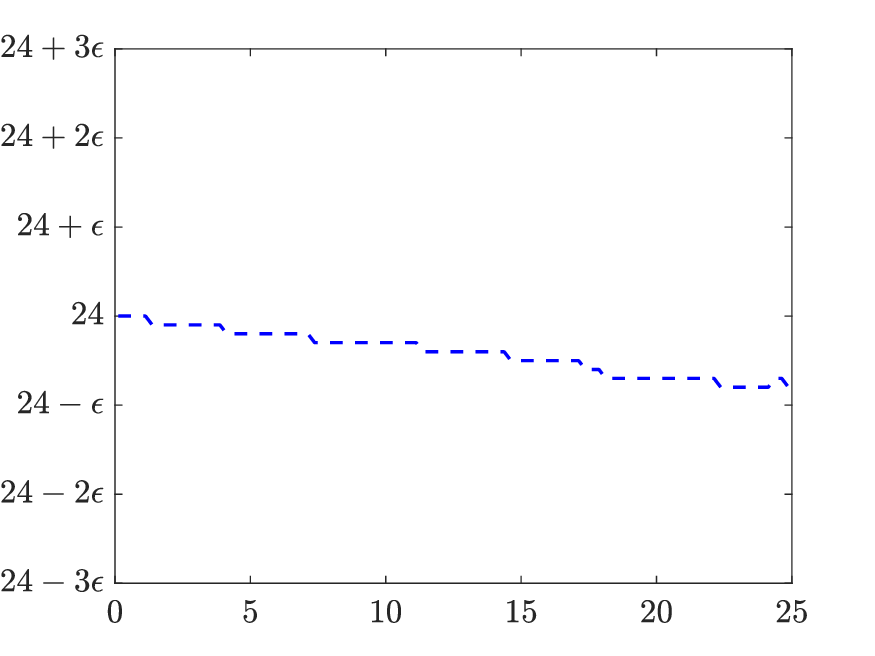}}\hspace*{0.05cm}
	\subfigure[Case (I): Global flux $\mbG^{(2)}$]{\includegraphics[trim=0.01cm 0.05cm 0.5cm 0.05cm,clip,width=4.15cm]{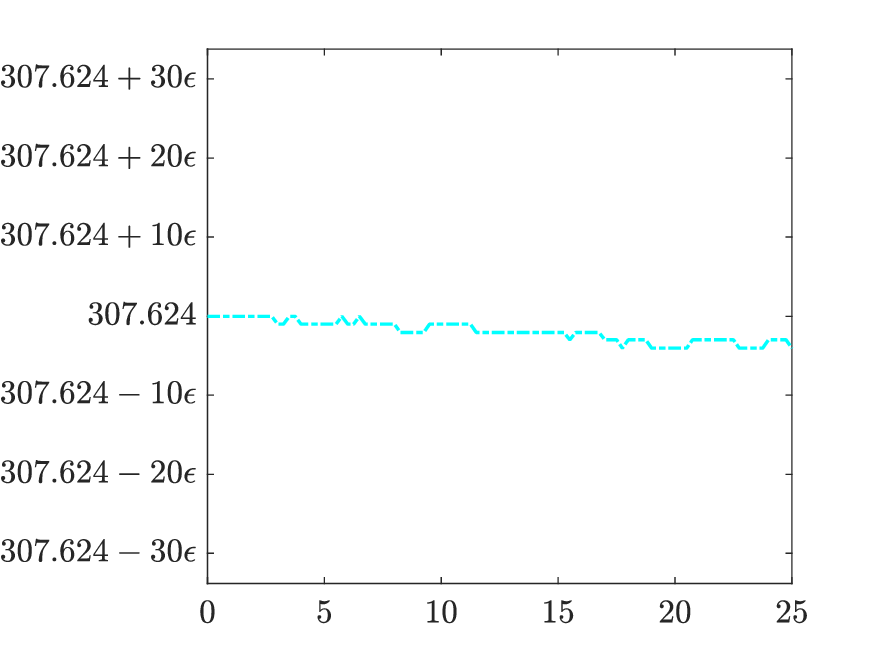}}}
\caption{\sf Example 4--Case (I): Convergent solutions for supercritical flow with friction. The deviation parameter $\epsilon =10^{-12}$.\label{Ex4a}}
\end{figure}

\begin{figure}[ht!]
\centerline{\subfigure[Case (II): surface $\xbar h+\xbar B$]{\includegraphics[trim=0.01cm 0.05cm 0.5cm 0.05cm,clip,width=4.15cm]{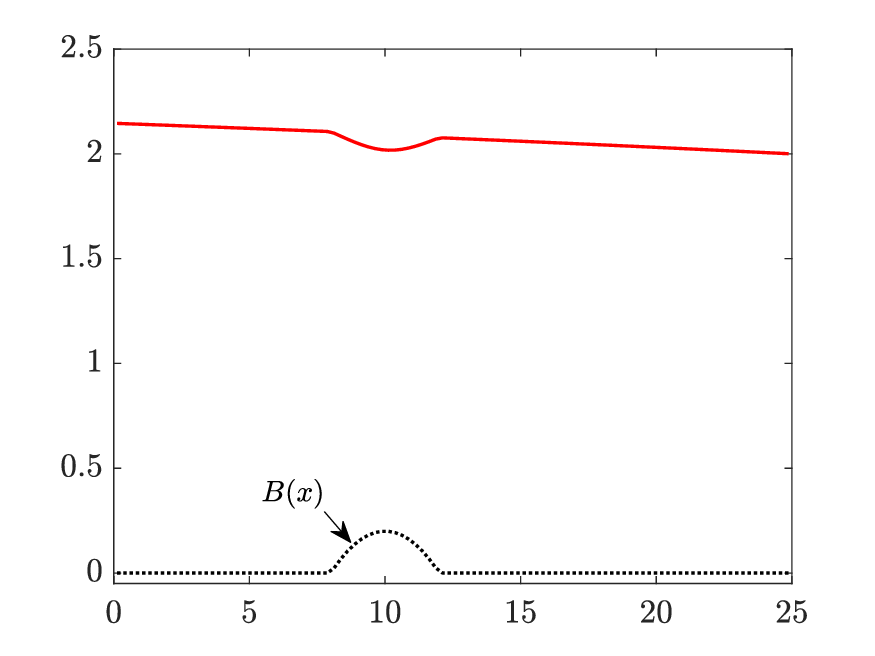}}\hspace*{0.05cm}
	\subfigure[Case (II): momentum $\xbar{(hu)}$]{\includegraphics[trim=0.01cm 0.05cm 0.5cm 0.05cm,clip,width=4.15cm]{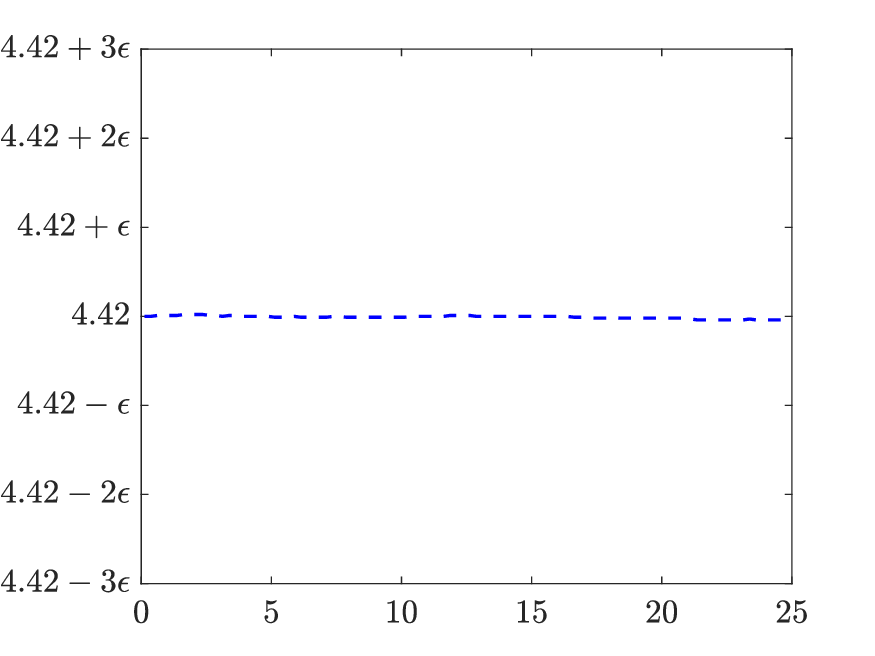}}\hspace*{0.05cm}
	\subfigure[Case (II): Global flux $\mbG^{(2)}$]{\includegraphics[trim=0.01cm 0.05cm 0.5cm 0.05cm,clip,width=4.15cm]{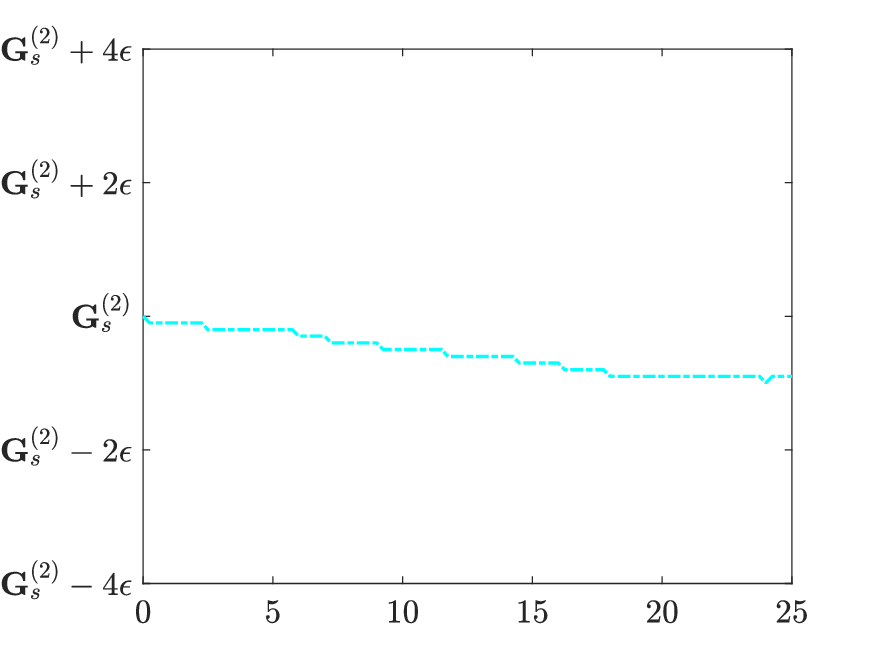}}}
\caption{\sf Example 4--Case (II): Convergent solutions for subcritical flow with friction. The deviation parameter $\epsilon =10^{-12}$ and the constant value $\mbG^{(2)}_s=31.700836562966$.\label{Ex4b}}
\end{figure}

In the fourth example, we consider a frictional case with Manning friction parameter $n=0.05$ and consider the same initial and boundary conditions as in Example 3. We plot the numerical solutions obtained by the proposed PP WB PAMPA scheme at $t=500$ using $100$ uniform grid cells in Figures \ref{Ex4a} and \ref{Ex4b}. One may notice that as in the frictionless case, both momentum $hu$ and global flux $\mbG^{(2)}$ computed by the proposed scheme are very accurate and have reached at the constant steady state. Moreover, with the influence of friction, the water depth and the water surface differ with those obtained in the frictionless cases. 

Next, we test the performance of the proposed PP WB PAMPA scheme in the presence of small perturbations of the moving water equilibria. We consider the following two sets of the steady-state data:
\begin{equation*}
\begin{aligned}
 &\mbox{Case (I):} \quad hu(x,0)=24,\quad \mbG^{(2)}=307.624;\\
 &\mbox{Case (II):} \quad hu(x,0)=4.42,\quad  \mbG^{(2)}=31.7008.
\end{aligned}
\end{equation*}
We note that the data are given in terms of the global quantities rather than in water depth. In order to start the computation at time $t=0$, one needs to obtain the values of $\xbar h_\jph$ and $h_j$. We, therefore, use Newton's method to numerically solve the following nonlinear equations:
\begin{equation}\label{eqs:nonlinear}
\begin{aligned}
  \mbG^{(2)}_j&=\frac{\big((hu)_j\big)^2}{h_j}+\frac{g}{2}(h_j)^2-\mbR^{(2)}_j,\\
   \mbG^{(2)}_\jph&=\frac{\big((hu)_\jph\big)^2}{(h_\jph)^2}+\frac{g}{2}(h_\jph)^2-\mbR^{(2)}_\jph,
\end{aligned}~\forall j,
\end{equation}
where $\mbR^{(2)}_j$ and $\mbR^{(2)}_\jph$ are given in \eqref{eq:scLIII1}. After obtaining the point values of water depth, it is easy to compute the initial
steady-state average value $\xbar h_\jph^{\mathrm eq}$ using Simpson's rule. Using these prepared data, we repeat the simulations from the still water case (Example 2) to verify the exact fully well-balanced property of the proposed HO-PAMPA scheme. We run the simulations using 100 cells until the final time $t=1000$, and compute the errors between the numerical and discrete prepared steady-state solutions. These errors are within machine precision and are omitted here for saving space. 

\begin{figure}[ht!]
\centerline{\subfigure[Case (I): $\xbar h(x,t)-\xbar h^{\mathrm eq}$]{\includegraphics[trim=1.0cm 0.05cm 0.5cm 0.05cm,clip,width=5.25cm]{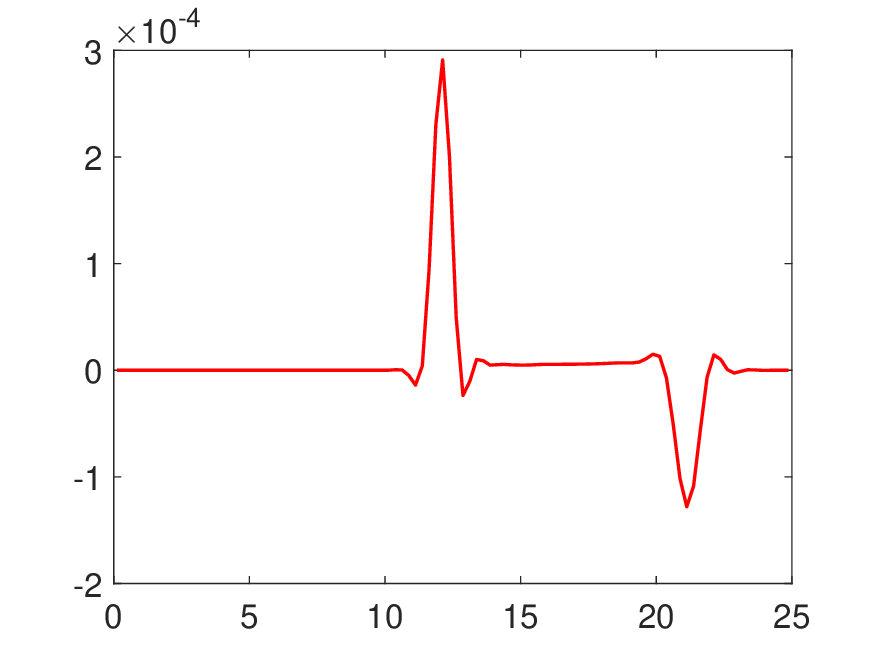}}\hspace*{0.5cm}
	\subfigure[Case (II): $\xbar h(x,t)-\xbar h^{\mathrm eq}$]{\includegraphics[trim=1.0cm 0.05cm 0.5cm 0.05cm,clip,width=5.25cm]{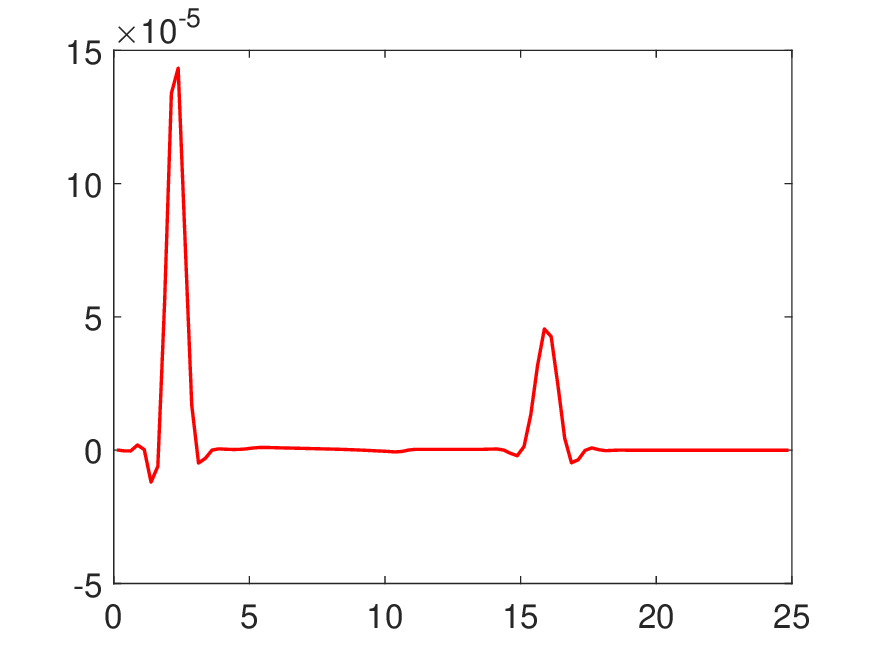}}}
\caption{\sf Example 4: The difference between the computed and the background discrete steady-state cell averages of water depth computed using 100 uniform cells.\label{Ex4na}}
\end{figure}

Finally, we add a small perturbation of the form $10^{-3}e^{-80(x-6)^2}$ to the steady-state water depth. We compute the numerical solutions at $t=1$ and $t=1.5$ for the Cases (I) and (II), respectively. The differences between the steady-state and numerical water depths are plotted in Figure \ref{Ex4na}. The results demonstrate that the proposed PP WB PAMPA scheme accurately captures the propagation of the small perturbation, in agreement with those reported in \cite[Example 8]{Abgrall2024}.  

\subsection*{Example 5---Riemann Problems over a Flat Bottom}
\begin{figure}[ht!]
\centerline{\subfigure[water height $\xbar h$ for Test 1]{\includegraphics[trim=1.0cm 0.05cm 0.5cm 0.05cm,clip,width=5.25cm]{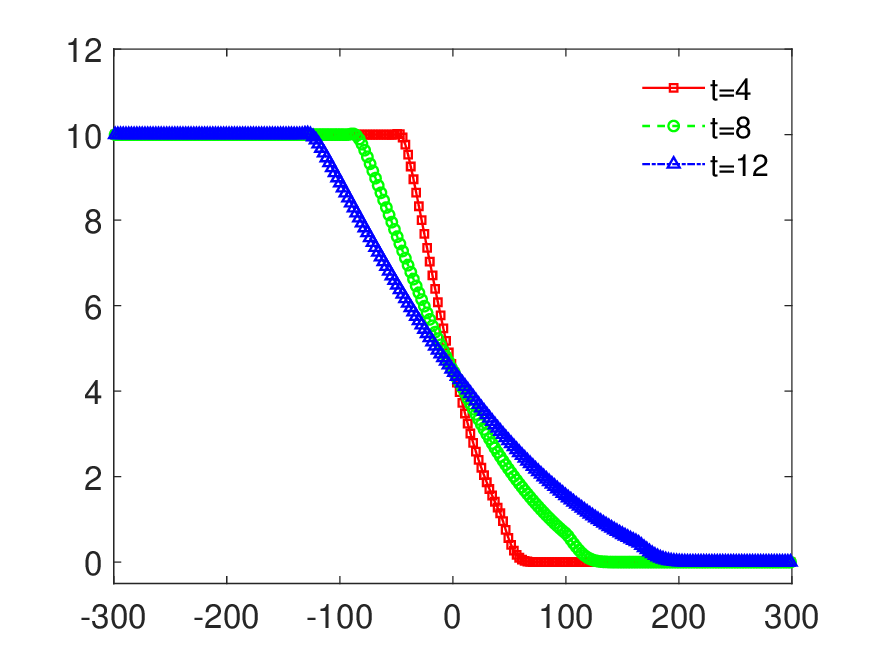}}\hspace*{0.5cm}
	\subfigure[momentum $\xbar{(hu)}$ for Test 1]{\includegraphics[trim=1.0cm 0.05cm 0.5cm 0.05cm,clip,width=5.25cm]{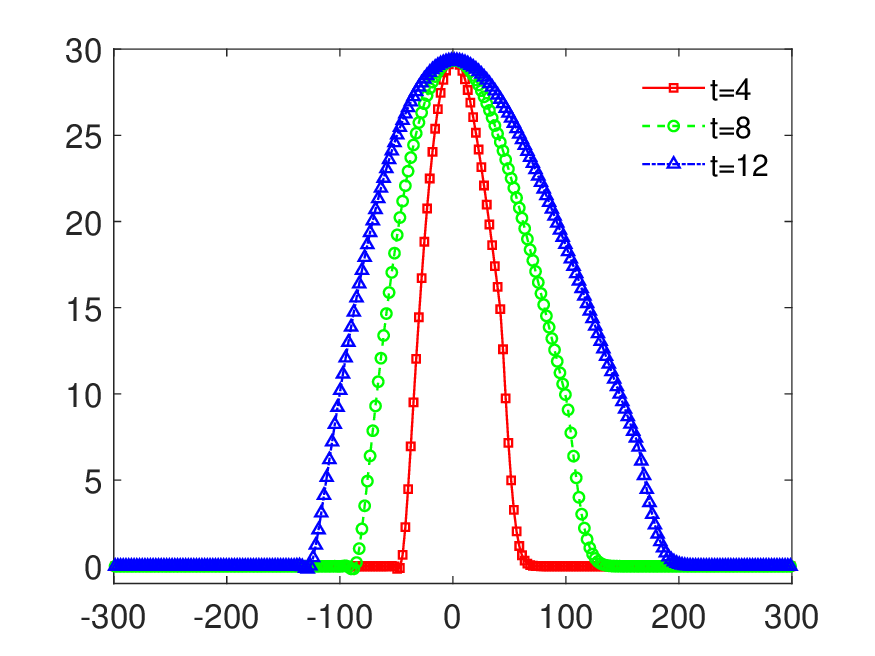}}}
\vskip5pt
\centerline{\subfigure[water height $\xbar h$ for Test 2]{\includegraphics[trim=1.0cm 0.05cm 0.5cm 0.05cm,clip,width=5.25cm]{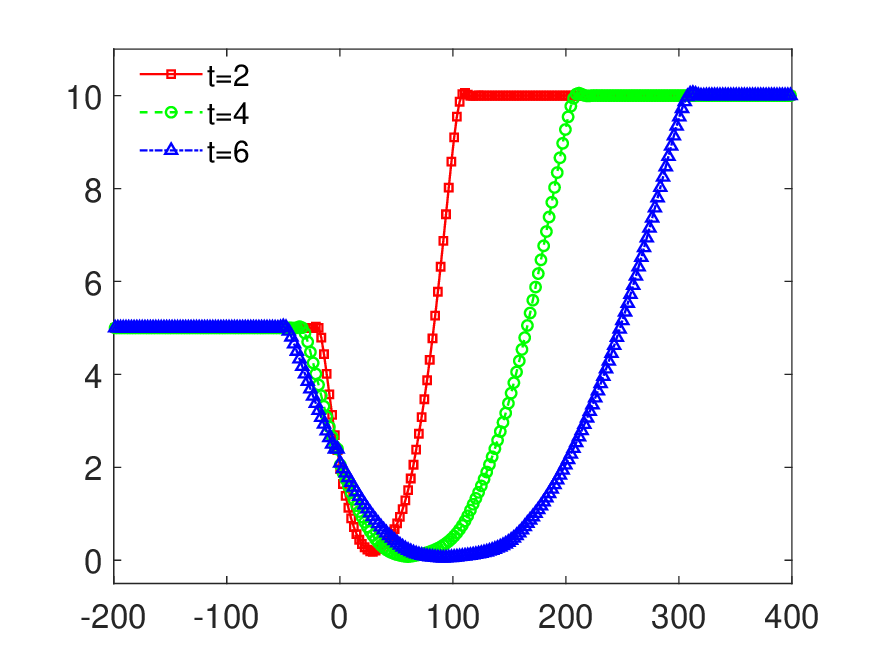}}\hspace*{0.5cm}
	\subfigure[momentum $\xbar{(hu)}$ for Test 2]{\includegraphics[trim=1.0cm 0.05cm 0.5cm 0.05cm,clip,width=5.25cm]{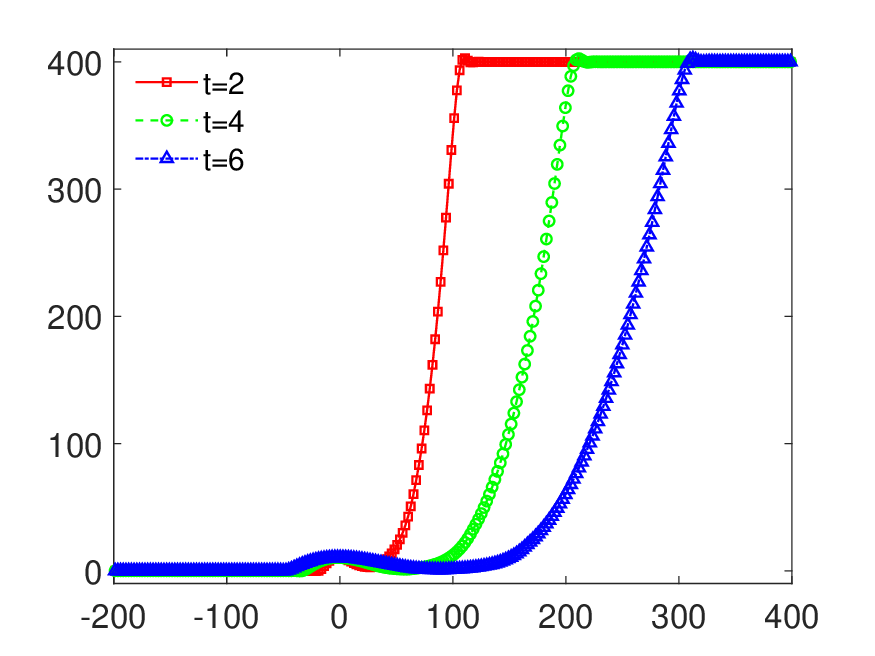}}}
\caption{\sf Example 5: Numerical solutions over a flat bottom.\label{Ex5a}}
\end{figure}

In the fifth example, we consider two test cases, taken from \cite{Xing2010}, to demonstrate the positivity-preserving property of the proposed PP WB PAMPA scheme. The bottom topography is flat ($B(x)=0$) and the initial data are given by 
\begin{equation*}
   \begin{aligned}
&\mbox{Test 1:}&& h(x,0)=\left\{\begin{array}{l}10,~~~~x\in[-300,0],\\0,~~~~ x\in(0,300],\end{array}\right. && hu(x,0)=0,\\
&\mbox{Test 2:}&& h(x,0)=\left\{\begin{array}{l}5,~~~~ x\in[-200,0],\\10,~~~~ x\in(0,400],\end{array}\right. && hu(x,0)=\left\{\begin{array}{l}0,~~~~ x\in[-200,0],\\400,~~~~ x\in(0,400],\end{array}\right.
\end{aligned}
\end{equation*}
We compute the numerical solution using the proposed PP WB PAMPA scheme on a uniform mesh with $250$ cells. The solutions at three different time snapshots are plotted in Figure \ref{Ex5a}. No negative water height is generated during the computation, and a good agreement with those obtained in \cite{Xing2010}.

\subsection*{Example 6---Riemann Problems over a Non-flat Bottom}
\begin{figure}[ht!]
\centerline{\subfigure[surface $\xbar h+\xbar B$ and topography $\xbar B$]{\includegraphics[trim=1.0cm 0.05cm 0.5cm 0.05cm,clip,width=5.25cm]{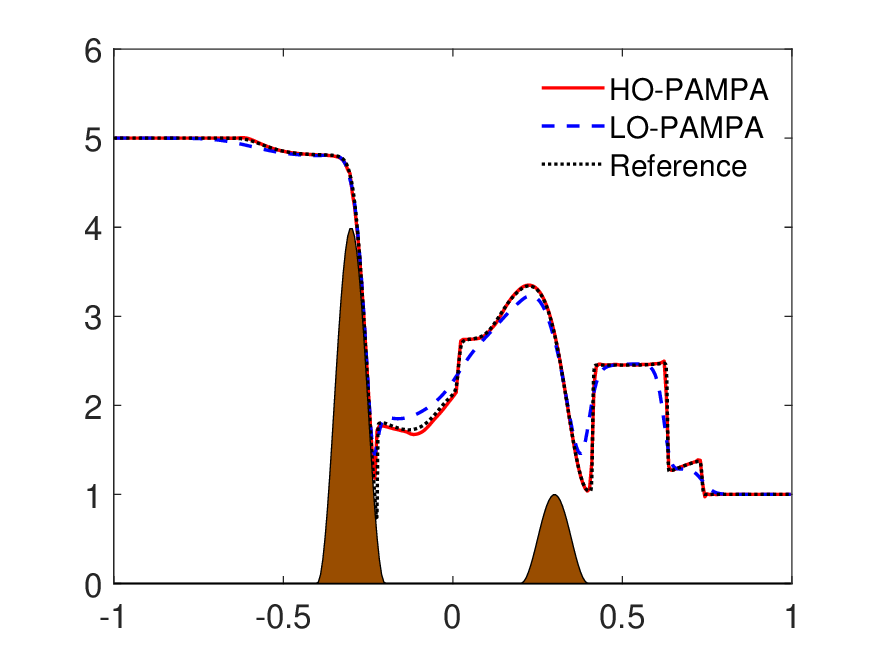}}\hspace*{0.5cm}
	\subfigure[momentum $\xbar{(hu)}$]{\includegraphics[trim=1.0cm 0.05cm 0.5cm 0.05cm,clip,width=5.25cm]{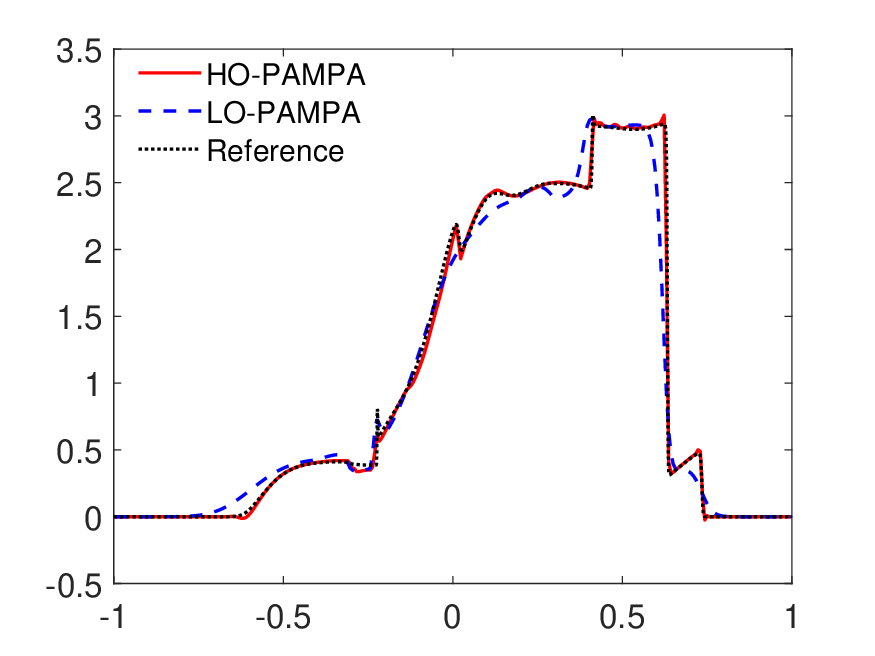}}}
\caption{\sf Example 6: Numerical solutions over a non-flat bottom.\label{Ex6a}}
\end{figure}

In the sixth example, we consider a dam-break over a non-flat bottom. The initial conditions are given by:
\begin{equation*}
  h(x,0)=\left\{\begin{aligned}
  &5-B(x), &&x<0,\\
  &1, &&x\geq0,
  \end{aligned}\right. \quad u(x,0)=0,
\end{equation*}
where the topography function is defined as
\begin{equation*}
  B(x)=\left\{\begin{aligned}
  &2\big(\cos(10\pi(x+0.3))+1\big),&&-\frac{2}{5}\leq x\leq-\frac{1}{5},\\
  &\frac{1}{2}\big(\cos(10\pi(x-0.3))+1\big), &&\frac{1}{5}\leq x\leq\frac{2}{5},\\
  &0,&&\mbox{otherwise}.
  \end{aligned}\right.
\end{equation*}
The computational domain is $[-1,1]$, the gravitational constant is $g=1$, and the numerical solutions are computed on a uniform mesh with 300 cells until  the final time $t=0.3$. We plot the obtained solution in Figure \ref{Ex6a}, together with a reference solution computed by the LO-PAMPA scheme using 3000 cells. We can see that, the results obtained by the HO-PAMPA and LO-PAMPA are consistent with the reference solution and the HO-PAMPA scheme outperforms the LO-PAMPA scheme. Moreover, the PP property of the proposed schemes is also verified, as one can clearly observe that the area near $x=-0.3$ is now almost dry and well resolved.   

\subsection*{Example 7---Parabolic Bowl}
\begin{figure}[ht!]
\centerline{\subfigure[surface at $t=1000$]{\includegraphics[trim=0.85cm 0.05cm 0.8cm 0.05cm,clip,width=4.15cm]{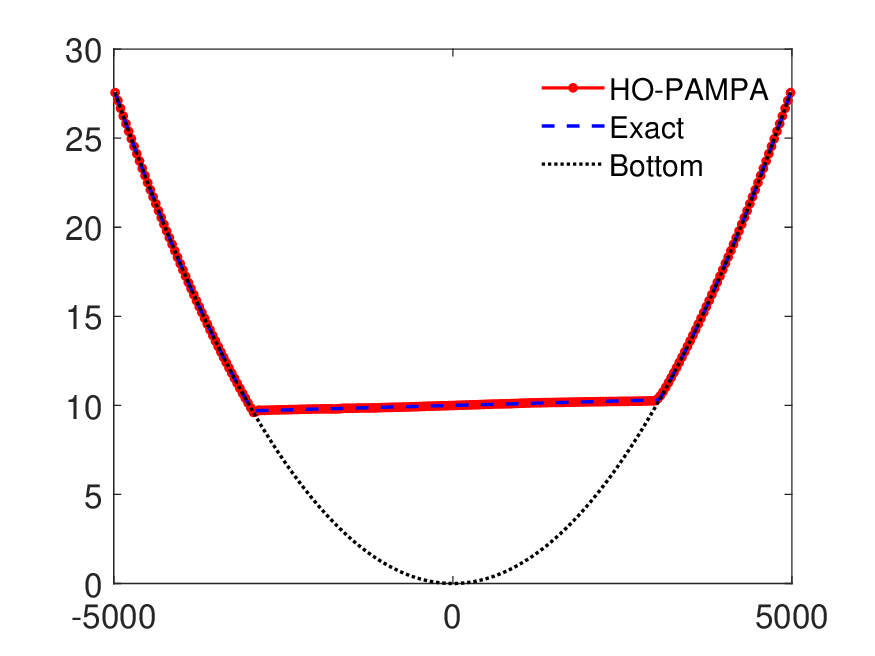}}\hspace*{0.05cm}
	\subfigure[surface at $t=2000$]{\includegraphics[trim=0.85cm 0.05cm 0.8cm 0.05cm,clip,width=4.15cm]{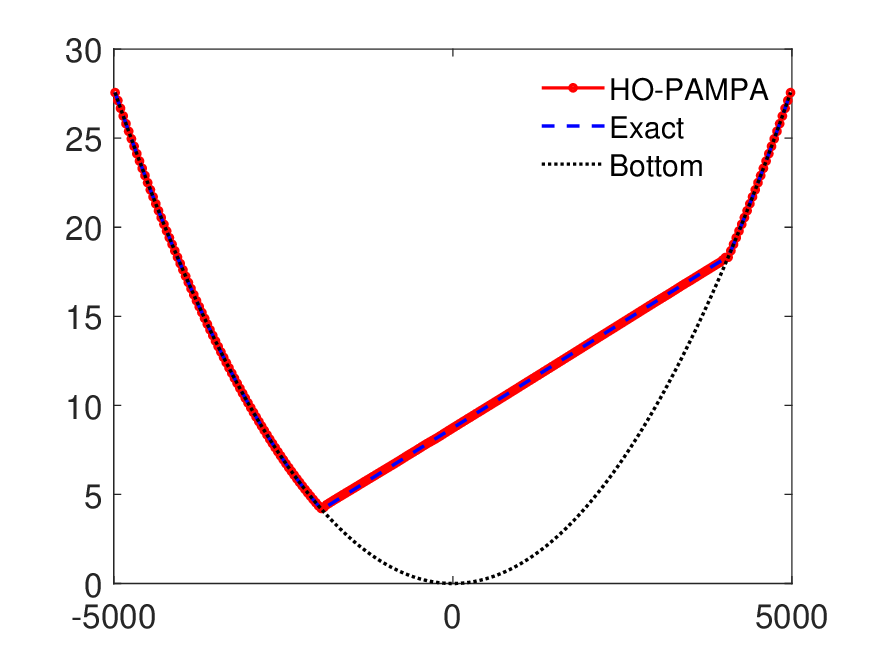}}\hspace*{0.05cm}
	\subfigure[surface at $t=3000$]{\includegraphics[trim=0.85cm 0.05cm 0.8cm 0.05cm,clip,width=4.15cm]{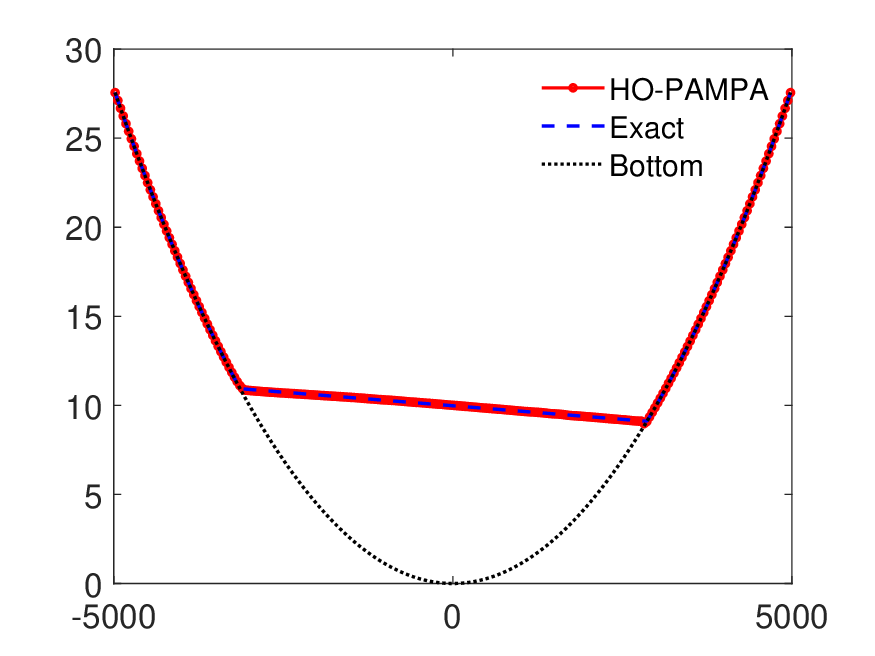}}}
\vskip5pt
\centerline{\subfigure[surface at $t=4000$]{\includegraphics[trim=0.85cm 0.05cm 0.8cm 0.05cm,clip,width=4.15cm]{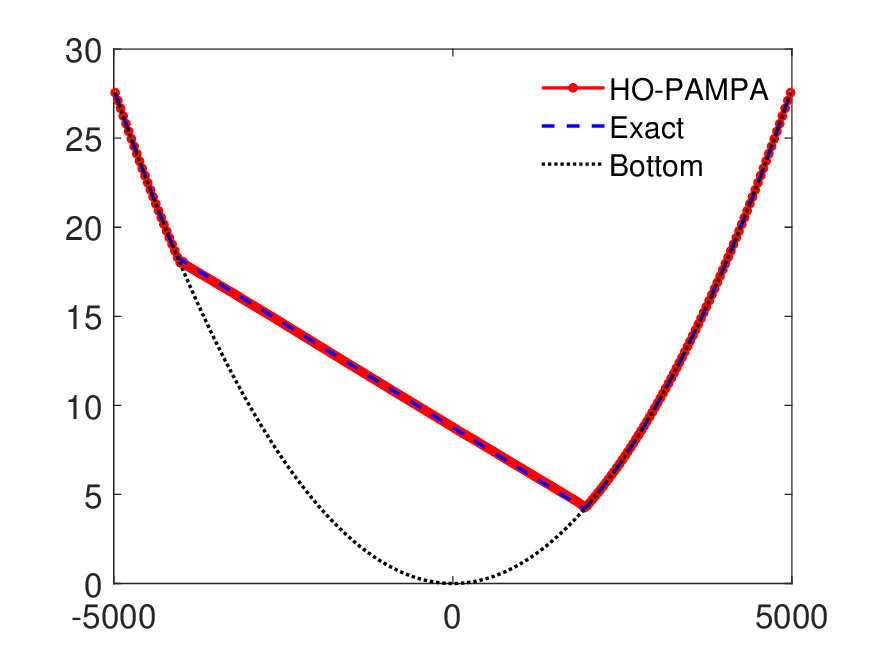}}\hspace*{0.05cm}
	\subfigure[surface at $t=5000$]{\includegraphics[trim=0.85cm 0.05cm 0.8cm 0.05cm,clip,width=4.15cm]{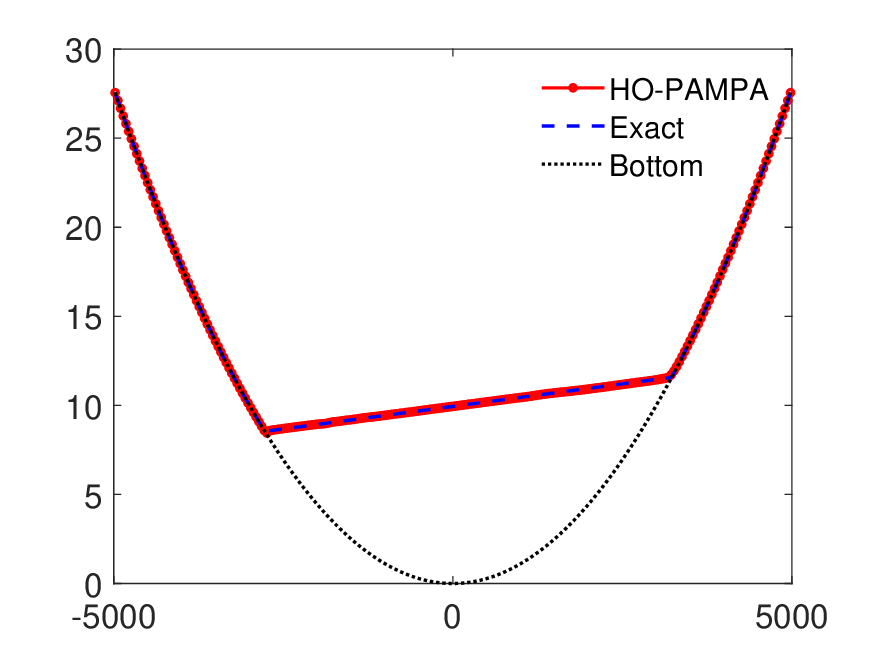}}\hspace*{0.05cm}
	\subfigure[surface at $t=6000$]{\includegraphics[trim=0.85cm 0.05cm 0.8cm 0.05cm,clip,width=4.15cm]{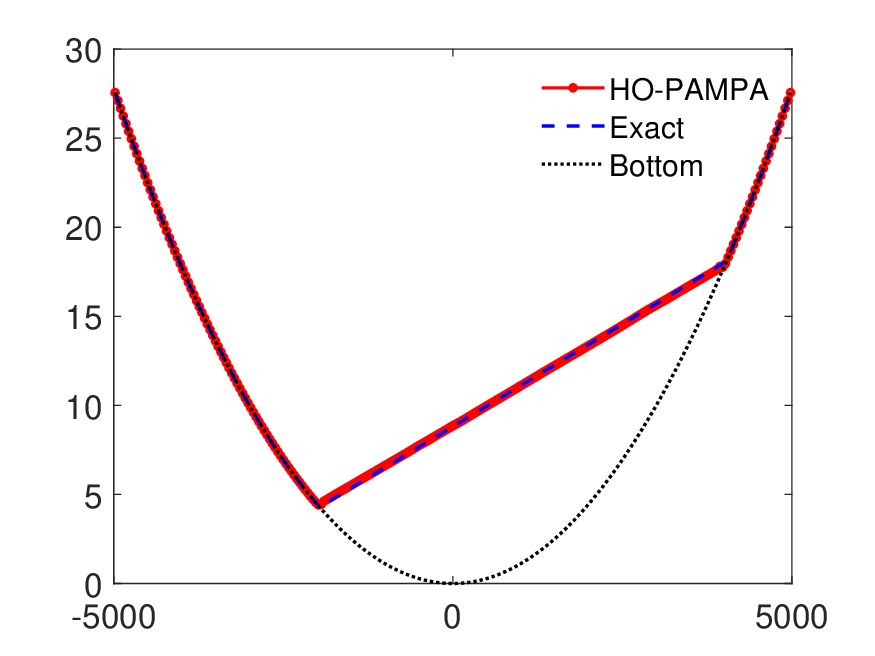}}}
\caption{\sf Example 7: The water surface level in the parabolic bowl problem at different times.\label{Ex7_para}}
\end{figure}

In the seventh example, taken from \cite{Xing2011,Xing2010}, we consider the flow with a parabolic bottom topography, 
\begin{equation*}
  B(x)=h_0(x/a)^2,\quad h_0=10,\quad a=3000,
\end{equation*}
prescribed in the domain $[-5000,5000]$. The analytical water surface for \eqref{1.3} without the friction term ($n=0$) is given by
\begin{equation}\label{Ex7-data1}
  h(x,t)+B(x)=h_0-\frac{b^2}{4g}\cos(2\omega t)-\frac{b^2}{4g}-\frac{bx}{2a}\sqrt{\frac{8h_0}{g}}\cos(\omega t),
\end{equation}
where $\omega=\sqrt{2gh_0/a}$ and $b=5$. The exact location of the wet/dry front takes the form:
\begin{equation}\label{Ex7-data2}
  x_0=-\frac{b\omega a^2}{2gh_0}\cos(\omega t)\pm a.
\end{equation}
The initial condition is then given by \eqref{Ex7-data1}--\eqref{Ex7-data2} and $hu(x,0)=0$. We run the simulations until the final time $t=6000$ using 250 uniform cells. The averaged water surface level ($\xbar h_\jph+\xbar B_\jph$, $\forall j$) at different times are shown in Figure \ref{Ex7_para}, together with the exact solution given by \eqref{Ex7-data1}--\eqref{Ex7-data2} for a comparison. As one can clearly see, the numerical solutions present a nice agreement with the exact ones and also further confirms the PP property of the proposed PP WB PAMPA schemes.

\subsection*{Example 8---Convergence to Moving-water equilibrium with Coriolis forces}
\begin{figure}[ht!]
\centerline{\subfigure[surface $\xbar h+\xbar B$]{\includegraphics[trim=0.01cm 0.05cm 0.8cm 0.05cm,clip,width=4.15cm]{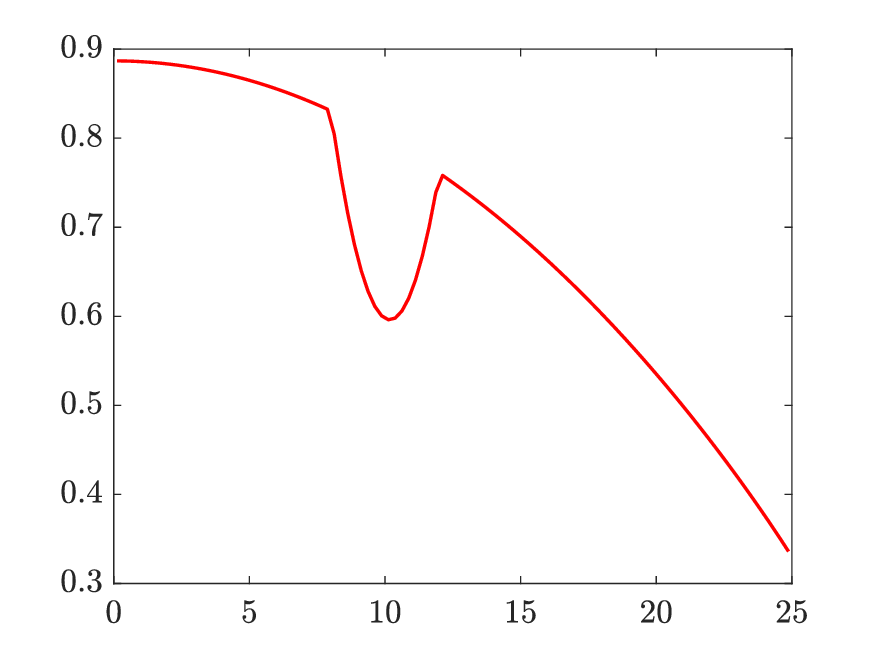}}\hspace*{0.05cm}
	\subfigure[momentum $\xbar {(hu)}$]{\includegraphics[trim=0.01cm 0.05cm 0.8cm 0.05cm,clip,width=4.15cm]{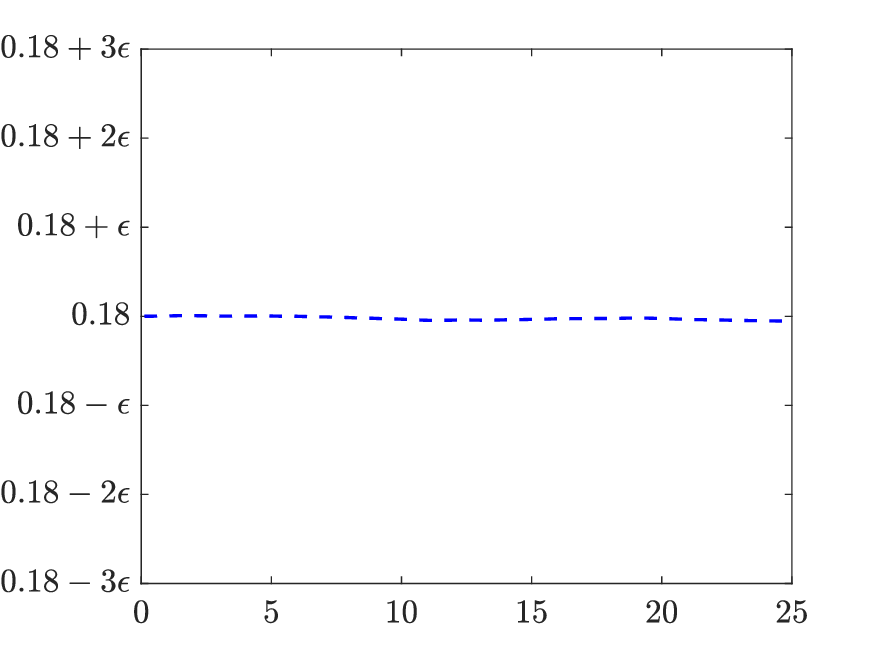}}\hspace*{0.05cm}
	\subfigure[momentum $\xbar {(hv)}$]{\includegraphics[trim=0.01cm 0.05cm 0.8cm 0.05cm,clip,width=4.15cm]{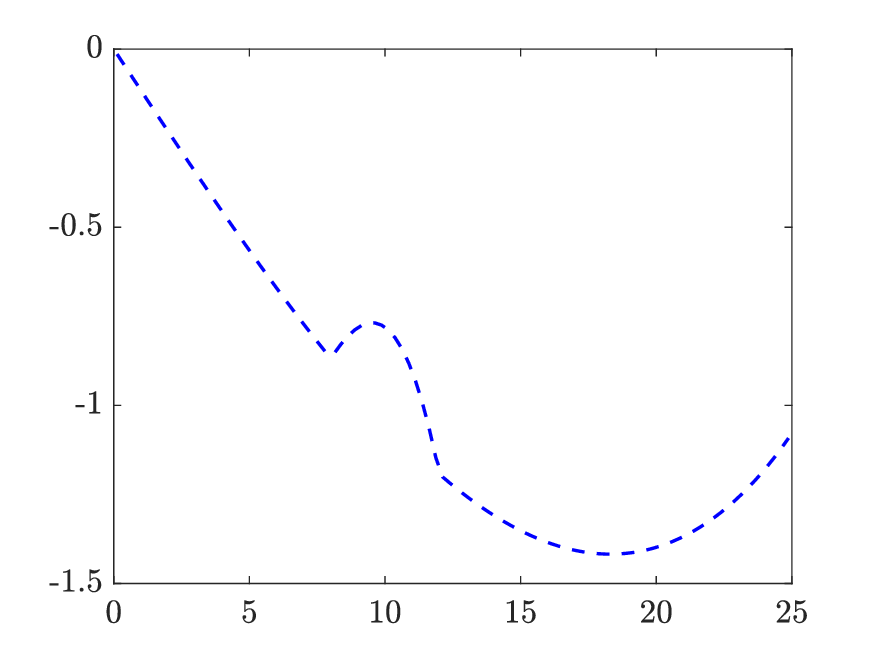}}}
\vskip5pt
\centerline{\subfigure[Global flux $\mbG^{(2)}$]{\includegraphics[trim=0.01cm 0.05cm 0.5cm 0.05cm,clip,width=5.25cm]{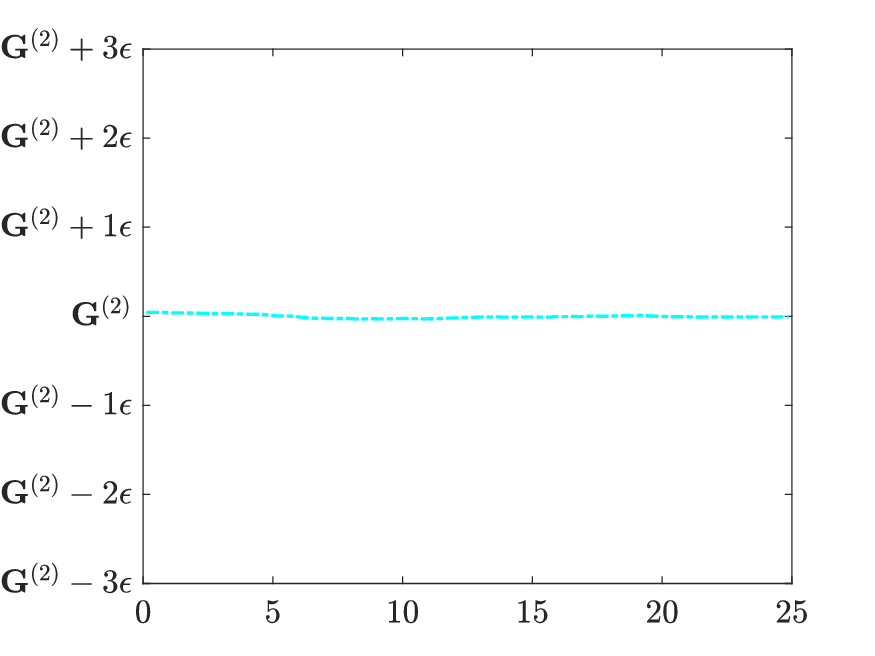}}\hspace*{0.5cm}
	\subfigure[Global flux $\mbG^{(3)}$]{\includegraphics[trim=0.01cm 0.05cm 0.5cm 0.05cm,clip,width=5.25cm]{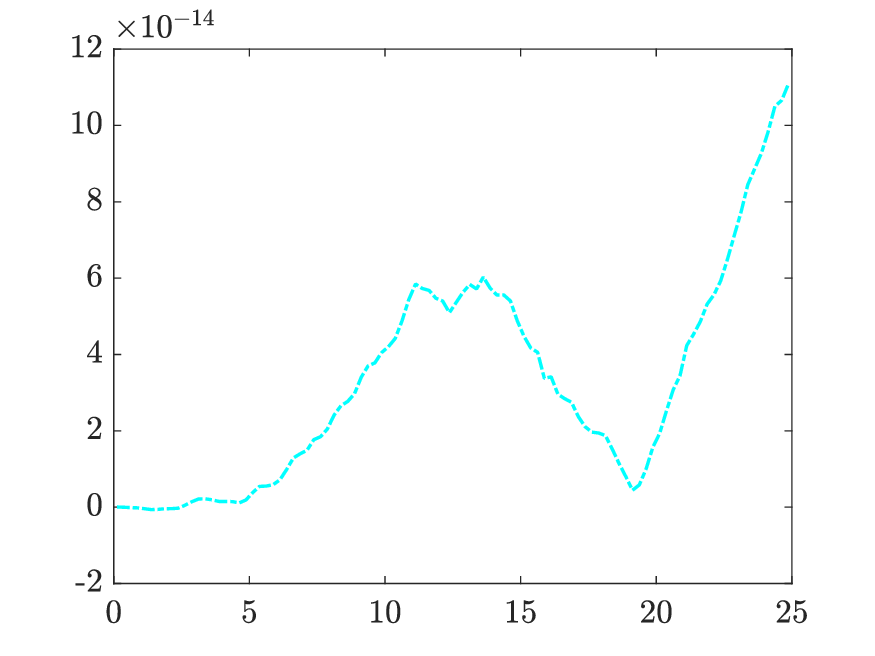}}}
\caption{\sf Example 8: Convergent solutions and discrete values of global fluxes for the flow with Coriolis forces. $\mbG^{(2)}=3.8948194772284$ and $\epsilon=10^{-2}$.\label{Ex8_sol}}
\end{figure}

In the eighth example, taken from \cite{Cao2022}, we consider the following initial data and boundary conditions:
\begin{equation*}
\left\{\begin{array}{l}h(x,0)=0.33,\quad ~hu(x,0)=0,\quad ~~~hv(x,0)=0,\\h(25,t)=0.33,\quad hu(0,t)=0.18,\quad hv(0,t)=0.\end{array}\right.
\end{equation*}
The bottom topography is the same one given by \eqref{Ex3_bottom} and the computational domain is also $[0,25]$. We set the Coriolis force $f=\frac{2\pi}{50}+0.01x$ and run the simulations on a uniform mesh with 100 cells until a final time $t=1000$. The obtained numerical solutions and the discrete values of global flux ($\mbG^{(2)}$ and $\mbG^{(3)}$, where $\mbG^{(1)}=hu$) are reported in Figure \ref{Ex8_sol}. From these results, we clearly see that, the discrete steady states are reached. We can also perform the same test as in Example 4 by adding a small perturbation to the steady-state solution and to examine the propagation captured by the proposed PP WB PAMPA scheme. Whereas, for saving space, we omit the results here.

%\subsection*{Example 9---Super-convergent to a Moving Equilibrium}
%In the ninth example, we consider the manufactured state used in \cite{Mantri2024,Cao2022,Desveaux2021} corresponding to a moving steady state of \eqref{1.3a}:
%\begin{equation*}
%  h(x,0)=\mathrm e^{2x},\quad hu(x,0)=1, \quad hv(x,0)=-fx\mathrm e^{2x},
%\end{equation*}
%and a non-flat bottom topography
%\begin{equation*}
%  B(x)=-\frac{1}{2}f^2x^2-\mathrm e^{2x}-\frac{1}{2}\mathrm e^{-4x},
%\end{equation*}
%prescribed in a domain $[0,1]$. The Coriolis force is $f=1$ and the gravitational constant is $g=1$.
\subsection*{Example 9---Geostrophic Equilibrium over a Flat Bottom}
\begin{figure}[ht!]
\centerline{\subfigure[cell averages of $h$]{\includegraphics[trim=1.0cm 0.05cm 0.5cm 0.05cm,clip,width=5.25cm]{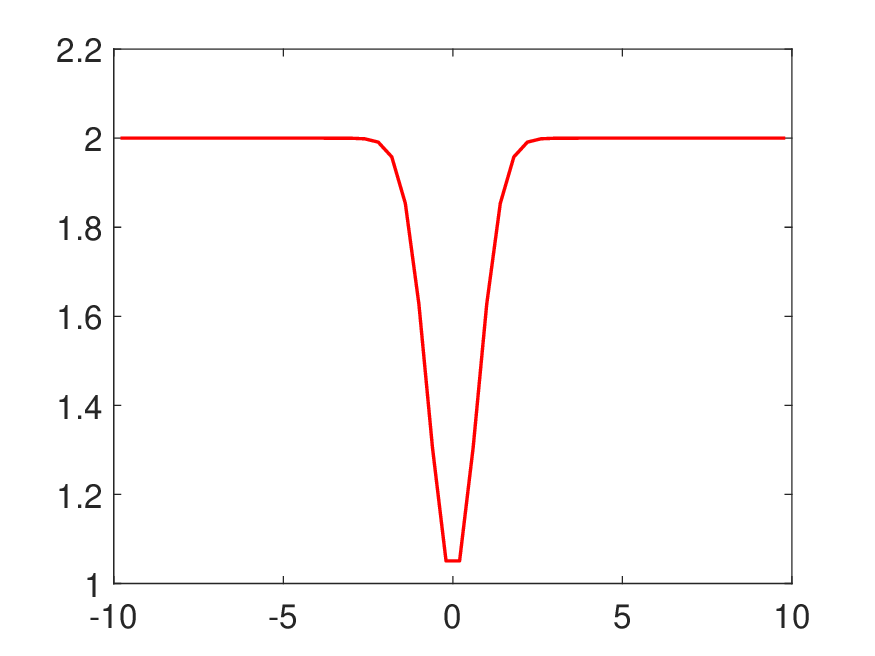}}\hspace*{0.5cm}
	\subfigure[point values of $h$]{\includegraphics[trim=1.0cm 0.05cm 0.5cm 0.05cm,clip,width=5.25cm]{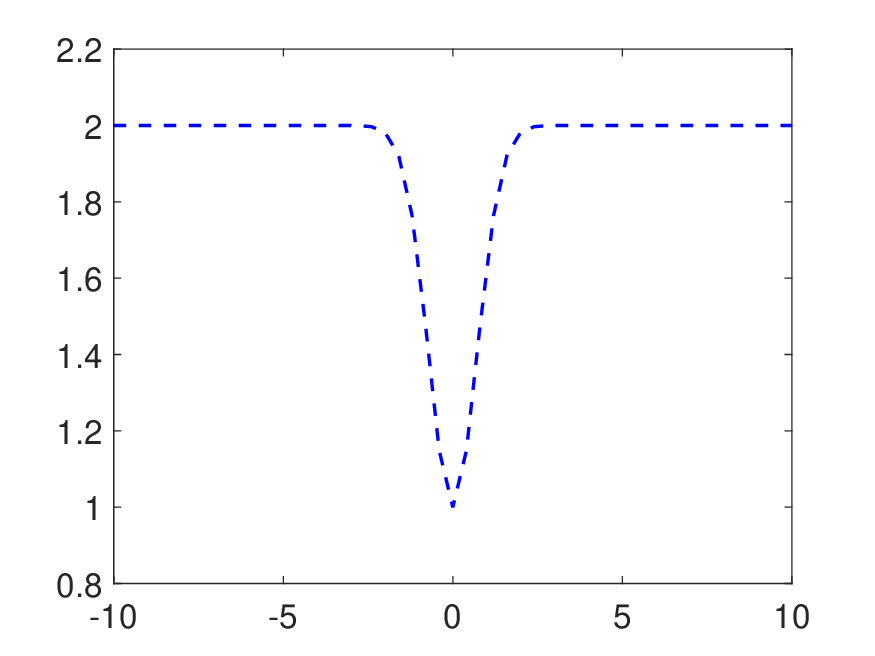}}}
\caption{\sf Example 9: Prepared discrete steady-state water height.\label{Ex9a}}
\end{figure}

In the ninth example, we consider the following initial conditions, taken from \cite{Chertock2018a}:
\begin{equation*}
  u(x,0)=0,\quad v(x,0)=\frac{2g}{f}x{\mathrm e}^{-x^2},
\end{equation*}
where $g=1$ and $f=10$. The discrete values for $h(x)^{\rm eq}$ are calculated analytically by solving the equations in \eqref{eqs:nonlinear} with $\mbG^{(2)}=2$ in the domain $[-10,10]$ covered by $50$ cells. In Figure \ref{Ex9a}, we plot the prepared discrete steady-state water height. We then use these initial setting to run the simulation until the final time $t=100$ and report the errors between the numerical solutions and the initial data in Table \ref{tab9}. Again, the WB property of the proposed scheme is verified.

\begin{table}[!ht]
\caption{\sf Example 9: Errors in $h$, $hu$, and $hv$ computed by high-order scheme.\label{tab9}}
\begin{center}
\begin{tabular}{ c| c c c}\hline 
\multicolumn{1}{c|}{\multirow{1}{*}{Variables}} &$L^1$-error &$L^2$-error&$L^\infty$-error\\ \hline
$h$ &$0$ & $0$  & $0$ \\ \hline
$hu$ &$6.63\times10^{-16}$ & $2.49\times10^{-16}$  & $2.29\times10^{-16}$ \\ \hline
$hv$ &$5.66\times10^{-15}$ & $1.97\times10^{-15}$  & $1.09\times10^{-15}$ \\ \hline
\end{tabular}
\end{center}
\end{table}

\subsection*{Example 10---Geostrophic Equilibrium over a Non-flat Bottom}
In the tenth example, we consider a non-flat bottom topography,
\begin{equation*}
  B(x)=\left\{
  \begin{aligned}
  &0.25(\cos(10\pi(x-0.8))+1),&&\mbox{if}~0.7\leq x\leq 0.9,\\
  &0,&&\mbox{otherwise}.
  \end{aligned}\right.
\end{equation*}
and the initial conditions again satisfies a geostrophic equilibrium given by
\begin{equation*}
  \mbG^{(2)}=2, \quad u(x,0)=0,\quad v(x,0)=0.05\sin(2\pi x),
\end{equation*}
over the domain $[0,1]$. We use $g=1$ and $f=10$ here. The initial values of water depth is prepared by solving the nonlinear equations \eqref{eqs:nonlinear}. We first run the simulation using 20 cells until a final time $t=20$ and report the errors in Table \ref{tab10}. As one can see again, the proposed PP WB PAMPA scheme is capable of exactly preserving the geostrophic equilibrium within machine accuracy even with a non-flat bottom topography.   
\begin{table}[!ht]
\caption{\sf Example 10: Errors in $h$, $hu$, and $hv$ computed by high-order scheme.\label{tab10}}
\begin{center}
\begin{tabular}{ c| c c c}\hline 
\multicolumn{1}{c|}{\multirow{1}{*}{Variables}} &$L^1$-error &$L^2$-error&$L^\infty$-error\\ \hline
$h$ &$6.48\times 10^{-15}$ & $8.04\times 10^{-15}$  & $1.02\times 10^{-14}$ \\ \hline
$hu$ &$1.89\times10^{-16}$ & $2.70\times10^{-16}$  & $7.39\times10^{-16}$ \\ \hline
$hv$ &$4.62\times10^{-15}$ & $6.53\times10^{-15}$  & $1.58\times10^{-14}$ \\ \hline
\end{tabular}
\end{center}
\end{table}

\begin{figure}[ht!]
\centerline{\subfigure[$\xbar h(x,t)-\xbar h(x,0)$]{\includegraphics[trim=1.0cm 0.05cm 0.5cm 0.05cm,clip,width=5.25cm]{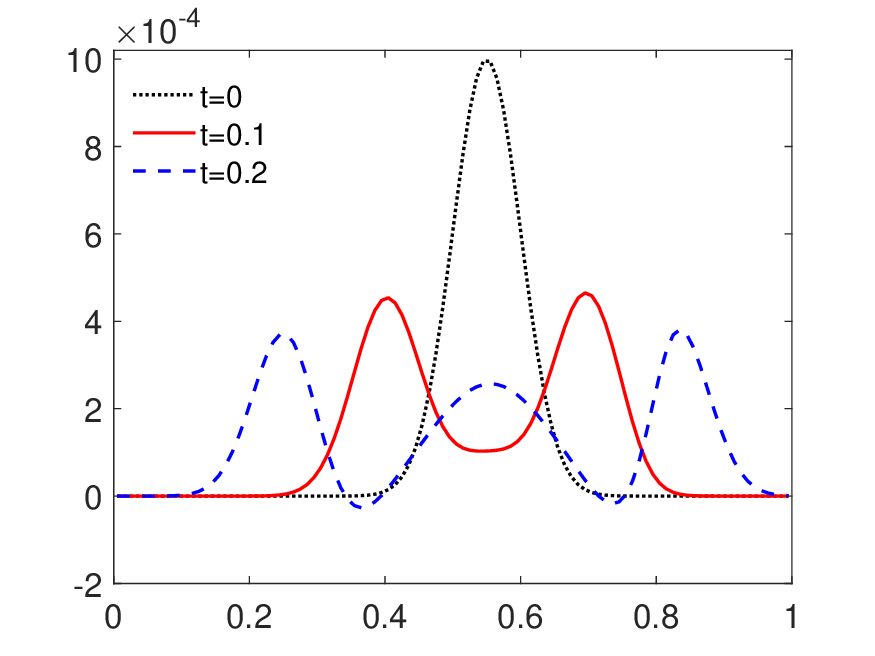}}\hspace*{0.5cm}
	\subfigure[$\xbar {(hu)}(x,t)-\xbar {(hu)}(x,0)$]{\includegraphics[trim=1.0cm 0.05cm 0.5cm 0.05cm,clip,width=5.25cm]{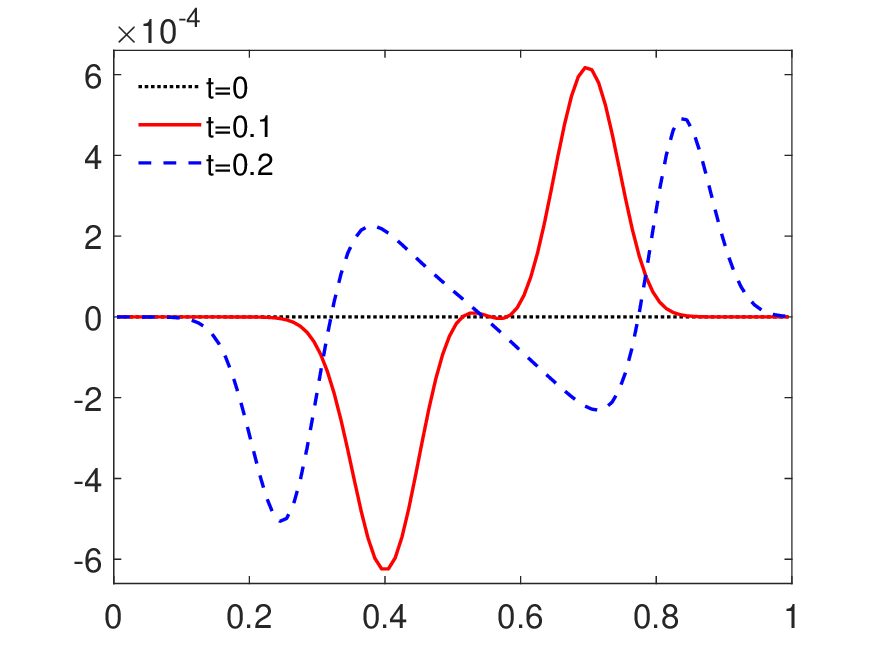}}}
\caption{\sf Example 10: Difference between the numerical solution and the background steady states.\label{Ex10diff}}
\end{figure}

To further examine the WB property of the scheme on handling with the geostrophic equilibrium, we study the propagation of a small perturbation. To this end, we add a small perturbation, $10^{-3}\mathrm{e}^{-200(x-0.55)^2}$, to the discrete steady-state water depth and compute the difference between the numerical solution and the background steady states. The results at times $t=0$, $0.1$, and $0.2$ are shown in Figure \ref{Ex10diff}. We can clearly see that the proposed PP WB PAMPA scheme well captures a small perturbation.
\subsection*{Example 11---Rossby Adjustment in an Open Domain}
\begin{figure}[ht!]
\centerline{\subfigure[$\xbar h $ at $t=0$]{\includegraphics[trim=0.85cm 0.05cm 0.8cm 0.05cm,clip,width=4.15cm]{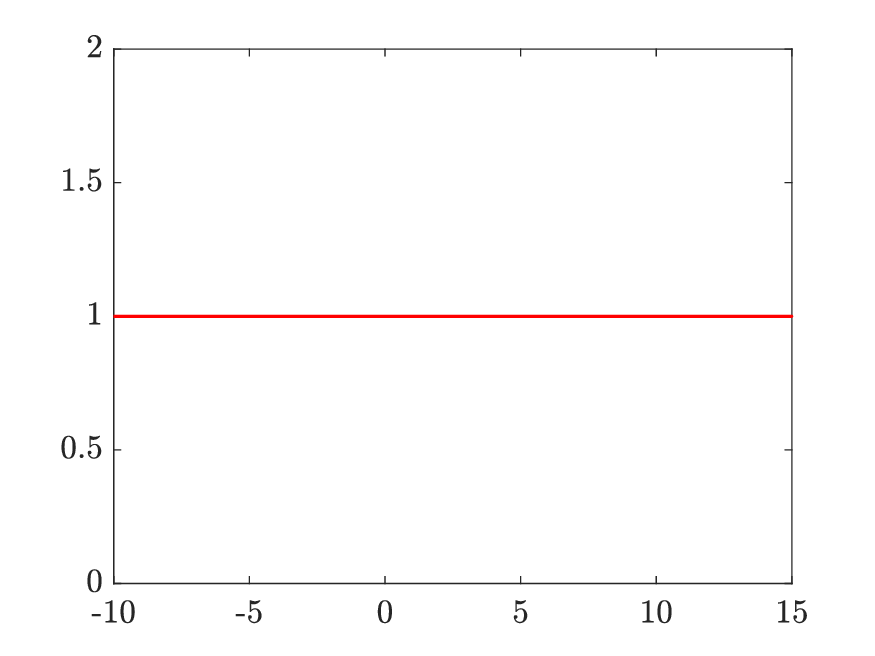}}\hspace*{0.05cm}
	\subfigure[$\xbar h $ at $t=0.2T_f$]{\includegraphics[trim=0.85cm 0.05cm 0.8cm 0.05cm,clip,width=4.15cm]{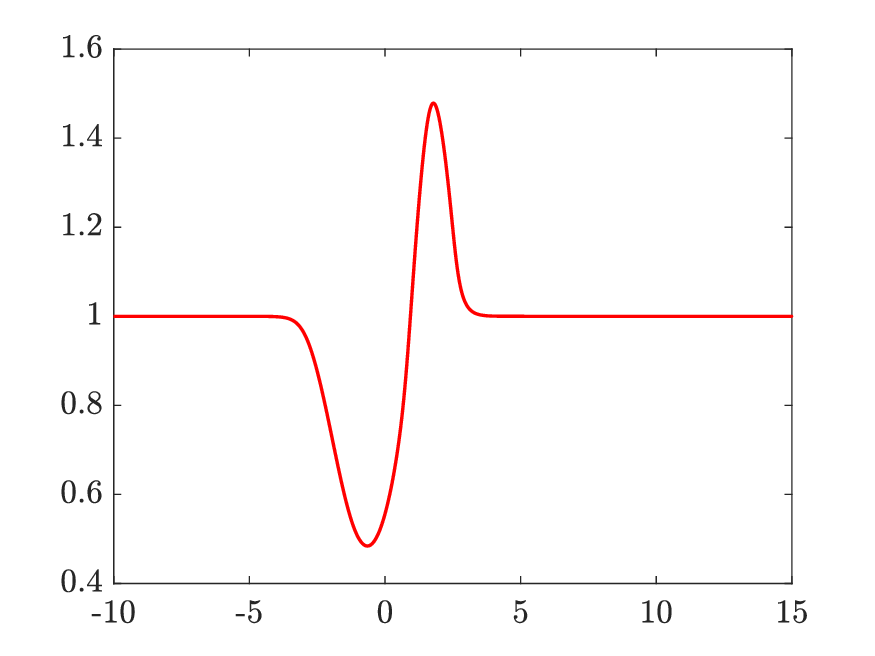}}\hspace*{0.05cm}
	\subfigure[$\xbar h $ at $t=0.4T_f$]{\includegraphics[trim=0.85cm 0.05cm 0.8cm 0.05cm,clip,width=4.15cm]{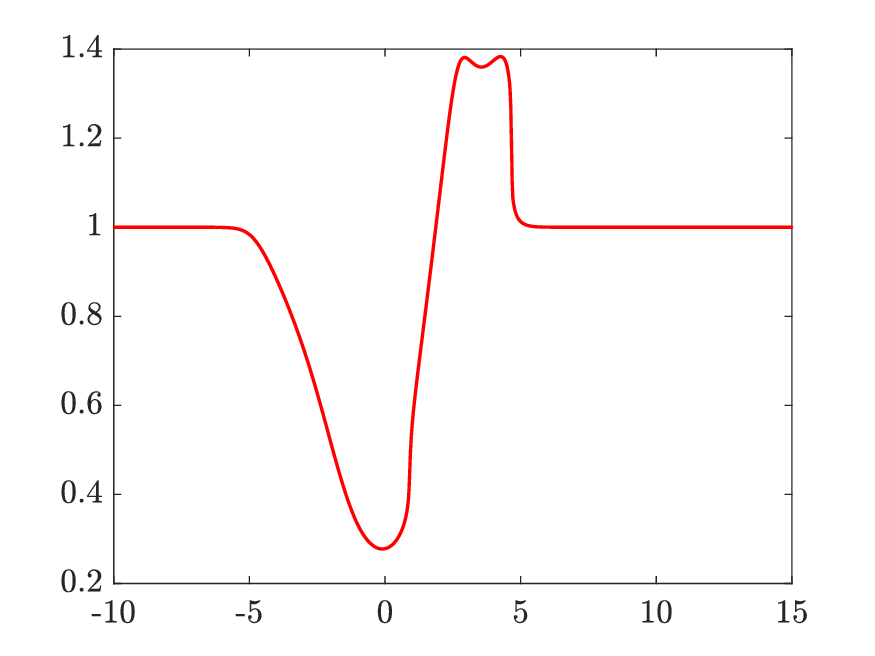}}}
\vskip5pt
\centerline{\subfigure[$\xbar h $ at $t=0.6T_f$]{\includegraphics[trim=0.85cm 0.05cm 0.8cm 0.05cm,clip,width=4.15cm]{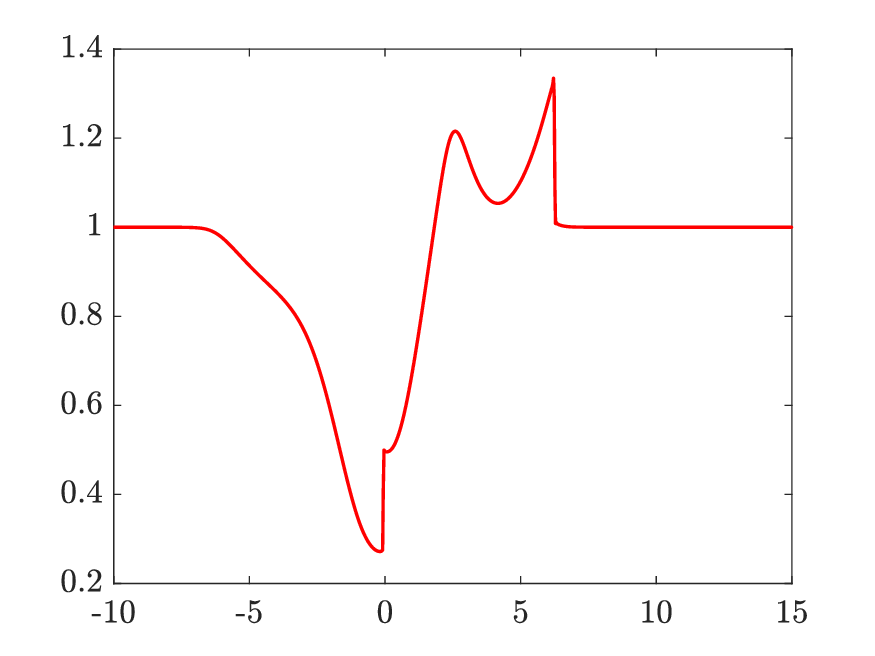}}\hspace*{0.05cm}
	\subfigure[$\xbar h $ at $t=0.8T_f$]{\includegraphics[trim=0.85cm 0.05cm 0.8cm 0.05cm,clip,width=4.15cm]{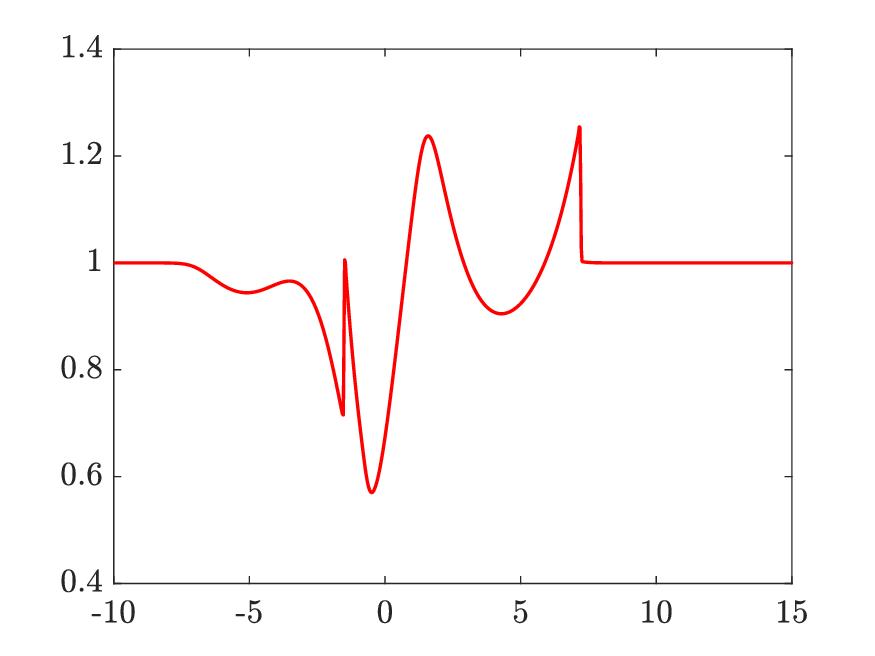}}\hspace*{0.05cm}
	\subfigure[$\xbar h $ at $t=T_f$]{\includegraphics[trim=0.85cm 0.05cm 0.8cm 0.05cm,clip,width=4.15cm]{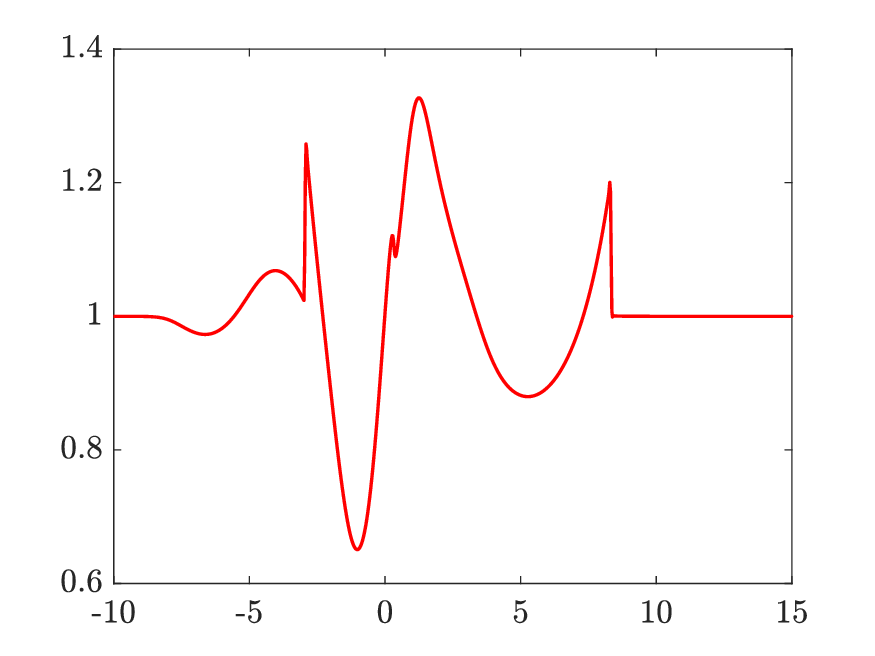}}}
\caption{\sf Example 11: time evolution of $h$ at different times.\label{Ex12}}
\end{figure}

In the last example, we numerically investigate the Rossby adjustment problem with the constant Coriolis parameter $f=1$ and $g=1$, which was studied previously in \cite{Chertock2018a,Castro2008}. We consider the following initial conditions, which correspond to a jet over a flat $h$:
\begin{equation*}
  h(x,0)=1,\quad u(x,0)=0,\quad v(x,0)=2N_2(x),
\end{equation*}
where
\begin{equation*}
  N_2(x)=\frac{(1+\tanh(2x+2))(1-\tanh(2x-2))}{(1+\tanh(2))^2}
\end{equation*}
and the bottom is flat $B(x)=0$. The simulation is run on a uniform mesh with 20000 cells. Time evolution of the water depth along an inertial period $T_f=\frac{2\pi}{f}$ is shown in Figure \ref{Ex12}. As one can see, fast inertia gravity waves are emitted as a consequence of the initial momentum disturbance, and two shocks arise in the wave front.
%We also show long time evolution of the values of $fv$ and $gh_x$ at the cell centers in Figure .  In this example, the solution is expected to converge to a geostrophic steady state, at which these quantities are supposed to be equal. However, as it was pointed out in \cite{Chertock2018a, Kurganov2020,Bouchut2004}, some wave modes have almost zero group-velocity and stay for long time in the core of the jet. Therefore, the convergence to the geostrophic equilibrium is quite slow. 

\section{Conclusion}
In this work, we have developed a new formulation of the Point-Average-Moment PolynomiAl-interpreted (PAMPA) method via the flux globalization to construct positivity-preserving and fully well-balanced numerical schemes for one-dimensional shallow water models. The discrete fully well-balanced condition is achieved by incorporating the source terms into the global flux, resulting in a quasi-conservative system that enables a clear characterization of the discrete steady states and provable superconvergence toward general analytical equilibria. In practice, the required integrals are evaluated using direct applications of Gauss--Lobatto quadrature. To ensure positivity of the water depth, capture wet–dry fronts, and suppress spurious oscillations near strong shocks, first-order schemes based on the hydrostatic reconstruction techniques are introduced. Furthermore, a monolithic convex limiting strategy is employed to determine the blending coefficients for the high- and low-order fluxes and residuals. The proposed schemes are proven to exactly preserve still-water steady states, approximately preserve moving-water and geostrophic equilibria with superconvergent accuracy, and maintain the positivity of the water depth.

Extensive numerical tests for the shallow water equations, with/without Manning friction and Coriolis effect, validate all the theoretical properties on various equilibria, including still-water, moving-water, and geostrophic steady states. Applications to complex problems involving sharp fronts and different wave interactions further demonstrate the accuracy and robustness of the proposed method.

Compared with existing flux globalization based well-balanced methods for one-dimensional shallow water models, the proposed flux globalization based PP WB PAMPA scheme offers several advantages: i) The PP WB PAMPA scheme, using global flux quadrature, is exactly well-balanced for still-water steady states, as proven in \cref{prop2}. The sub-cell LobattoIII method has been used for the source term integral in the half cell. We have shown that for the shallow water equations, it automatically exactly preserve the still-water states. In contrast, the flux globalization based DGSEM methods in \cite{Mantri2024,Ciallella2023}, where the LobbatoIIIA has been used to compute the integral on the half cell, require an additional source term reformulation following the approach of \cite{Xing2005}, while the flux globalization based central-upwind scheme of \cite{Cheng2019} does not achieve this property. ii) The PP WB PAMPA scheme is positivity-preserving and capable of handling wet–dry fronts through blending with a positivity- and still-water preserving first-order scheme. This capability is absent in other flux globalization methods, such as those in \cite{Mantri2024,Cheng2019,Ciallella2023}. iii) A simple oscillation-eliminating parameter is introduced to adaptively adjust the blending coefficients near shock waves, enabling the PP WB PAMPA scheme to robustly capture genuinely discontinuous solutions. This improvement is not realized in the DGSEM framework of \cite{Mantri2024}. iv) The proposed formulation for preserving discrete moving-water steady states is straightforward, fully local,  and computationally efficient, whereas the flux globalization based central-upwind schemes of \cite{Cheng2019,Cao2022,Chertock2022} require solving additional nonlinear equations to recover well-balanced reconstructed point values (a generalized hydrostatic reconstruction step). v) We note that, despite the use of integrals to define the global fluxes
(see e.g. equation \eqref{2.5}), \cref{compactness} shows that  the final numerical method   is actually compact. This property enables their blending with the first-order schemes. For the evolution of point values, the computation in \eqref{2.6c} similarly confirms this compactness. A current limitation of the PP WB PAMPA scheme is its inability to handle discontinuous bottom topography. This could be done with some     estimation of the jump in $\mbR$; see, e.g., \cite{Cao2022,Ciallella2023,Kazolea2025}.  

The extensions to multidimensional problems and to more complex systems will be investigated in the near future.

\appendix 

\section{Consistency  of sub-cell LobattoIII collocation method}\label{app:sc-lobattoIII}

We prove here the statement of proposition \ref{pp:sc-lobattoIII}.
We use explicitly the fact that the main unknowns are approximations of 
$$
\mbu_{j}\approx  \mbu(x_{j}),\;
\mbu_{j+1/2}\approx  \mbu(x_{j}+\dx/2),\;
\mbu_{j+1}\approx  \mbu(x_{j}+\dx),
$$
and we formally replace the above with sampled values of a smooth exact solution in the discrete equations to evaluate the consistency error. Without loss of generality, we consider the problem
$$
\mbf'= \mbS(\mbu(\mbf),x)
$$
which is the steady ODE of the balance law, with classical assumptions that $\mbf(\mbu)$ is a one-to-one smooth mapping, and $\mbS$ has bounded derivatives.

Using a truncated Taylor series expansion we can readily show that
 $$
 \begin{aligned}
 &\mbf_{j+1}-\mbf_{j}
 -\dfrac{\dx}{6}\mbS(\mbu_j,x_j)
 -2\dfrac{\dx}{3}\mbS(\mbu_\jph,x_\jph)
 -\dfrac{\dx}{6}\mbS(\mbu_{j+1},x_{j+1})
 \\&= \mbf_{j+1}-\mbf_{j} 
 -\dfrac{\dx}{6}\mbf'_j
 -2\dfrac{\dx}{3}\mbf'_\jph 
 -\dfrac{\dx}{6}\mbf'_{j+1} \\
 & =\mbf_{j+1}-\mbf_{j}-\mbf'_j\dx  -\mbf''_j\dfrac{\dx^2}{2}
 -\mbf'''_j\dfrac{\dx^3}{6}
 -\mbf''''_j\dfrac{\dx^4}{24}   
 -\mbf'''''_j\dfrac{5\dx^5}{576}=\mathcal{O}(\dx^5).
 \end{aligned}
 $$
 This shows that 
$$
\sum_j\{ \mbf_{j+1}-\mbf_{j}
 -\dfrac{\dx}{6}\mbS(\mbu_j,x_j)
 -2\dfrac{\dx}{3}\mbS(\mbu_\jph,x_\jph)
 -\dfrac{\dx}{6}\mbS(\mbu_{j+1},x_{j+1})\} = \mathcal{O}(\dx^4),
$$
provided we have at least $\mbu_\jph-\mbu(x_\jph) =\mathcal{O}(\dx^4)$. To show this we first   compute (using again 
a truncated Taylor series expansion) 
$$
\mbu_{j+1/4}=\dfrac{3}{8}\mbu_{j}+\dfrac{3}{4}\mbu_{j+1/2}
- \dfrac{1}{8}\mbu_{j+1} = \mbu(x_j + \dx/4) -\dfrac{\dx^3}{128} \mbu''' + \mathcal{O}(\dx^4).
$$

This allows to show that, as long as $\mbS$ has bounded derivatives, then 
$$
\begin{aligned}
&\mbf_{\jph} -\mbf_j   
 -\dfrac{\dx}{12}\mbS(\mbu_j,x_j)
 -\dfrac{\dx}{3}\mbS(\mbu_{j+1/4},x_{j+1/4})
 -\dfrac{\dx}{12}\mbS(\mbu_\jph,x_\jph)\\
 &=
 \mbf_{\jph} -\mbf_j   
 -\dfrac{\dx}{12}\mbf'_j 
 -\dfrac{\dx}{3}(\mbf'_{j+1/4} -  \mbS_\mbu(\mbu_{j+1/4},x_{j+1/4})
 \dfrac{\dx^3}{128} \mbu''')
 -\dfrac{\dx}{12}\mbf'_\jph  \\
 &=
 \mbf_{\jph} -\mbf_j   
 -\dfrac{\dx}{2}\mbf'_j 
 -\dfrac{\dx^2}{8}\mbf''_j 
 -\dfrac{\dx^3}{48}\mbf'''_j 
 -\dfrac{5 \dx^4}{1728}\mbf''''_j\\ 
 &\quad  +  \mbS_\mbu(\mbu_{j+1/4},x_{j+1/4})
 \dfrac{\dx^4}{384} \mbu'''
 +\mathcal{O}(\dx^5)
 \\&= \big(-\dfrac{5 }{1728}\mbf''''_j 
   +  \mbS_\mbu(\mbu_{j+1/4},x_{j+1/4})
 \dfrac{\mbu'''}{384} \big)\dx^4
 +\mathcal{O}(\dx^5),
 \end{aligned}
$$
which achieves the proof.

\bibliographystyle{siamplain}
\bibliography{refer_list}

\end{document}